\documentstyle[amscd,amssymb,verbatim,diagrams,12pt]{amsart}
\renewcommand{\span}{{\operatorname{span}}}
\newcommand{\univ}{{\operatorname{univ}}}
\newcommand{\tr}{{\operatorname{tr}}}
\newcommand{\he}{{\operatorname{he}}}
\newcommand{\Si}{\Sigma}
\newcommand{\bC}{{\bf C}}
\newcommand{\bn}{{\bf n}}
\newcommand{\HU}{{\cal HU}}

\renewcommand{\Bar}{\operatorname{Bar}}
\newcommand{\Art}{\operatorname{Art}}
\newcommand{\Sets}{\operatorname{Sets}}
\newcommand{\wgt}{\operatorname{wt}}

\newcommand{\vareps}{\varepsilon}

\newcommand{\del}{\partial}

\newcommand{\fG}{{\frak G}}

\renewcommand{\mod}{\operatorname{mod}}

\newcommand{\gr}{\operatorname{gr}}
\newcommand{\cusp}{\operatorname{cusp}}

\renewcommand{\int}{\operatorname{int}}
\newcommand{\Tor}{\operatorname{Tor}}
\newcommand{\und}{\underline}

\newcommand{\OO}{{\cal O}}
\newcommand{\RR}{{\cal R}}

\let\liminv\varprojlim

\newcommand{\NN}{{\cal N}}
\newcommand{\KK}{{\cal K}}
\newcommand{\bL}{{\bf L}}
\newcommand{\be}{{\bf e}}
\newcommand{\bu}{{\bf u}}

\newcommand{\BB}{{\cal B}}

\newcommand{\SL}{\operatorname{SL}}
\newcommand{\unit}{{\bf 1}}
\newcommand{\G}{{\Bbb G}}

\newcommand{\mg}{{\frak m}}
\newcommand{\hra}{\hookrightarrow}
\newcommand{\lan}{\langle}
\newcommand{\ran}{\rangle}

\newcommand{\CC}{{\cal C}}
\newcommand{\UU}{{\cal U}}

\newcommand{\Spec}{\operatorname{Spec}}

\newcommand{\Proj}{\operatorname{Proj}}
\renewcommand{\P}{{\Bbb P}}

\newcommand{\si}{\sigma}

\newcommand{\Pic}{\operatorname{Pic}}

\newcommand{\ga}{\gamma}
\newcommand{\de}{\delta}
\newcommand{\eps}{\epsilon}

\renewcommand{\ker}{\operatorname{ker}}
\newcommand{\im}{\operatorname{im}}

\newcommand{\A}{{\Bbb A}}

\numberwithin{equation}{subsection}

\newtheorem{thm}{Theorem}[subsection]
\newtheorem{prop}[thm]{Proposition}
\newtheorem{lem}[thm]{Lemma}
\newtheorem{cor}[thm]{Corollary}
{  \theoremstyle{definition}
\newtheorem{defi}[thm]{Definition}
\newtheorem{ex}[thm]{Example}

\newtheorem{rem}[thm]{Remark}
\newtheorem{rems}[thm]{Remarks}
}

\newcommand{\Pf}{\noindent {\it Proof}}
\newcommand{\id}{\operatorname{id}}

\newcommand{\ov}{\overline}

\newcommand{\can}{\operatorname{can}}

\renewcommand{\AA}{{\cal A}}

\newcommand{\FF}{{\cal F}}

\newcommand{\MM}{{\cal M}}
\newcommand{\TT}{{\cal T}}
\newcommand{\XX}{{\cal X}}

\newcommand{\SS}{{\cal S}}

\newcommand{\Om}{\Omega}

\newcommand{\MA}{{\cal MA}}

\newcommand{\Hom}{\operatorname{Hom}}

\newcommand{\Ext}{\operatorname{Ext}}
\newcommand{\End}{\operatorname{End}}
\newcommand{\Res}{\operatorname{Res}}

\renewcommand{\a}{\alpha}
\renewcommand{\b}{\beta}
\newcommand{\om}{\omega}

\newcommand{\la}{\lambda}

\newcommand{\R}{{\Bbb R}}
\newcommand{\Z}{{\Bbb Z}}
\newcommand{\Q}{{\Bbb Q}}
\newcommand{\La}{\Lambda}
\newcommand{\Ga}{\Gamma}

\newcommand{\wt}{\widetilde}
\newcommand{\ot}{\otimes}

\newcommand{\sub}{\subset}
\newcommand{\ed}{\qed\vspace{3mm}}

\newcommand{\Qcoh}{\operatorname{Qcoh}}

\newcommand{\cha}{\operatorname{char}}
\newcommand{\Per}{\operatorname{Perf}}

\newcommand{\sslash}{\mathbin{/\mkern-6mu/}}

\title{Moduli of curves as moduli of $A_\infty$-structures}
\author{Alexander Polishchuk}
\thanks{Supported in part by NSF grant}

\begin{document}
\begin{abstract}
We define and study the stack $\UU^{ns,a}_{g,g}$ of (possibly singular) projective
curves of arithmetic genus $g$ with $g$ smooth marked points forming an ample non-special divisor.
We define an explicit closed embedding of a natural $\G_m^g$-torsor $\wt{\UU}^{ns,a}_{g,g}$
over $\UU^{ns,a}_{g,g}$
into an affine space and give explicit equations of the universal curve (away from characteristics $2$ and $3$).
This construction can be viewed as a generalization of the Weierstrass cubic
and the $j$-invariant of an elliptic curve to the case $g>1$.
Our main result is that in characteristics different from $2$ and $3$ the moduli space $\wt{\UU}^{ns,a}_{g,g}$ is isomorphic to the moduli space of minimal $A_\infty$-structures on a certain finite-dimensional graded associative algebra $E_g$
(introduced in \cite{FP}). We show how to compute explicitly the $A_\infty$-structure
associated with a curve $(C,p_1,\ldots,p_g)$ in terms of certain canonical generators of the algebra
$\OO(C\setminus\{p_1,\ldots,p_g\})$ and canonical formal parameters at the marked points.
We study the GIT quotients associated with our representation of $\UU^{ns,a}_{g,g}$
as the quotient of an affine scheme by $\G_m^g$ and show that some of the corresponding quotient
stacks give modular compactifications of
$\MM_{g,g}$ in the sense of \cite{Smyth}.
%projective birational models of $M_{g,g}$ in the sense of \cite{FS}. 
We also consider an analogous picture for curves
of arithmetic genus $0$ with $n$ marked points which gives a new presentation of the moduli space of
$\psi$-stable curves (also known as Boggi-stable curves) and its interpretation in terms of $A_\infty$-structures.
\end{abstract}

\maketitle

\bigskip

\centerline{\sc Introduction}

The idea to study algebraic varieties in terms of their derived categories of coherent sheaves
goes back at least to the work of Bondal-Orlov \cite{BO}. An example that stimulated the present
work is the fact that a smooth projective curve of genus $g\ge 2$
can be recovered from the corresponding derived category. 
In the subsequent development of the theory of derived categories it has been realized that a more
flexible framework is obtained by considering their enhancements to dg-categories
or $A_\infty$-categories (see \cite{BK}, \cite{Toen}). 
The main goal of our work, continuing \cite{FP}, is to get a computable invariant of a curve from the
corresponding enhanced derived category. The immediate motivation was the fact that
one can recover the $j$-invariant of an elliptic curve by studying the $A_\infty$-algebra associated 
with the generator $\OO\oplus L$ of the derived category, where $L$ is a line bundle of degree $1$
(see \cite[Sec.\ 5, Thm.\ 3]{Fisette-ell}, \cite[Prop.\ 9]{LP1}, \cite[Sec.\ 5, Thm.\ C]{LP2}). 
To get an analogous picture in the higher genus case
we considered in \cite{FP} the $A_\infty$-algebra coming from a certain generator of the derived category
for a curve of genus $g$ that depends on a choice of $g$ generic points (see \eqref{G-gen-eq} below). 
We showed that already triple products in this $A_\infty$-algebra
give a computable invariant which distinguishes generic curves for $g\ge 6$.
In the present paper we push this analysis further by considering the entire $A_\infty$-algebra structure,
by comparing the corresponding moduli space with the moduli space of curves, and by extending
this picture to singular curves.

In the first part of this paper (Sections \ref{curves-sec} and \ref{equations-sec})
we study by classical methods (not involving $A_\infty$-structures)
%new feature of this paper is the introduction of a more classical algebraic structure
%(that we call a {\it marked algebra of genus $g$}), associated
the moduli stack $\UU^{ns,a}_{g,g}$ of {\it non-special curves},
parametrizing (not necessarily nodal) projective curves $C$ of genus $g$
with $g$ smooth marked points $p_1,\ldots,p_g$, such that the divisor $p_1+\ldots+p_g$ is non-special and ample. 
%Note that this is an algebraic stack but not a Deligne-Mumford stack. 
The main idea is to study the algebra $\OO(C\setminus\{p_1,\ldots,p_g\})$
equipped with the filtration
recording the information about the poles at the marked points. We show (assuming the
characteristic is not $2$ or $3$) that this algebra has a canonical set of generators, 
so that the defining equations can be written explicitly, 
generalizing the Weierstrass cubic equation in genus $1$.
Thus, one can view the coefficients of these defining equations (that can be rescaled by elements of $\G_m^g$)
as the higher genus analog of the $j$-invariant. 
Note that this picture is reminiscent of Petri's analysis of the equations of the canonical curve (see \cite{Petri}),
however, we need fewer genericity assumptions and there are fewer coefficients in our defining equations.
More precisely, the number of coefficients in our picture
is a quadratic polynomial in $g$, whereas the number of coefficients in Petri's relations is a cubic
polynomial in $g$. Note that Petri's analysis has been extended
to some singular canonical curves in \cite{Sch} and to some noncanonical divisors in \cite{AS}.
One feature of our approach is that we impose no apriori restrictions on singularities of the curve (except that it should
be reduced).
Also, we are able to represent the relevant stack of pointed curves as a global quotient by the torus action and
study the corresponding GIT stability conditions (see Section \ref{GIT-sec}).

Now let us formulate our main result in more detail.
We start with a projective curve $C$ over $k$ (where $k$ is a field) of arithmetic genus $g$ with 
$g$ smooth marked points $p_1,\ldots,p_g$ such that $h^0\bigl(C,\OO_C(p_1+\ldots+p_g)\bigr)=1$,
i.e., the divisor $D=p_1+\ldots+p_g$ is non-special. As in \cite{FP}, we consider the
algebra $\Ext^*(G,G)$ for
\begin{equation}\label{G-gen-eq}
G=\OO_C\oplus\OO_{p_1}\oplus\ldots\oplus\OO_{p_g}
\end{equation}
The condition that $D$ is non-special implies that the algebra $\Ext^*(G,G)$ can be identified
with the graded associative algebra $E_g$ (in \cite{FP} we called this algebra $E_{g,g}$)
defined as follows:
\begin{equation}\label{Eg-eq}
E_g=k[Q]/J,
\end{equation}
where $k[Q]$ is the path algebra of the quiver $Q$ with $g+1$ vertices marked as $\OO$, $\OO_{p_1},\ldots,
\OO_{p_g}$ that has for each $i=1,\ldots,g$ one arrow $A_i$ of degree $0$ going from $\OO_{p_i}$ to $\OO$ and one arrow $B_i$ of degree $1$ going from $\OO$ to $\OO_{p_i}$. The ideal $J$ is generated by
the elements 
$$A_iB_iA_i, B_iA_iB_i, A_iB_j,$$
where $i\neq j$. More precisely, the isomorphism $\Ext^*(G,G)\simeq E_g$ depends on a
choice of trivializations of the tangent spaces at the marked ponts. 
Thus, the $A_\infty$-enhancement of the derived category provides a minimal $A_\infty$-structure on 
$E_g$ associated with $(C,p_1,\ldots,p_g)$. In fact, as we showed in \cite{FP}, any such $A_\infty$-structure
is determined (up to an equivalence) by $m_i$ with $i\le 6$, provided $g\ge 2$.

On the other hand, working over a field $k$,
we define the moduli stack of all minimal $A_\infty$-structures on $E_g$, compatible
with its algebra structure, viewed up to a gauge equivalence, and show that it is actually an affine scheme
of finite type over $k$ (provided $\cha(k)\neq 2,3$).
The above construction of $A_\infty$-structures associated with curves works well in families over any affine base. 
Our main result compares the moduli of non-special curves with the moduli of $A_\infty$-structures on $E_g$.

\medskip

\noindent
{\bf Theorem A}. {\it Assume that $g\ge 1$. 
The algebraic stack $\UU^{ns,a}_{g,g}\times\Spec(\Z[1/6])$ is equivalent to the quotient stack
$(\wt{\UU}^{ns,a}_{g,g}\times\Spec(\Z[1/6]))/\G_m^g$, where $\wt{\UU}^{ns,a}_{g,g}\times \Spec(\Z[1/6])$ 
is an affine scheme of finite type over $\Z[1/6]$.
Furthermore, for any field $k$ of characteristic $\neq 2,3$, the affine scheme
$\wt{\UU}^{ns,a}_{g,g}\times\Spec(k)$ is naturally isomorphic to the moduli scheme
of minimal $A_\infty$-structures on the algebra $E_g$ up to a gauge equivalence.}

\medskip

Note that in genus $1$ this result is very close to Theorem C in \cite{LP2}
(and we use some ideas of \cite{LP2} in the proof of our result).
The ampleness condition in the definition of $\UU^{ns,a}_{g,g}$ for $g=1$ becomes equivalent
to the irreducibility of the curve. For arbitrary $g\ge 1$ 
the construction of the $A_\infty$-structure on $E_g$ associated with $(C,p_1,\ldots,p_g)$ works
even if $p_1+\ldots+p_g$ is not ample, however, the ampleness condition is needed in order to have an isomorphism
of moduli spaces in Theorem A (see Remark \ref{ainf-map-rem}).

The proof of Theorem A consists of two parts. First, we construct an explicit embedding
of $\wt{\UU}^{ns,a}_{g,g}$ into an affine space, using a canonical basis in the algebra
$\OO(C\setminus\{p_1,\ldots,p_g\})$ for each curve in our moduli space. 
Then we study the formal neighborhood of the point
in $\wt{\UU}^{ns,a}_{g,g}$ corresponding to the most singular curve in it, the {\it cuspidal curve of
genus $g$}, which is the union of $g$ usual cuspidal curves of genus $1$, joined at their cusps
(see Section \ref{cusp-sec}).
We prove that the relevant deformation functors (for curves and for $A_\infty$-structures)
are isomorphic and then use the $\G_m$-action, where $\G_m\sub\G_m^g$ is the diagonal,
to deduce the result.

Thus, any minimal $A_\infty$-structure on $E_g$ corresponds to some (possibly quite singular) curve.
Since smoothness can be characterized intrinsically in terms of the dg-category of perfect complexes
(see \cite[Prop.\ 3.13]{Lunts-cat-res}),
we get the following compact description of the derived categories associated to smooth curves.

\medskip

\noindent
{\bf Corollary B}. {\it Assume that $k$ is a perfect field of characteristic $\neq 2,3$.
Let $(E^\infty, m_\bullet)$ be a minimal $A_\infty$-algebra such that
$(E^\infty,m_2)\simeq E_g$. Then $E^\infty$ is smooth as an $A_\infty$-algebra if and only if
there exists a smooth projective curve $C$ of genus $g$ such that
$$D^b(C)\simeq \Per(E^\infty),$$ 
where $\Per(E^\infty)$ is the perfect derived category of $E^\infty$.}

\medskip

The above characterization of the perfect derived categories associated to smooth curves
could be useful in trying to establish the homological mirror symmetry connecting such categories
with appropriate Fukaya categories (the other direction of the homological mirror symmetry involving
the Fukaya categories of higher genus curves has been established in \cite{Seidel-curve}, \cite{Efimov}).
However, at present there is no proposal for what should be considered on the symplectic side.

The equivalence of Theorem A is quite explicit and computable:
we show how to calculate the higher products of the $A_\infty$-structure associated with
a curve using explicit homotopies in a natural dg-model. The answer is given in terms of
the coefficients of expansions of the canonical generators of the algebra 
$\OO(C\setminus\{p_1,\ldots,p_g\})$ in terms of certain canonical formal parameters at the points
$p_1,\ldots,p_g$ (these same coefficients define an affine embedding of the moduli space
$\wt{\UU}^{ns,a}_{g,g}$).

The representation of the moduli space $\UU^{ns,a}_{g,g}$, 
in the form $\wt{\UU}^{ns,a}_{g,g}/\G_m^g$ allows to construct birational
models of the coarse moduli space $M_{g,g}$ by 
taking a GIT quotient $\wt{\UU}^{ns,a}_{g,g}\sslash_\chi \G_m^g$ and considering the
component of smoothable curves. More precisely, we show that we get a projective birational model of $M_{g,g}$ 
whenever the character $\chi$ of $\G_m^g$  
belongs to a certain explicit cone in $\R^g$ (see Proposition \ref{GIT-prop}). Furthermore, we show
that for $\chi$ inside a smaller cone all the points of $\wt{\UU}^{ns,a}_{g,g}$, such that the corresponding curve
is smooth, are $\chi$-stable. It would be interesting to study the natural birational maps from $\ov{M}_{g,g}$ to these 
GIT quotients, and for small $g$ compare them to the birational models considered in \cite{FV}.

Of independent interest could be our study of the moduli functor for minimal $A_\infty$-structures up
to a gauge equivalence. We prove (see Theorem \ref{a-inf-moduli-thm} and Corollary \ref{representability-cor}) 
that for any finite-dimensional graded associative
$k$-algebra $E$ such that $HH^1(E)_{<0}=0$, where $HH^i(E)_j$ denotes Hochschild cohomology
(see Section \ref{ainf-gen-sec} for our conventions on the bigrading), 
this moduli functor is representable by an affine scheme over $k$.

Our explicit embedding of the moduli space $\wt{\UU}^{ns,a}_{g,g}$ into an affine space 
leads us to an explicit description 
of the hyperelliptic locus, which is simply the set of fixed points of
the natural involution changing the sign in the trivializations of the tangent spaces at the marked points
(see Theorem \ref{hyperell-thm}).
Note that the open part of this locus corresponding to smooth curves is a $\G_m^g$-torsor over the
configuration space of $2g+2$ distinct points in $\P^1$, $g$ of which are ordered.

We also consider an analogous picture for curves of arithmetic genus $0$
with $n$ marked points. We consider the moduli stack of {\it $\psi$-prestable curves} (this condition
means that each component of a curve contains at least one marked point), and give a natural presentation
of the moduli space as a quotient of the explicit affine scheme by $\G_m^n$ 
(see Theorem \ref{genus-0-moduli-thm}), as well as an interpretation
in terms of moduli of $A_\infty$-structures (in Theorem \ref{genus-0-a-inf-moduli-thm}). 
We show that among the corresponding GIT quotients is the moduli scheme $\ov{M}_{0,n}[\psi]$
of {\it $\psi$-stable} curves, also known as {\it Boggi-stable} curves (see \cite{Boggi}, \cite[Sec.\ 4.2.1]{FS}, \cite[Sec.\ 7.2]{GJM}).
This moduli space was first constructed by Boggi in \cite{Boggi}, following a topological construction by
Kontsevich \cite{Kon}. In \cite{GJM} it was realized as a GIT quotient with respect to an action of $\SL(n)$ on a certain Chow variety. Our approach to $\ov{M}_{0,n}[\psi]$ is quite different and leads to its presentation by explicit equations
in $(\P^{n-3})^n$ (see Corollary \ref{Boggi-cor}).

%Recall that encoding a $g$-pointed curve by the products $(m_3,\ldots,m_6)$ is possible because
%of the Hochschild cohomology vanishing $HH^n_{2-n}(E_g)=0$ for $n>6$ (see \cite[???]{FP}).
%For more detailed understanding of the products $(m_3,\ldots,m_6)$ we need to know other Hochschild
%cohomology groups as well. The relevant computations were made in \cite{Fisette}. 
%Here we present a geometric way to compute these.

%Explicit equations for $\wt{\UU}^{ns,a}_{g,g}$???

The paper is organized as follows.
In Section \ref{curves-sec} we introduce the moduli stacks $\UU^{ns,a}_{g,g}$ and
$\wt{\UU}^{ns,a}_{g,g}$ and study their relation with marked algebras of genus $g$
(which arise as algebras of the form $\OO(C\setminus\{p_1,\ldots,p_g\})$).
We prove the equivalence of the corresponding moduli problems over $\Spec(\Z[1/6])$ and show that it
leads to an explicit embedding of $\wt{\UU}^{ns,a}_{g,g}$ into an affine space 
(see Theorem \ref{moduli-thm}).

In Section \ref{equations-sec} we find a minimal set of generators of the algebra of functions on
$\wt{\UU}^{ns,a}_{g,g}$ and some relations between them. 
Note that some parts of Section \ref{equations-sec} are quite computational, however, 
the results of this Section are not used in the proof of Theorem A.
Then in \ref{GIT-sec} we study the GIT picture for the $\G_m^g$-action on $\wt{\UU}^{ns,a}_{g,g}$,
and 
in \ref{hyperell-sec} we study the hyperelliptic locus. In Section \ref{Petri-sec} 
we work out the relation between our approach
and Petri's analysis of equations for the canonical embedding of a non-hyperelliptic curve.
In particular, for a generic curve $C$ we determine the quadratic equations of $C$ in its canonical embedding
in terms of some of the affine coordinates on $\wt{\UU}^{ns,a}_{g,g}$ (see Proposition \ref{quadratic-rel-prop}).

In Section \ref{Cech-ainf-sec} we give a method for computing $A_\infty$-structures associated with curves,
using a version of Cech resolutions (where open neighborhoods of $p_i$'s are replaced with the formal
neighborhoods). The explicit formulas for the triple products are quite simple. Those for
$m_4$ and $m_5$ are more involved and are given in the Appendix.

In Section \ref{moduli-sec} we establish an isomorphism between $\wt{\UU}^{ns,a}_{g,g}$ and
the moduli space of $A_\infty$-structures stated in Theorem A. First, in \ref{ainf-moduli-sec}
we prove the criterion of representability of the moduli functor of minimal $A_\infty$-structures.
Then in \ref{cusp-sec} we study the geometry of the cuspidal curve of genus $g$.
In \ref{def-theory-sec} we compare the deformation theories for our two
moduli problems at the point corresponding to the cuspidal curve, 
which allows us to finish the proof in \ref{proof-sec}.
In \ref{normal-sec} we point out some consequences for the Hochschild cohomology of $E_g$ and
for the normal forms of $A_n$-structures on $E_g$ for small $n$.

Finally, in Section \ref{genus-0-sec} we consider an analogous picture for the moduli of $\psi$-prestable curves of genus $0$.

{\it Acknowledgments}. I am grateful to Robert Fisette for useful discussions of normal forms
of $A_n$-structures on $E_g$ and to Nick Proudfoot for answering my questions about toric GIT quotients.

\section{Moduli of non-special curves}\label{curves-sec}

\subsection{Non-special curves and marked algebras}

Our main object of interest is the moduli stack $\UU^{ns,a}_{g,g}$ classifying (not necessarily smooth)
projective curves $C$ of genus $g$ with $g$ marked points $p_1,\ldots,p_g$, forming a non-special ample 
divisor.\footnote{The superscripts $ns$ and $a$ stand for {\it non-special} and {\it ample}, since
this is what we require of the divisor $p_1+\ldots+p_g$.} 
Here is a more precise definition.

\begin{defi}\label{moduli-stack-def-1} 
The stack $\UU^{ns,a}_{g,g}$ associates with a noetherian scheme $S$ the groupoid of the following data
and their isomorphisms: a flat proper morphism $\pi:C\to S$ together with $n$ disjoint sections 
$p_1,\ldots,p_g:S\to C$ of $\pi$, such that the following conditions hold for each geometric fiber $C_s$ (where
we denote still by $p_i$ the induced points on $C_s$): 

\noindent
(i) $C_s$ is a connected reduced curve of arithmetic genus $g$, smooth at each marked point $p_i$;

\noindent
(ii) one has $H^1\bigl(C_s,\OO(p_1+\ldots+p_g)\bigr)=0$;

\noindent
(iii) the divisor $D=p_1+\ldots+p_g$ on $C_s$ is ample.
\end{defi}

%Remarks: 
Note that since $\chi\bigl(\OO_{C_s}(p_1+\ldots+p_g)\bigr)=\chi(\OO_{C_s})+g=1$, the condition (ii) is equivalent to
$h^0(C_s,p_1+\ldots+p_g)=1$.
The condition (iii) is equivalent to requiring that there is at least one marked point on each irreducible component of $C$.
By the base change theorem, this implies that the natural map
\begin{equation}\label{ns-map-eq}
\OO_S\to R\pi_*\bigl(\OO_C(D)\bigr),
\end{equation}
where $D=\im(p_1)+\ldots+\im(p_g)$, is an isomorphism. 
Note that in the case of an affine base $S=\Spec(R)$, \eqref{ns-map-eq} becomes an isomorphism
$$R\simeq H^*\bigl(C,\OO(D)\bigr).$$
Below we will often denote the relative divisors $\im(p_i)\sub C$ simply by
$p_i$.

%Thus, if $R$ is a commutative ring then $\UU^{ns,a}_{g,g}(R)$ classifies flat proper families
%$\pi:C\to\Spec(R)$ of curves of arithmetic genus $g$ with $g$ marked
%points $p_i:\Spec(R)\to C$, such that $\pi$ is smooth near the images of $p_i$, the effective divisor
%$D=\im(p_1)+\ldots+\im(p_g)\sub C$ is ample and the natural map 

\begin{defi}\label{moduli-stack-def-2}
We define $\wt{\UU}^{ns,a}_{g,g}\to \UU^{ns,a}_{g,g}$ to be the
$\G_m^g$-torsor corresponding to choices of nonzero tangent vectors at the marked points.
In terms of families $\pi:C\to S$
%$\Spec(R)$ 
this corresponds to choosing a trivialization of each line bundle $L_i:=p_i^*\OO_C(p_i)=\pi_*\bigl(\OO_C(p_i)/\OO_C\bigr)$.
\end{defi}

In other words, for a family in $\wt{\UU}^{ns,a}_{g,g}$ we have in addition to the data $(C,p_1,\ldots,p_n)$ as in
Definition \ref{moduli-stack-def-1}, global sections
$v_i$ of $\pi_*\bigl(\OO_C(p_i)/\OO_C\bigr)$, for each $i=1,\ldots,g$, inducing isomorphisms
\begin{equation}\label{tangent-trivialization-eq}
\OO_S\simeq \pi_*\bigl(\OO_C(p_i)/\OO_C\bigr).
%H^0(C,\OO_C(p_i)/\OO_C).
\end{equation}

We are going to show (away from characteristics $2$ and $3$) 
that $\wt{\UU}^{ns,a}_{g,g}$ is an affine scheme of finite type by relating it to a different kind of data.

For a base commutative ring $R$ 
consider the algebra $R^g=R\times\ldots\times R$ (the direct product of $g$ copies of $R$).
We denote by $e_i\in R^g$, $i=1,\ldots,g$, the natural idempotents.
For $g\ge 1$ let us consider the subalgebra
$$C_{R,g}\sub R^g[t]$$
consisting of polynomials with coefficents in $R^g$ with no linear term in $t$ and with constant
term in $R\sub R^g$. 
We view $C_{R,g}$ as a graded $R$-algebra, where $\deg t=1$.

Let $S_g$ denote the semigroup $\{0\}\cup (\Z_{>0})^g$. We denote elements of $S_g$ as formal effective divisors
$n_1p_1+\ldots+n_gp_g$ with $n_i\in\Z_{>0}$ or $n_1=\ldots=n_g=0$ (where $p_i$ are formal symbols). 
We also view $S_g$ as a lattice with respect to the partial order
$$n_1p_1+\ldots+n_gp_g\le n'_1p_1+\ldots+n'_gp_g \ \text{ if } n_1\le n'_1,\ldots, n_g\le n'_g.$$

\begin{defi}\label{marked-alg-defi} 
A {\it marked algebra of genus $g$} over $R$ is a commutative $R$-algebra $A$ equipped
with the following data:

\noindent
(i) an exhaustive $S_g$-valued $R$-algebra filtration 
$(F_\bullet A)$ on $A$ 
%such that $F_{n_1p_1+\ldots+n_gp_g}A$ are $R$-submodules, and 
such that $F_0 A=R$ and the map
$$n_1p_1+\ldots+n_gp_g\mapsto F_{n_1p_1+\ldots+n_gp_g}A$$
is a morphism of lattices, i.e., 
\begin{equation}\label{lattice-conditions-eq}
F_{\min(\bn,\bn')}A=F_{\bn}A\cap F_{\bn'}A, \ \ F_{\max(\bn,\bn')}A=F_{\bn}A+F_{\bn'}A.
\end{equation}
%is a free $R$-module of rank $1$.

\noindent
(ii) Consider an increasing filtration on $A$ given by 
$$\ov{F}_nA=F_{nD}A \ \text{ for } n\ge 0, \ \text{ where } D=p_1+\ldots+p_g.$$
Then there should be fixed an isomorphism of graded $R$-algebras
\begin{equation}\label{algebra-marking-eq}
\gr^\bullet_{\ov{F}} A\simeq C_{R,g}
\end{equation}
where $\gr^\bullet_{\ov{F}} A$ is the associated graded algebra with respect to $\ov{F}_\bullet A$.
In addition, we require that for $i=1,\ldots,g$ and any $n\ge 1$ the image of the embedding
$$F_{nD+p_i}A/F_{nD}A\hra F_{(n+1)D}A/F_{nD}A$$
gets identified under \eqref{algebra-marking-eq} with $Re_it^{n+1}\sub (C_{R,g})_n$. 
\end{defi}

Note that for any marked algebra $A$ for each $m\ge 1$ we get an identification
\begin{equation}\label{ass-graded-2-eq}
F_{nD+mp_i}A/F_{nD+(m-1)p_i}A\rTo{\sim} F_{(n+m-1)D+p_i}A/F_{(n+m-1)D}A\simeq Re_it^{n+m},
\end{equation}
where the first isomorphism is induced by the natural inclusions.

As the reader might have anticipated, marked algebras of genus $g$ arise from non-special curves.

\begin{lem}\label{marked-alg-curve-lem} 
For a flat proper family of curves $\pi:C\to\Spec(R)$ with $g$ marked points $p_i:\Spec(R)\to C$, where
$\pi$ is smooth along $p_i$'s, 
and trivializations \eqref{tangent-trivialization-eq}, such that the map \eqref{ns-map-eq} is an isomorphism,
consider the algebra $A:=H^0(C\setminus D,\OO)$ with the filtration
$$F_{n_1p_1+\ldots+n_gp_g}=H^0\bigl(C,\OO_C(n_1p_1+\ldots+n_gp_g)\bigr).$$
Then the isomorphisms \eqref{tangent-trivialization-eq} induce an isomorphism \eqref{algebra-marking-eq}
making $A$ into a marked algebra of genus $g$ over $R$.
\end{lem}

\Pf . The vanishing of $H^1\bigl(C,\OO_C(D)\bigr)$ implies the vanishing of $H^1\bigl(C,\OO_C(n_1p_1+\ldots+n_gp_g)\bigr)$
for any $n_1\ge 1,\ldots,n_g\ge 1$. Hence we have exact sequences
\begin{equation}\label{F-n-p-i-eq}
0\to F_{\bn}\to F_{\bn+p_i}\to H^0\bigl(C,\OO_C((n_i+1)p_i)/\OO_C(n_ip_i)\bigr)\to 0
\end{equation}
for every $\bn=n_1p_1+\ldots+n_gp_g\in S_g$, where the third term is isomorphic to $R$.
In particular, for every $\bn\in S_g$ with $n_i>1$ the natural map
$$F_{\bn}\to H^0\bigl(C,\OO_C(n_ip_i)/\OO_C((n_i-1)p_i)\bigr)\simeq R$$
is surjective. Using this we can check the conditions \eqref{lattice-conditions-eq}. 
Namely, the first condition follows from the fact that if we choose $N$ bigger than all $n_i$'s then
$F_\bn$ is the kernel of the natural map
$$F_{Np_1+\ldots+Np_g}\to \bigoplus_{i=1}^g H^0\bigl(C,\OO_C(Np_i)/\OO_C(n_ip_i)\bigr).$$
The second condition follows from the equality 
$$F_\bn=\sum_{i=1}^g F_{\bn^{(i)}},$$ 
where $\bn^{(i)}$ has $n_i$ as $i$th component and $1$'s as other components. The latter equality is deduced
by induction using \eqref{F-n-p-i-eq}.
%This easily implies that $\bn\mapsto F_{\bn}A$ is a morphism of lattices. 

Next, for every $n>1$ we have an exact sequence
$$0\to F_{(n-1)D}\to F_{nD}\to H^0\bigl(C,\OO_C(nD)/\OO_C((n-1)D)\bigr)\to 0,$$
while $F_D=R$. Together with the isomorphisms \eqref{tangent-trivialization-eq} this
leads to an identification of $\gr^\bullet_{\ov{F}}A$ with $C_{R,g}$, satisfying the required properties.
\ed

\subsection{Canonical generators for marked algebras and the moduli spaces}\label{marked-moduli-sec}

\begin{lem}\label{marked-alg-main-lem} 
Let $A$ be a marked algebra of genus $g$ over $R$. Assume that $6$ is invertible 
in $R$.

\noindent
(i) There exist unique elements
$f_i\in F_{D+p_i}A$ and $h_i\in F_{D+2p_i}A$, where $i=1,\ldots,g$, called {\em canonical generators of} $A$,
such that 
$$f_i\mod F_DA = e_it^2, \ \ h_i\mod F_{2D}A= e_it^3,$$
$$h_i^2-f_i^3\in F_{3D}A, \ \ f_ih_i^2-f_i^4\in F_{4D}A.$$

\noindent
(ii) The elements $(f_i^n, f_i^nh_i)_{i=1,\ldots, g; n\ge 0}$ form an $R$-basis in $A$.

\noindent
(iii) The algebra $A$ is generated over $R$ by $(f_i)$ and $(h_i)$ and has defining relations of the form
\begin{equation}\label{main-eq-ff}
f_if_j=\a_{ji}h_i+\a_{ij}h_j+\ga_{ji}f_i+\ga_{ij}f_j+\sum_{k\neq i,j}c_{ij}^k f_k+a_{ij},
\end{equation}
\begin{equation}\label{main-eq-fh}
f_ih_j=d_{ij}f_j^2+t_{ji}h_i+v_{ij}h_j+r_{ji}f_i+\de_{ij}f_j+\sum_{k\neq i,j} e_{ij}^kf_k+b_{ij},
\end{equation}
\begin{equation}\label{main-eq-hh}
h_ih_j=\b_{ji}f_i^2+\b_{ij}f_j^2+\vareps_{ji}h_i+\vareps_{ij}h_j+\psi_{ji}f_i+\psi_{ij}f_j+
\sum_{k\neq i,j}l_{ij}^kf_k+u_{ij},
\end{equation}
\begin{equation}\label{main-eq-h2f3}
h_i^2=f_i^3+\sum_{j\neq i}g_i^jh_j+\pi_if_i+\sum_{j\neq i}k_i^jf_j+s_i,
\end{equation}
where $i\neq j$ (the coefficients are some elements of $R$).

\noindent
(iv) If $(f_i,h_i)_{i=1,\ldots,g}$ is a collection of elements, such that $f_i\in F_{2D}A$, $h_i\in F_{3D}A$, 
$$f_i\mod F_DA = e_it^2, \ \ h_i\mod F_{2D}A= e_it^3,$$
and such that relations of the form \eqref{main-eq-ff}--\eqref{main-eq-h2f3} are satisfied, then
$(f_i,h_i)$ are the canonical generators of $A$.

\noindent
(v) Let $A$ and $A'$ be marked algebras of genus $g$ over $R$, where $6$ is invertible in $R$.
Then any isomorphism of $R$-algebras $A\rTo{\sim}A'$, sending the canonical generators of $A$ to those of $A'$,
is an isomorphism of marked algebras.
\end{lem}

\Pf . (i) Let $\wt{f}_i\in F_{D+p_i}A$ and $\wt{h}_i\in F_{D+2p_i}$ 
be some liftings of $e_it^2$ and of $e_it^3$, respectively (with respect to identifications \eqref{ass-graded-2-eq}).
Then we have 
$$\wt{f}_i^3\equiv e_it^6\equiv \wt{h}_i^2 \mod F_{5D}A.$$
Hence, 
$$\wt{h}_i^2-\wt{f}_i^3\in F_{5D}A\cap F_{3D+3p_i}A=F_{3D+2p_i}.$$
The isomorphism \eqref{ass-graded-2-eq} shows that $\wt{h_i}$ projects to a generator of 
$F_{D+2p_i}A/F_{D+p_i}A\simeq R$. Hence, $F_{D+2p_i}A$ (resp., $F_{D+p_i}$) has an $R$-basis 
$$1, \wt{f}_i, \wt{h}_i \ (\text{resp.,} \ 1, \wt{f}_i).$$
Thus, any other liftings of $e_it^2$ and $e_it^3$ have form
$$f_i=\wt{f}_i+a_i, \ \ h_i=\wt{h}_i+b_i\wt{f}_i+c_i.$$
We have
$$h_i^2-f_i^3\equiv \wt{h}_i^2-\wt{f}_i^3+2b_i\wt{h}_i\wt{f}_i+(b_i^2-3a_i)\wt{f}_i^2\mod F_{3D}A.$$
Using \eqref{ass-graded-2-eq} we see that the elements $(\wt{h}_i\wt{f_i}, \wt{f}_i^2)$ project 
to a basis of $F_{3D+2p_i}A/F_{3D}A$.
It follows that we can choose $a_i$ and $b_i$ uniquely so that
$h_i^2-f_i^3\in F_{3D}A$. There still remains an ambiguity in adding a constant to $h_i$. 
We claim that the condition $f_ih_i^2-f_i^4\in F_{4D}A$ uniquely resolves this ambiguity.
Indeed, if we replace $h_i$ by $h_i+c$ then 
$f_ih_i^2-f_i^4\in F_{4D+p_i}A$ modulo $F_{4D}A$ changes to
$$f_ih_i^2-f_i^4+2cf_ih_i\mod F_{4D}A.$$
Now our claim follows from the fact that $f_ih_i$ projects to a generator of $F_{4D+p_i}A/F_{4D}A\simeq R$.

\noindent
(ii) This follows immediately from the fact that the initial parts of these elements with respect to the filtration $\ov{F}$ form a basis
in $\gr^\bullet_{\ov{F}} A$.

\noindent
(iii) By (ii), $A$ is generated by $(f_i)$ and $(h_i)$. To see that relations of the form \eqref{main-eq-ff},
\eqref{main-eq-fh} and \eqref{main-eq-hh} hold for any $i\neq j$, we just observe that
$$f_if_j\in F_{2D+p_i+p_j}A, \ f_ih_j\in F_{2D+p_i+2p_j}A, \ h_ih_j\in F_{2D+2p_i+2p_j}A$$
and use the fact that the $g+5$ elements
$$1, f_1, \ldots, f_g, h_i, h_j, f_j^2, f_i^2$$
form an $R$-basis of $F_{2D+2p_i+2p_j}A$ with the first $g+4$ (resp., $g+3$) giving an $R$-basis
of $F_{2D+p_i+2p_j}A$ (resp., $F_{2D+p_i+p_j}A$). 
Similarly, the condition $h_i^2-f_i^3\in F_{3D}A$ implies a relation of the form
$$h_i^2-f_i^3=c_ih_i+\sum_{j\neq i}g_i^jh_j+\pi_if_i+\sum_{j\neq i}k_i^jf_j+s_i.$$
Further, we have 
$$0\equiv f_ih_i^2-f_i^4\equiv c_if_ih_i \mod F_{4D}A.$$
Since $f_ih_i$ projects to the basis element $e_it^5$ in $F_{5D}A/F_{4D}A$, this implies that $c_i=0$.

Let $A'$ be the quotient of the polynomial algebra $R[f_1,\ldots,f_g,h_1,\ldots,h_g]$ by
the relations \eqref{main-eq-ff}--\eqref{main-eq-h2f3}. We claim that the elements of
$$\BB=\{f_i^n, f_i^nh_i \ |\ i=1,\ldots,g;n\ge 0\}$$ 
span $A'$ as an $R$-module.
Indeed, let us consider the degree lexicographical order on monomials in $(f_i)$, $(h_i)$,
given by 
\begin{equation}\label{deg-lex-order}
\deg f_i=2, \ \deg h_i=3, \ h_1>\ldots>h_g>f_1>\ldots>f_g.
\end{equation}
Then our relations show that any monomial except for the elements of $\BB$ can be expressed in $A'$
in terms of smaller ones, which proves our claim.
Since the natural homomorphism $A'\to A$ sends elements of $\BB$ to an $R$-basis of $A$,
it follows that this homomorphism is an isomorphism.

\noindent
(iv) The relation \eqref{main-eq-h2f3} immediately shows that $h_i^2-f_i^3\in F_{3D}A$. We also see from 
\eqref{main-eq-fh} that $f_ih_j\in F_{4D}A$ for $i\neq j$. Hence, using \eqref{main-eq-h2f3} again we see
that $f_i(h_i^2-f_i^3)\in F_{4D}A$. It remains to check that $h_j\in F_{D+2p_j}A$ for each $j$. 
Let $(f^{\can}_i, h^{\can}_i)$ be the canonical generators of $A$. It is enough to show that $h_j-h^{\can}_j\in F_{D+p_j}A$.
Since we already know that $h_j-h^{can}_j\in F_{2D}A$, it suffices to check that the $i$th component of 
$$h_j-h^{\can}_j\mod F_DA\in F_{2D}A/F_DA\simeq (C_{R,g})_2=\bigoplus_{i=1}^g R\cdot e_it^2$$
vanishes for every $i\neq j$. Since $f_i$ projects to $e_it^2\in (C_{R,g})_2$,
it is enough to check that $f_i(h_j-h^{\can}_j)\in F_{3D+p_j}A$ for each $i\neq j$.
But \eqref{main-eq-fh} shows that $f_ih_j\in F_{3D+p_j}A$ and similarly $f_i^{\can}h_j^{\can}\in F_{3D+p_j}A$.
Since $f_i$ differs from $f_i^{\can}$ by a constant, we also have $f_ih^{\can}_j\in F_{3D+p_j}A$. Hence,
$f_ih_j-f_ih^{\can}_j\in F_{3D+p_j}A$, as required.

\noindent
(v) It is enough to show that the filtration $F_\bullet A$ can be recovered from the canonical generators $(f_i,h_i)$.
By the lattice property it is enough to recover the submodules $F_{nD+mp_i}A$ for all $n\ge 1, m\ge 0$. 

Let us rename the elements of the basis $\BB$,  different from $1$, as follows:
\begin{equation}\label{b-i-n-eq}
b_i[n]:=\begin{cases} f_i^{n/2}, & n \text{ is even},\\ f_i^{(n-3)/2}h_i, & n \text{ is odd},\end{cases}
\end{equation}
where $n\ge 2$, $i=1,\ldots,g$. Note that $b_i[2k]\in F_{kD+kp_i}A$ while $b_i[2k+1]\in F_{kD+(k+1)p_i}A$.

First, as in (ii) we see that $(1, b_i[m])$, where $m\le n$, $i=1,\ldots,g$, is an $R$-basis of 
$F_{nD}A$. Next, using the isomorphisms \eqref{ass-graded-2-eq}
we see that the elements $b_i[n+1],\ldots,b_i[n+m]$ project to a basis of $F_{nD+mp_i}A/F_{nD}A$, provided
$m\le n$. Thus, we recovered $F_{nD+mp_i}A$ from the canonical generators, in the range $m\le n$.
For general $(m,n)$ and any $N>1$ we claim that 
$$F_{nD+mp_i}A=\{x\in F_{(n+m)D}A \ |\ x\cdot F_{ND}A\sub F_{(N+n)D+mp_i}A\}.$$
Then taking $N\ge m$ would allow to recover $F_{nD+mp_i}A$ from the canonical generators.
To prove the claim we use the induction on $m$. The case $m=0$ is a tautology. For $m>0$ let 
$x\in  F_{(n+m)D}A$ be such that $xF_{ND}\sub F_{(N+n)D+mp_i}A$. Considering the image
of $x$ in $F_{(n+m)D}A/F_{(n+m-1)D}A$ and using \eqref{algebra-marking-eq}, we see that
$x\in F_{(n+m-1)D+p_i}A$. Using the isomorphism \eqref{ass-graded-2-eq} we can choose
$x'\in F_{nD+mp_i}A$ such that $x\equiv x'\mod F_{(n+m-1)D}A$.
Then 
$$(x-x')F_{ND}A\sub F_{(N+n)D+mp_i}A\cap F_{(N+n+m-1)D}A=F_{(N+n)D+(m-1)p_i}A.$$
By the induction assumption, this implies that $x-x'\in F_{nD+(m-1)p_i}A$, hence
$x\in F_{nD+mp_i}A$, as claimed.
\ed

\begin{lem}\label{grobner-lem} 
Let $R$ be a commutative ring.

\noindent
(i) A commutative $R$-algebra $B$ generated by $f_1,\ldots,f_g,h_1,\ldots,h_g$ with
defining relations \eqref{main-eq-ff}--\eqref{main-eq-h2f3} has elements
$(f_i^n, f_i^nh_i)_{i=1,\ldots,g, n\ge 0}$ as an $R$-basis if and only if the coefficients
\begin{equation}\label{big-generators-eq}
(\a_{ij},\b_{ij},\ga_{ij},\de_{ij},\vareps_{ij},\psi_{ij},\pi_i,d_{ij},t_{ij},v_{ij},r_{ij},c_{ij}^k,
e_{ij}^k,l_{ij}^k,g_i^j,k_i^j,a_{ij},b_{ij},u_{ij},s_i)
\end{equation}
satisfy certain universal polynomial relations with integer coefficients.

\noindent
(ii) For an $R$-algebra $B$ satisfying the equivalent conditions in (i)
the $R$-submodules 
$$\ov{F}_nB:=R+\sum_{m\le n, i=1,\ldots,g} R\cdot b_i[m],$$
where $b_i[m]$ are given by \eqref{b-i-n-eq},
define an algebra filtration, i.e., $\ov{F}_nB\cdot \ov{F}_{n'}B\sub \ov{F}_{n+n'}B$.
\end{lem}

\Pf . (i) Note that the property in question holds if and only if
the four types of elements in $R[f_1,\ldots,f_g,h_1,\ldots,h_g]$, given by the difference of the left and right
sides in \eqref{main-eq-ff}--\eqref{main-eq-h2f3},
form a Gr\"obner basis of the ideal they generate
with respect to the degree lexicographical order defined by \eqref{deg-lex-order}.
By Buchberger's Criterion, to get the corresponding relations between the coefficients one has to look at two ways to apply the relations to $$f_if_jf_k, \ f_if_jh_k, \ f_if_jh_j,\ f_ih_j^2,\  f_ih_jh_k, \ f_ih_ih_j,\ h_ih_jh_k, \ h_ih_j^2$$
and then compare the coefficients of the elements $(f_i^n, f_i^nh_i)$ (see e.g, \cite[Thm.\ 15.8]{Eis-CA}).

\noindent
(ii) It is enough to check that for any $n\ge 2$ and any $i,j$ one has
$f_jb_i[n]\in \ov{F}_{n+2}B$ and $h_jb_i[n-1]\in \ov{F}_{n+2}B$ (for $n\ge 3$). 
For $n\le 3$ this follows from the relations \eqref{main-eq-ff}--\eqref{main-eq-h2f3}.

For $n>3$ we have $b_i[n]=f_ib_i[n-2]$. 
Hence, expressing $f_jf_i$, where $i\neq j$, 
via \eqref{main-eq-ff} and using the induction assumption
we get 
$$f_jb_i[n]=f_jf_ib_i[n-2]\in \ov{F}_{n+1}B$$
(note that for $j=i$ we have $f_ib_i[n]=b_i[n+2]\in \ov{F}_{n+2}B$). 

Now let us check that $h_jb_i[n-1]\in  \ov{F}_{n+2}B$.
If $n=4$ then $h_jb_i[3]=h_jh_i$ is in $\ov{F}_6B$
by \eqref{main-eq-hh}, \eqref{main-eq-h2f3}, so we can assume that $n>4$.
For $j\neq i$ we write $b_i[n-1]=f_ib_i[n-3]$ and express $h_jf_i$ in terms of $f_j^2$ and $(f_k,h_k)$, using
\eqref{main-eq-fh}. Note that by the induction assumption, 
$$f_j^2b_i[n-3]\in f_j\ov{F}_{n-1}B\sub \ov{F}_{n+1}B.$$ 
Hence, we derive that $h_jf_ib_i[n-3]\in \ov{F}_{n+2}B$, as required.
It remains to check that $h_ib_i[n-1]\in \ov{F}_{n+2}B$. If $n$ is odd then $h_ib_i[n-1]=b_i[n+2]\in\ov{F}_{n+2}B$.
Now assume that $n$ is even. Then we have
$h_ib_i[n-1]=h_i^2b_i[n-4]$. The relation \eqref{main-eq-h2f3} expresses $h_i^2$ in terms of $f_i^3$ and $(f_k,h_k)$.
By the induction assumption, we have
$$f_i^3b_i[n-4]\in f_i^2\ov{F}_{n-2}B\sub f_i\ov{F}_nB,$$
which is in $\ov{F}_{n+2}B$, as we just proved. Thus, we deduce that $h_ib_i[n-1]\in \ov{F}_{n+2}B$, as required.
\ed

\begin{defi}
We define $\SS_g$ to be the affine scheme over $\Z[1/6]$ defined by the universal
polynomial relations of Lemma \ref{grobner-lem} between the coordinates \eqref{big-generators-eq}.
\end{defi}
%We will show later (see ???) that these relations are such that all the generators can be expressed
%in terms of $\a_{ij}, \b_{ij}, \ga_{ij}, \vareps_{ij}$ and $\pi_i$ (which satsify further relations).

\begin{thm}\label{moduli-thm} 
The following three stacks over $\Spec(\Z[1/6])$ are isomorphic: 

\noindent
(i) $\wt{\UU}^{ns,a}_{g,g}\times \Spec(\Z[1/6])$;

\noindent
(ii) the moduli stack $\MA_g$ of marked algebras of genus $g$ (over $\Spec(\Z[1/6])$); and 

\noindent
(iii) the affine scheme $\SS_g$.

Moreover, the isomorphisms are compatible with the following $\G_m^g$-actions: 

\noindent
(i) on $\wt{\UU}^{ns,a}_{g,g}$---
by rescaling the tangent vectors $v_i$ at the marked points by
$$v_i\mapsto \la_i^{-1}v_i;$$ 

\noindent
(ii) on $\MA_g$---by composing the
isomorphism \eqref{algebra-marking-eq} with the $R$-automorphisms of $C_{R,g}$ sending $te_i$ to 
$\la_ite_i$; 

\noindent
(iii) on $\SS_g$---by rescalings of the structure constants induced by the substitutions
\begin{equation}\label{Gm-fh-eq}
f_i\mapsto \la_i^{-2}f_i, \ \ h_i\mapsto \la_i^{-3}h_i.
\end{equation}
\end{thm}

\Pf . Let $R$ be a Noetherian commutative ring such that $6$ is invertible in $R$.
By Lemma \ref{marked-alg-curve-lem}, we have a functor of groupoids
\begin{equation}\label{H0-functor-eq}
H^0:\wt{\UU}^{ns,a}_{g,g}(R)\to \MA_g(R)
\end{equation}
associating with a family of curves $\pi:C\to\Spec(R)$ 
%with $g$ marked points and trivializations \eqref{tangent-trivialization-eq}, which is a family
in  $\wt{\UU}^{ns,a}_{g,g}(R)$, the marked algebra $H^0(C\setminus\{p_1,\ldots,p_g\},\OO)$.
Also, the construction of the canonical generators (see Lemma \ref{marked-alg-main-lem}) gives a functor 
\begin{equation}\label{can-functor-eq}
\can: \MA_g(R)\to \SS_g(R).
\end{equation}

Next, we are going to construct a functor (really a map, since $\SS_g(R)$ has only identity morphisms)
\begin{equation}\label{Proj-functor-eq}
\Proj: \SS_g(R)\to \wt{\UU}^{ns,a}_{g,g}(R).
\end{equation}
An $R$-point of $\SS_g$ corresponds to an $R$-algebra $A$ with defining 
relations of the form  \eqref{main-eq-ff}--\eqref{main-eq-h2f3},
such that $(1, b_i[n])_{n\ge 2,i=1,\ldots,g}$ (see \eqref{b-i-n-eq}) is an $R$-basis in $R$. 
Let us consider the algebra filtration $(\ov{F}_\bullet A)$ defined in Lemma \ref{grobner-lem}(ii).
The relations \eqref{main-eq-ff}--\eqref{main-eq-h2f3} imply that we have a natural isomorphism
$$\gr^\bullet_{\ov{F}}A\simeq C_{R,g},$$
sending the image of $f_i$ (resp., $h_i$) to $e_it^2$ (resp., $e_it^3$).
Now let us consider the Rees algebra associated with the filtration $(\ov{F}_nA)$:
$$\RR A:=\bigoplus_{n\ge 0} \ov{F}_n A.$$
This is a graded algebra, so we can consider the corresponding projective scheme 
$$C:=\Proj(\RR A)$$
over $\Spec(R)$. Let $T$ denote the element $1\in \ov{F}_1 A=(\RR A)_1$. Then the open subset $T\neq 0$
in $C=\Proj(\RR A)$ is the affine scheme $\Spec(A)$, while the complementary closed subset is the subscheme
$\Proj(\gr^\bullet_{\ov{F}}A)$ which is isomorphic to $\Proj(C_{R,g})$. 
We have natural sections $p_i:\Spec(R)\to \Proj(C_{R,g})$, $i=1,\ldots,g$, such that
$\Proj(C_{R,g})$ is the disjoint union of the images of $p_1,\ldots,p_g$. 

Let $F_i$ (resp., $H_i$) be the element $f_i\in \ov{F}_2 A=(\RR A)_2$ (resp.,
$h_i\in \ov{F}_3 A=(\RR A)_3$).
The algebra $\RR A$ is generated by the elements $(F_i, H_i, T)$, where $i=1,\ldots,g$,
with the following defining relations,
which are homogeneous versions of relations \eqref{main-eq-ff}--\eqref{main-eq-h2f3}:
\begin{equation}\label{main-eq-hom-ff}
F_iF_j=\a_{ji}H_iT+\a_{ij}H_jT+\ga_{ji}F_iT^2+\ga_{ij}F_jT^2+\sum_{k\neq i,j}c_{ij}^k F_kT^2+a_{ij}T^4,
\end{equation}
\begin{equation}\label{main-eq-hom-fh}
F_iH_j=d_{ij}F_j^2T+t_{ji}H_iT^2+v_{ij}H_jT^2+r_{ji}F_iT^3+\de_{ij}F_jT^3+\sum_{k\neq i,j} e_{ij}^kF_kT^3+b_{ij}T^5,
\end{equation}
\begin{equation}\label{main-eq-hom-hh}
H_iH_j=\b_{ji}F_i^2T^2+\b_{ij}F_j^2T^2+\vareps_{ji}H_iT^3+\vareps_{ij}H_jT^3+\psi_{ji}F_iT^4+
\psi_{ij}F_jT^4+\sum_{k\neq i,j}l_{ij}^kF_kT^4+u_{ij}T^6,
\end{equation}
\begin{equation}\label{main-eq-hom-h2f3}
H_i^2=F_i^3+\sum_{j\neq i}g_i^jH_jT^3+\pi_iF_iT^4+\sum_{j\neq i}k_i^jF_jT^4+s_iT^6,
\end{equation}
where $i\neq j$.
On the image $p_i$ we have $T=F_j=H_j=0$ for $j\neq i$ and $F_i\neq 0$, $H_i\neq 0$.
Hence, the section $H_i/F_i$ trivializes $\OO_C(1)$ near the image of $p_i$. On the other hand,
the section $T$ trivializes it on the subset $(T\neq 0)$. Hence, the sheaf
$\OO_C(1)$ is invertible. Similarly we see that $\OO_C(n)\simeq \OO_C(1)^{\ot n}$, so
$\OO_C(1)$ is ample. Thus, the divisor $(T=0)=D=\im(p_1)+\ldots+\im(p_g)$ on $C$ is ample.

Next, we want to check that the projection $\pi:C\to\Spec(R)$ is flat of relative dimension $1$
and smooth near the images of $p_1,\ldots,p_g$.
We have $C\setminus D=\Spec(A)$, and $A$ is a free $R$-module, 
so $\pi$ is flat on $C\setminus D$. To see that $\pi$ is flat near $D$ consider the open subset
$F_iH_i\neq 0$. Let $A_i$ be the degree $0$ part in the localization $(\RR A)_{F_iH_i}$.
On this open subset $D$ is given as the vanishing locus of the function 
$$t_i:=TF_i/H_i\in A_i,$$
which is not a zero divisor in $A_i$, since $T$ is not a zero divisor in $\RR A$.
We already know that the localization $(A_i)_{t_i}$ is flat over $R$. On the other hand, we have
an isomorphism $A_i/(t_i)\simeq R$. This implies that $A_i$ is flat over $R$ by Lemma \ref{flat-lem}
below. This also shows that $\pi$ has relative dimension $1$ near $p_1,\ldots,p_g$.
Next, let us show that $\pi:C\to \Spec(R)$ is smooth near $p_1,\ldots,p_g$.
Consider the open subset $F_iH_i\neq 0$ which is a neighborhood of $p_i$ in $C$. Then $p_i$ is the
intersection of $D$ with this open subset, so $t_i$ generates the ideal of $p_i$ on this open subset.
Since $t_i$ is a nonzerodivisor, this implies the required smoothness near $p_i$.

Next, we claim that viewing elements of the algebra $A$ as functions on $C\setminus D$ and considering
the maximal order of poles at $p_1,\ldots,p_g$ we recover the filtration $(\ov{F}_\bullet A)$.
Indeed, as we have seen above we can use $t_i=TF_i/H_i$ as a local parameter near $p_i$.
The equation 
$$\frac{T^2}{F_i}=t_i^2\cdot \frac{H_i^2}{F_i^3}$$
shows that $f_i=F_i/T^2$ has a pole of order $2$ near $p_i$.
Since $h_i^2/f_i^3=H_i^2/F_i^3$ is invertible near $p_i$, it follows that $h_i$ has a pole of order $3$
near $p_i$. Note that $f_j/f_i=F_j/F_i$ vanishes at $p_i$ while $h_j/h_i=H_j/H_i\in (t_i)^2$ (the latter follows from 
\eqref{main-eq-hom-fh}). Hence, $f_j$ and $h_j$ for $j\neq i$
have poles of order at most $1$ at $p_i$. Next, the equation \eqref{main-eq-hom-h2f3} implies that
$$\frac{h_i^2}{f_i^3}(p_i)=\frac{H_i^2}{F_i^3}(p_i)=1,$$
so we deduce that the Laurent decompositions of $f_i$ and $h_i$ in terms of the parameter $t_i$
start with $t_i^{-2}$ and $t_i^{-3}$, respectively. 
%Of course, the same argument works for any point $p_i$.
Since $(f_i^n, f_i^nh_i)$ form an $R$-basis of $A$, this gives an identification of the filtration
$(\ov{F}_\bullet A)$ with the filtration given
by polar conditions at $p_1,\ldots,p_g$. 

Since we can characterize $H^0\bigl(C,\OO(nD)\bigr)$ inside the algebra $A$ by polar conditions
at $p_1,\ldots,p_g$, we deduce that $R=H^0(C,\OO)=H^0\bigl(C,\OO(D)\bigr)$ and 
that $H^0\bigl(C,\OO(nD)\bigr)$ is a free $R$-module of rank $(n-1)g+1$ for $n\ge 1$.
By flatness of the family the analogous assertion is true for geometric fibers $C_s$ over $\Spec(k)$
(with $R$ replaced by $k$).
It follows that these fibers $C_s$ are projective curves of arithmetic genus $g$ (since $\OO(D)$ is ample),
and hence they have $H^1\bigl(C_s,\OO(D)\bigr)=0$.
Finally, the local parameters $t_1,\ldots,t_g$ induce the required trivializations 
\eqref{tangent-trivialization-eq},
so our family is an object of $\wt{\UU}^{ns,a}_{g,g}(R)$.
This finishes the construction of the map \eqref{Proj-functor-eq}.
Note that this construction also equips $A$ with the structure of a marked algebra, and
Lemma \ref{marked-alg-main-lem}(iv) implies that $(f_i,h_i)$ are the canonical generators with respect to this structure.
Thus, we deduce that the composition
$$\can\circ H^0 \circ\Proj:\SS_g(R)\to\SS_g(R)$$
is the identity.

On the other hand, if the algebra $A$ is associated with a family of curves $(C,p_1,\ldots,p_g)$
in $\UU^{ns,a}_{g,g}(R)$ then we claim that the above construction recovers the original family of curves. Indeed, 
we have $\RR A=\bigoplus_n H^0(C,\OO(nD))$, where $D=p_1+\ldots+p_g$. Since $D$ is ample, the natural
morphism $C\to \Proj(\RR A)$ is an isomorphism. One can easily check that it is compatible with the marked points.
Thus, we get that the composition
$$\Proj\circ\can\circ H^0: \UU^{ns,a}_{g,g}(R)\to \UU^{ns,a}_{g,g}(R)$$
is isomorphic to the identity functor.

Finally, we observe that by Lemma \ref{marked-alg-main-lem}(v), the functor $\can:\MA_g(R)\to \SS_g(R)$
is fully faithful. Hence, we deduce that the composition 
$$H^0\circ\Proj\circ\can:\MA_g\to \MA_g$$
is also the identity. Thus, all three functors $H^0$, $\can$ and $\Proj$ are equivalences.

The compatibility with the $\G_m^g$-actions is a straightforward check.
\ed

We have used the following result.

\begin{lem}\label{flat-lem} Let $A$ be a commutative $R$-algebra, and let $t\in A$ be a nonzero divisor.
Assume that $A_t$ and $A/(t)$ are flat over $R$. Then $A$ is also flat over $R$.
\end{lem}

\Pf . Consider the exact sequence
$$0\to A\rTo{t} A\rTo{} A/(t)\to 0.$$
Let $M$ be any $R$-module.
Since $A/(t)$ is a flat $R$-module, the long exact sequence of $\Tor(?,M)$ shows
that the multiplication by $t$ induces an isomorphism
$\Tor_i^R(A,M)\to \Tor_i^R(A,M)$ for any $i\ge 1$.
Since $A_t=\varinjlim (A\stackrel{t}{\to} A\stackrel{t}{\to} A\to\ldots)$, this implies the isomorphism
$$\Tor_i^R(A,M)\simeq \Tor_i^R(A_t,M)=0$$
for $i\ge 1$ since $A_t$ is flat over $R$.
\ed

\begin{rems} 1. We normalize our $\G_m^g$-actions
so that the induced action of the diagonal $\G_m\sub\G_m^g$ has non-negative weights on the algebra
of functions on $\SS_g$, where the action on functions is given by the operators $(\la^{-1})^*$. This diagonal
$\G_m$-action will play an important role in the proof of Theorem A.

\noindent
2. Let $C$ be the universal curve over $\wt{\UU}^{ns,a}_{g,g}$.
The action of $\G_m^g$ on $C\setminus D$ given by \eqref{Gm-fh-eq}
extends to a $\G_m^g$-action on $C$, compatible with the
$\G_m^g$-action on the moduli space, so that the projection and the sections $p_i$ are $\G_m^g$-equivariant. 

\noindent
3. It is easy to see from the proof of Theorem \ref{moduli-thm} that for any $(C,p_1,\ldots,p_g)\in \UU^{ns,a}_{g,g}(k)$
the divisor $3D$ is actually very ample on $C$.
\end{rems}

\begin{defi}\label{cusp-curve-def}
Let $k$ be a field. We define the {\it cuspidal curve of genus} $g$ over $k$, 
$C^{\cusp}_g$ as the union of $g$ usual (projective)
cuspidal curves $C_i$ of arithmetic genus $1$, glued along their singular points $q_i\in C_i$ in such a way
that the local ring of $C^{\cusp}_g$ at the singular point is the subring of 
$\prod_{i=1}^g\OO_{C_i,q_i}$ consisting of $(f_i)$ such that $f_i(q_i)=f_j(q_j)$. 
\end{defi}

This curve can be viewed as the $k$-point of the moduli space $\wt{\UU}^{ns,a}_{g,g}$ corresponding to the marked algebra
$A^{\cusp}$ over $k$ with generators
$f_1,\ldots,f_g$, $h_1,\ldots,h_g$ and defining relations
$$h_i^2=f_i^3, \ \ f_if_j=f_ih_j=h_ih_j=0, \ i\neq j.$$
Note that $A^{\cusp}$ is the algebra of functions on the affine part of $C^{\cusp}_g$ obtained by deleting points at infinity $p_1,\ldots,p_g$, where 
$p_i\in C_i$ (as in Theorem \ref{moduli-thm}, we can recover $C^{\cusp}_g$ from $A^{\cusp}$ as $\Proj$ of the Rees 
algebra of $A^{\cusp}$). We will refer to this point of $\wt{\UU}^{ns,a}_{g,g}$ as the point corresponding to
the cuspidal curve $C^{\cusp}_g$ and denote it by $[C^{\cusp}_g]$. Note that under the isomorphism
$\wt{\UU}^{ns,a}_{g,g}\simeq \SS_g$ it corresponds to the point where all the coordinates \eqref{big-generators-eq} vanish.
In particular, this point is invariant under $\G_m^g$, and we obtain the induced $\G_m^g$-action on 
$C^{\cusp}_g$ preserving the $g$ points at infinity.

%We will use the following presentation of the cuspidal curve as $\Proj$: 
%$C^{\cusp}_g=\Proj(\RR A^{\cusp})$, where $A^{\cusp}$ is the commutative $k$-algebra 
%and $\RR A^{\cusp}$ is the Rees algebra associated with the filtration $F_\bullet A^{\cusp}$, where
%$F_n$ is spanned by $f_i^m$ with $2m\le n$ and $h_if_i^m$ with $3+2m\le n$.
%In other words, this is the curve that corresponds to the algebra $A^{\cusp}$, viewed as a marked algebra of
%genus $g$ (see Section \ref{curves-sec}). 
%In this presentation $C^{\cusp}_g$ comes equipped
%with $g$ smooth points $p_1,\ldots,p_g$ such that $D=p_1+\ldots+p_g$ is the zero locus
%of the section $1\in F_1$ of $\OO(1)$ (in particular, $D$ is ample). Furthermore, 
%we can view $f_i$ and $h_i$ as functions on the complement $C^{\cusp}_g\setminus D$ and
%$f_i/h_i$ induces a trivialization $v_i$ of the tangent space at $p_i$. 

%Corollary about automorphisms.
%Let $\varphi:C\to C$ be an automorphism fixing each of the marked points $p_1,\ldots,p_g$.
%Then $\varphi$ is determined by its action on the tangent spaces at $p_1,\ldots,p_g$.
%Assume $C$ is smooth. Suppose $\varphi$ acts on the tangent space at $p_i$ as rescaling by $\la_i$.
%Whenever $\a_{ij}\neq 0$, we have $\la_i^2=\la_j$. Whenever $\b_{ij}\neq 0$, we have $\la_i^3=\la_j$.
%Recall that $\a_{ij}=0$ means $h^0(2p_i+D_{ij})>1$.

\section{Further analysis of the algebra of functions on $\wt{\UU}^{ns,a}_{g,g}$}\label{equations-sec}

\subsection{Canonical parameters near marked points and expansions}

Let $(C,p_1,\ldots,p_g)$ be a family of non-special curves in $\wt{\UU}^{ns,a}_{g,g}(R)$, where $R$ is a
commutative ring,
and let $v_i\in H^0\bigl(\OO_C(p_i)/\OO_C\bigr)$ be the corresponding trivializations of normal bundles to sections
$p_i$. For example, we could take the universal family over $\Z[1/6]$ (see Theorem \ref{moduli-thm}).
By a parameter of order $N\ge 1$ at $p_i$, compatible with $v_i$, we understand an element 
$t_i\in H^0\bigl(C,\OO_C(-p_i)/\OO_C(-(N+1)p_i)\bigr)$ such that $\lan v_i, t_i\mod \OO_C(-2p_i)\ran=1$.
Note that since $H^1\bigl(C,\OO_C(-2p_i)/\OO_C(-(N+1)p_i)\bigr)=0$, we have an exact sequence
\begin{align*}
&0\to H^0\bigl(\OO_C(-2p_i)/\OO_C(-(N+1)p_i)\bigr)\to
H^0\bigl(C,\OO_C(-p_i)/\OO_C(-(N+1)p_i)\bigr)\to \\
&H^0\bigl(C,\OO_C(-p_i)/\OO_C(-2p_i)\bigr)\to 0
\end{align*}
which implies the existence of parameters of any order $N\ge 1$ at $p_i$.
For any such parameter the induced map
$$R[t_i]/(t_i^{N+1})\to H^0\bigl(C,\OO_C/\OO_C(-(N+1)p_i)\bigr)$$
is an isomorphism. We also have the induced isomorphism
$$t_i^{-N}R[t_i]/t_iR[t_i]\rTo{\sim} H^0\bigl(C,\OO_C(Np_i)/\OO_C(-p_i)\bigr),$$
so for a section $f\in H^0\bigl(C,\OO_C(Np_i)/\OO_C(-p_i)\bigr)$ we can define a polar part at $p_i$
as the corresponding element of $t_i^{-N}R[t_i]/R[t_i]$.

Set $D=p_1+\ldots+p_g\sub C$, $D_i=D-p_i$. We have the following relative version of
\cite[Lem.\ 4.1.3]{FP}.

\begin{lem}\label{parameter-lem} 
Assume that $N>1$ is such that $N!$ is invertible in $R$. 
Then there exist unique parameters $t_i$ of order $N$ at $p_i$, for $i=1,\ldots,g$, compatible
with $v_i$ and such that
for every $n$, $1<n\le N$, there exists a section
$$f_i[n]\in H^0(C,\OO(D_i+np_i)) \ \text{ with the polar part } t_i^{-n} \text{ at } p_i.$$
\end{lem}

\Pf . The argument is essentially the same as in \cite[Lem.\ 4.1.3]{FP}.
We start by picking arbitrary parameters $t_i$ at $p_i$ and then improve them step by step.
Since $H^1\bigl(C,\OO(D)\bigr)=0$, we have an exact sequence 
$$0\to R=H^0\bigl(C,\OO(D)\bigr)\to H^0\bigl(C,\OO(D_i+2p_i)\bigr)\to H^0\bigl(C,\OO(2p_i)/\OO(p_i)\bigr)\to 0.$$ 
Thus, we can choose $f_i[2]\in H^0\bigl(C,\OO(D_i+2p_i)\bigr)$ with the polar part at $p_i$ of the form
$\frac{1}{t_i^2}+\frac{c}{t_i}$ for some $c\in R$. Thus, replacing $t_i$ by $t_i+\frac{c}{2}t_i^2$
we obtain that the polar part of $f_i[2]$ becomes just $1/t_i^2$. Then we proceed by induction
as in \cite[Lem.\ 4.1.3]{FP}.
\ed

Assume that $(n+2)!$ is invertible in $R$. Then using the parameters of order $n+2$ constructed in the
above lemma we get partial Laurent expansions of sections of $\OO_C(np_i)/\OO_C(-3p_i)$ in
$t_i^{-n}R[t_i]/t_i^3R[t_i]$. In particular, for each $i$ the section $f_i[n]\in H^0(C,\OO(D_i+np_i))$ defined in Lemma
\ref{parameter-lem} has expansions
$$f_i[n]=\frac{1}{t_i^n}+c_{ii}[n]+l_{ii}[n]t_i+q_{ii}[n]t_i^2 \mod t_i^3,$$
$$f_i[n]=\frac{p_{ij}[n]}{t_j}+c_{ij}[n]+l_{ij}[n]t_j+q_{ij}[n]t_j^2 \mod t_j^3$$
at $p_i$ and $p_j$ (where $j\neq i$) respectively, for some constants
$p_{ij}[n]$, $c_{ij}[n]$, $l_{ij}[n]$ and $q_{ij}[n]$ in $R$.
Note that $f_i[n]$ are determined uniquely up to adding a constant, hence, all of these constants
except for $c_{ij}[n]$ are determined uniquely (and $c_{ij}[n]$ can be modified by
$c_{ij}[n]\mapsto c_{ij}[n]+c_i$ for some collection of constants $(c_i)$). Note that if we only need expansions
up to linear terms in $t_i$ then it is enough to consider the canonical parameters of order $n+1$,
defined whenever $(n+1)!$ is invertible.
%Thus, we can assume that
 %$$c_{ii}[n]=0 \text{ for } n\ge 4$$
%(we want to keep the flexibility of adding a constant to $f_i[2]$ and $f_i[3]$).

Now assume that $\cha(k)\neq 2$ or $3$, and let us apply the construction of Lemma \ref{parameter-lem} for $N=4$,
which gives canonical parameters $t_i$ of order $4$ and the 
corresponding sections $f_i[2]$, $f_i[3]$ and $f_i[4]$. Let us fix a choice of such $f_i[2]$, $f_i[3]$, $f_i[4]$
(there is a freedom in adding a constant to each).
Let us set $f_i=f_i[2]$, $h_i=f_i[3]$.  
We will see below that
$f_i$ and $h_i$ can always be adjusted by adding constants so as to satisfy conditions
of Lemma \ref{marked-alg-main-lem}(i)
for the marked algebra $A=H^0(C-D,\OO)$.
We are going to express the structure constants in the defining relations \eqref{main-eq-ff}--\eqref{main-eq-h2f3} of the marked algebra $A$ in terms of some of the constants 
$(p_{ij}[n], c_{ij}[n], l_{ij}[n], q_{ij}[n])$ (see Proposition \ref{constants-prop} below).
 
\subsection{Some relations}\label{relations-sec}

As we have observed above, the expansions of $f_i=f_i[2]$ up to quadratic terms 
(resp., of $h_i=f_i[3]$ up to linear terms, and of $f_i[4]$ up to constant terms) 
at all the marked points are well defined, so the following constants are well defined.
\begin{equation}
\label{some-constants-eq}
\begin{array}{l}
\a_{ij}:=p_{ij}[2], \\ 
\ga_{ij}:=c_{ij}[2], \ \b_{ij}:=p_{ij}[3],\\ 
\de_{ij}:=l_{ij}[2], \ \vareps_{ij}:=c_{ij}[3], \ \eta_{ij}:=p_{ij}[4], \\
 \pi_{ij}:=q_{ij}[2], \ \vartheta_{ij}:=l_{ij}[3], \ \zeta_{ij}:=c_{ij}[4].
\end{array}
\end{equation}
%where $\zeta_{ii}=0$. 
%As we will see later, these constants determine the products $m_3$, $m_4$, $m_5$
%and $m_6$ of the $A_\infty$-structure on $E_g$
%associated with $(C,p_1,\ldots,p_g)$.

For $i\neq j$ we have $f_if_j\in H^0\bigl(C,\OO(2D+p_i+p_j)\bigr)$, and 
analyzing the polar parts at $p_i$ and $p_j$ (we only need to look at the terms $\frac{1}{t_i^n}$
and $\frac{1}{t_j^n}$ for $n\ge 2$) we get the following more precise version of \eqref{main-eq-ff}
(see \cite[(4.1.4)]{FP}):
\begin{equation}\label{ff-prod-eq}
f_if_j=\sum_{k\neq i,j}\a_{ik}\a_{jk}f_k+\a_{ji}h_i+\a_{ij}h_j+\ga_{ji}f_i+\ga_{ij}f_j+a_{ij},
\end{equation}
for some constants $a_{ij}=a_{ji}$.  
By comparing the polar parts of both sides at $p_j$ and $p_k$ with $k\neq i,j$ we
get (see \cite[(4.1.5), (4.1.6)]{FP})
\begin{equation}\label{bc-rel}
\de_{ij}=\sum_{k\neq i,j}\a_{ik}\a_{jk}\a_{kj}+\a_{ji}\b_{ij}+(\ga_{ji}-\ga_{jj})\a_{ij}, 
\end{equation}
\begin{equation}\label{ct-rel}
\a_{ik}(\ga_{jk}-\ga_{ji})+\a_{jk}(\ga_{ik}-\ga_{ij})=\sum_{l\neq i,j,k}\a_{il}\a_{jl}\a_{lk}+\a_{ji}\b_{ik}+\a_{ij}\b_{jk}.
\end{equation}
%Note that \eqref{bc-rel} expresses $(\de_{ij})$ for $i\neq j$ in terms of $(\a_*,\b_*,\ga_*)$.
By comparing the constant terms in \eqref{ff-prod-eq} 
at $p_k$ with $k\neq i,j$ we get (see \cite[(4.1.7)]{FP})
\begin{equation}\label{abcd-rel}
%\begin{array}{l}
\a_{ik}\de_{jk}+\a_{jk}\de_{ik}+\ga_{ik}\ga_{jk}=
\sum_{l\neq i,j}\a_{il}\a_{jl}\ga_{lk}+\a_{ji}\varepsilon_{ik}+\a_{ij}\varepsilon_{jk}+\ga_{ji}\ga_{ik}+\ga_{ij}\ga_{jk}+a_{ij}.
%\end{array}
\end{equation}
In addition, comparing the constant terms in \eqref{ff-prod-eq} at $p_j$ we get
\begin{equation}\label{pi-a-rel}
\a_{ij}\de_{jj}+\pi_{ij}=\sum_{k\neq i,j}\a_{ik}\a_{jk}\ga_{kj}+\a_{ji}\vareps_{ij}+\a_{ij}\vareps_{jj}+
\ga_{ji}\ga_{ij}+a_{ij}.
\end{equation}
Finally, comparing the linear terms in \eqref{ff-prod-eq} at $p_k$, where $k\neq i,j$, and at $p_j$ we get
\begin{equation}\label{agdpt-rel}
\pi_{ik}\a_{jk}+\pi_{jk}\a_{ik}+(\ga_{ik}-\ga_{ij})\de_{jk}+(\ga_{jk}-\ga_{ji})\de_{ik}=\sum_{l\neq i,j}\a_{il}\a_{jl}\de_{lk}+\a_{ji}\vartheta_{ik}+\a_{ij}\vartheta_{jk},
\end{equation}
\begin{equation}\label{rho-eq}
\rho_{ij}+\de_{ij}\ga_{jj}+\a_{ij}\pi_{jj}=\sum_{k\neq i,j}\a_{ik}\a_{jk}\de_{kj}+\a_{ji}\vartheta_{ij}+
\a_{ij}\vartheta_{jj}+\ga_{ji}\de_{ij},
\end{equation}
where $\rho_{ij}$ is the coefficient of $t_j^3$ in the expansion of $f_i$ at $p_j$.

Similarly, for $i\neq j$ we have $f_ih_j-\a_{ij}f_j[4]\in H^0\bigl(C,\OO(2D+p_i+p_j)\bigr)$, so analyzing the polar 
parts we deduce the identity
%(see \cite[(4.1.8)]{FP}):
\begin{equation}\label{f-i-h-j-eq}
f_ih_j=\a_{ij}f_j[4]+\ga_{ij}h_j+\de_{ij}f_j+\b_{ji}h_i+\vareps_{ji}f_i+\sum_{k\neq i,j} \a_{ik}\b_{jk}f_k+b_{ij}
\end{equation}
for some constants $b_{ij}$. By comparing the polar parts of both sides at $p_i$ and $p_k$, where $k\neq i,j$,
we get 
%(see \cite[(4.1.9), (4.1.10)]{FP})
\begin{equation}\label{pr-rel}
\vartheta_{ji}=\sum_{k\neq i,j}\a_{ik}\b_{jk}\a_{ki}+(\ga_{ij}-\ga_{ii})\b_{ji}+\de_{ij}\a_{ji}
+\a_{ij}\eta_{ji},
\end{equation}
\begin{equation}\label{r-rel}
\a_{ik}(\varepsilon_{jk}-\varepsilon_{ji})+(\ga_{ik}-\ga_{ij})\b_{jk}=\sum_{l\neq i,j,k}\a_{il}\b_{jl}\a_{lk}+\b_{ji}\b_{ik}+
\de_{ij}\a_{jk}+\a_{ij}\eta_{jk}.
\end{equation}
%Note that \eqref{pr-rel} expresses $\vartheta_{ij}$ in terms of $(\a_*,\b_*,\ga_*)$ (where we use
%\eqref{bc-rel} to get rid of $\de_{ij}$).
Comparing the constant terms in \eqref{f-i-h-j-eq} at $p_k$, where $k\neq i,j$, we get
%(see \cite[(4.1.11)]{FP})
\begin{equation}\label{pqe-rel}
\a_{ik}\vartheta_{jk}+\ga_{ik}\varepsilon_{jk}+\de_{ik}\b_{jk}=
\sum_{l\neq i,j}\a_{il}\b_{jl}\ga_{lk}+\b_{ji}\varepsilon_{ik}+\ga_{ij}\varepsilon_{jk}+\varepsilon_{ji}\ga_{ik}+\de_{ij}\ga_{jk}+\a_{ij}\zeta_{jk}+b_{ij}.
\end{equation}
Comparing the polar and constant terms in \eqref{f-i-h-j-eq} at $p_j$ we get
\begin{equation}\label{q-ij-2-eq}
\pi_{ij}=\b_{ji}\b_{ij}+(\vareps_{ji}-\vareps_{jj})\a_{ij}+\sum_{k\neq i,j}\a_{ik}\a_{kj}\b_{jk},
\end{equation}
\begin{equation}\label{b-ij-main-eq}
\a_{ij}\vartheta_{jj}+\rho_{ij}=\a_{ij}\zeta_{jj}+\de_{ij}\ga_{jj}+\b_{ji}\vareps_{ij}+\ga_{ij}\vareps_{ji}
+\sum_{k\neq i,j}\a_{ik}\b_{jk}\ga_{kj}+b_{ij}.
\end{equation}
Let us denote by $\nu_{ji}$ the coefficient of $t_i^3$ in the expansion of $h_j$ at $p_i$.
Looking at the linear terms in \eqref{f-i-h-j-eq} at $p_i$ we get
\begin{equation}\label{nu-l[4]-eq}
\nu_{ji}+\ga_{ii}\vartheta_{ji}+\pi_{ii}\b_{ji}=\a_{ij}l_{ji}[4]+\ga_{ij}\vartheta_{ji}+\de_{ij}\de_{ji}+
\b_{ji}\vartheta_{ii}+\sum_{k\neq i,j}\a_{ik}\b_{jk}\de_{ki}.
\end{equation}

Let us fix the unique choice of $f_i[4]$ (by adding a constant) so that we have
\begin{equation}\label{fi-2-eq}
f_i^2=f_i[4]+\sum_{j\neq i}\a_{ij}^2f_j+2\ga_{ii}f_i.
\end{equation}
From this, looking at the polar and constant parts at $p_i$ and $p_j$ we get
\begin{equation}\label{l-ii-eq}
2\de_{ii}=\sum_{j\neq i}\a_{ij}^2\a_{ji},
\end{equation}
\begin{equation}\label{q-ii-u-eq}
2\pi_{ii}=\sum_{j\neq i}\a_{ij}^2\ga_{ji}+\ga_{ii}^2+\zeta_{ii},
\end{equation}
\begin{equation}\label{p-ij-4-eq}
\eta_{ij}=2\a_{ij}(\ga_{ij}-\ga_{ii})-\sum_{k\neq i,j}\a_{ik}^2\a_{kj},
\end{equation}
\begin{equation}\label{c-ij-4-eq}
\ga_{ij}^2+2\a_{ij}\de_{ij}=\zeta_{ij}+\sum_{k\neq i}\a_{ik}^2\ga_{kj}+
2\ga_{ii}\ga_{ij}.
\end{equation}
In addition, looking at the linear parts at $p_j$ we get
\begin{equation}\label{l[4]-eq}
2\a_{ij}\pi_{ij}+2\ga_{ij}\de_{ij}=
l_{ij}[4]+\sum_{k\neq i} \a_{ik}^2\de_{kj}+2\ga_{ii}\de_{ij}.
\end{equation}
%Using \eqref{q-ii-u-eq} we express $u_i$ in terms of $(\a_{ij})$, $(c_{ij}[2])$, $(c_{ii}[2])$ and $(q_{ii}[2])$.
%In addition, we get equations ???
%Finally, using \eqref{q-ii-u-eq} and \eqref{c-ij-4-eq}, we express $c_{ij}[4]$ in terms of 
%$(\a_{ij})$, $(c_{ij}[2])$, $(c_{ii}[2])$, $(l_{ij}[2])$ and $(q_{ii}[2])$.

Next, comparing the polar parts of $h_i^2$ and $f_i^3$ we get the relation
\begin{equation}\label{h-i-2-f-i-3-eq}
\begin{array}{l}
h_i^2=f_i^3-3\ga_{ii}f_i[4]+(2\vareps_{ii}-3\de_{ii})h_i+(2\vartheta_{ii}-3\pi_{ii}-3\ga_{ii}^2)f_i\\
-\sum_{j\neq i}\a_{ij}^3h_j+\sum_{j\neq i}(\b_{ij}^2-3\a_{ij}^2\ga_{ij})f_j+s_i
\end{array}
\end{equation}
for some constant $s_i$.
Considering the polar and constant terms at $p_j$ we get
\begin{equation}\label{pi-theta-eq}
\begin{array}{l}
2\b_{ij}\vareps_{ij}=3\a_{ij}^2\de_{ij}+3\a_{ij}\ga_{ij}^2
-3\ga_{ii}\eta_{ij}+(2\vareps_{ii}-3\de_{ii})\b_{ij}+
(2\vartheta_{ii}-3\pi_{ii}-3\ga_{ii}^2)\a_{ij}-\\
\sum_{k\neq i,j}\a_{ik}^3\b_{kj}+
\sum_{k\neq i,j}(\b_{ik}^2-3\a_{ik}^2\ga_{ik})\a_{kj},
\end{array}
\end{equation}
\begin{equation}\label{s-i-eq}
\begin{array}{l}
\vareps_{ij}^2+2\b_{ij}\vartheta_{ij}=3\a_{ij}^2\pi_{ij}+6\a_{ij}\ga_{ij}\eta_{ij}+\ga_{ij}^3-
3\ga_{ii}\zeta_{ij}+(2\vareps_{ii}-3\de_{ii})\vareps_{ij}\\
+(2\vartheta_{ii}-3\pi_{ii}-3\ga_{ii}^2)\ga_{ij}
-\sum_{k\neq i}\a_{ik}^3\vareps_{kj}+\sum_{k\neq i}(\b_{ik}^2-3\a_{ik}^2\ga_{ik})\ga_{kj}+s_i.
\end{array}
\end{equation}

Next, considering the product $h_ih_j$ for $i\neq j$ we get
\begin{equation}\label{h-i-h-j-eq}
h_ih_j=\b_{ij}f_j[4]+\b_{ji}f_i[4]+\vareps_{ij}h_j+\vareps_{ji}h_i+\vartheta_{ij}f_j+\vartheta_{ji}f_i+
\sum_{k\neq i,j}\b_{ik}\b_{jk}f_k+u_{ij}
\end{equation}
for some constants $u_{ij}=u_{ji}$.
Looking at the polar and constant terms at $p_k$, where $k\neq i,j$, we get
\begin{equation}\label{ab-eps-th-eta-eq}
\b_{ik}(\vareps_{jk}-\vareps_{ji})+\b_{jk}(\vareps_{ik}-\vareps_{ij})=
\b_{ij}\eta_{jk}+\b_{ji}\eta_{ik}+\vartheta_{ij}\a_{jk}+\vartheta_{ji}\a_{ik}+
\sum_{l\neq i,j,k}\b_{il}\b_{jl}\a_{lk},
\end{equation}
\begin{equation}\label{u-i-j-eq}
\begin{array}{l}
\b_{ik}\vartheta_{jk}+\b_{jk}\vartheta_{ik}+\vareps_{ik}\vareps_{jk}=
\b_{ij}\zeta_{jk}+\b_{ji}\zeta_{ik}+\vareps_{ij}\vareps_{jk}+\vareps_{ji}\vareps_{ik}+
\vartheta_{ij}\ga_{jk}+\vartheta_{ji}\ga_{ik}+\\
\sum_{l\neq i,j}\b_{il}\b_{jl}\ga_{lk}+u_{ij}.
\end{array}
\end{equation}
Also, comparing the constant terms at $p_j$ in \eqref{h-i-h-j-eq} we get
\begin{equation}\label{nu-u-eq}
\nu_{ij}+\vareps_{ij}\vareps_{jj}+\b_{ij}\vartheta_{jj}=\b_{ij}\zeta_{jj}+\b_{ji}\zeta_{ij}+
\vareps_{ij}\vareps_{jj}+\vareps_{ji}\vareps_{ij}+\vartheta_{ij}\ga_{jj}+\vartheta_{ji}\ga_{ij}+
\sum_{k\neq i,j}\b_{ik}\b_{jk}\ga_{kj}+u_{ij}
\end{equation}
(recall that $\nu_{ij}$ is the coefficient of $t_j^3$ in the expansion of $h_i$ at $p_j$).

%Finally, assume in addition that $\cha(k)>5$. Then we can apply the construction of Lemma 
%\ref{parameter-lem} for $N=5$ extending $t_i$'s to parameters of order $5$
%and obtaining the corresponding sections $f_i[5]\in H^0\bigl(C,\OO(D_i+5p_i)\bigr)$. 
%Looking at the products $f_ih_i\in H^0\bigl(C,\OO(5D)\bigr)$ and adding an appropriate constant to $f_i[5]$, we
%get
%\begin{equation}\label{f-i-h-i-eq}
%f_ih_i=f_i[5]+\ga_{ii}h_i+(\de_{ii}+\vareps_{ii})f_i+\sum_{j\neq i}\a_{ij}\b_{ij}f_j.
%\end{equation}
%Comparing the polar parts at $p_i$ and $p_j$ we get
%\begin{equation}\label{l-ii-3-eq}
%\pi_{ii}+\vartheta_{ii}=\sum_{j\neq i}\a_{ij}\a_{ji}\b_{ij}.
%\end{equation}

\subsection{Minimal set of generators}\label{min-gen-sec}

%Let us make a few more observations about the relations obtained above.

We have seen that the constants (elements of $R$) introduced in the previous section, such as
the constants \eqref{some-constants-eq}, as well as
$a_{ij}$, $b_{ij}$, $u_{ij}$ and $s_i$, satisfy many polynomial relations with integer coefficients.
We are going to use these relations to find a minimal set of generators of the algebra of functions on 
$\wt{\UU}^{ns,a}_{g,g}$.

First, we will relate the constants \eqref{some-constants-eq}
%calculated for the universal family over $\wt{\UU}^{ns,a}_{g,g}$, 
to the coefficients \eqref{big-generators-eq} that determine the equations of the affine curve $C\setminus D$ over $R$.

\begin{lem}\label{two-coef-lem} 
Assume that $6$ is invertible in $R$. 

\noindent
(i) Let us normalize (by adding constants)
$f_i=f_i[2]$ and $h_i=f_i[3]$ so that $\ga_{ii}=0$ and $2\vareps_{ii}=3\de_{ii}$. Then 
the elements $(f_i,h_i)_{i=1,\ldots,g}$ are the canonical generators of the algebra $\OO(C\setminus D)$ 
considered in Lemma \ref{grobner-lem}. One has the following relations between the constants 
\eqref{some-constants-eq}, defined using canonical parameters of order $4$ at $p_i$ as in Lemma \ref{parameter-lem}, 
and the coefficients \eqref{big-generators-eq}:
\begin{equation}\label{pi-i-th-pi-ii-eq}
\pi_i=2\vartheta_{ii}-3\pi_{ii}.
\end{equation}
$$d_{ij}=\a_{ij}, \ t_{ji}=\b_{ji}, \ v_{ij}=\ga_{ij}, \ c_{ij}^k=\a_{ik}\a_{jk}, \ g_i^j=-\a_{ij}^3,$$
$$r_{ji}=\vareps_{ji}-\a_{ij}\a_{ji}^2,$$
$$\de_{ij}=\a_{ji}\b_{ij}+\a_{ij}\ga_{ji}+\sum_{k\neq i,j}\a_{ik}\a_{jk}\a_{kj},$$
$$e_{ij}^k=\a_{ik}\b_{jk}-\a_{ij}\a_{jk}^2,$$
$$\psi_{ij}=\ga_{ji}\b_{ij}+3\a_{ij}\a_{ji}\ga_{ij}
+\sum_{k\neq i,j}\b_{ik}\a_{jk}\a_{kj}+\sum_{k\neq i,j}\a_{ik}(\a_{ij}\a_{jk}\a_{ki}
-\a_{ji}\a_{kj}\a_{ik}),$$
$$l_{ij}^k=\b_{ik}\b_{jk}-\b_{ij}\a_{jk}^2-\b_{ji}\a_{ik}^2,$$
$$k_i^j=\b_{ij}^2-3\a_{ij}^2\ga_{ij}.$$

\noindent
(ii) 
%\begin{lem}\label{express-lem} 
Assume that $g\ge 2$ and that $f_i$ and $h_i$ are normalized as in (i).
Let us denote by $R[C]\sub R$ the $\Z[1/6]$-subalgebra generated by $\a_{ij}$,
$\b_{ij}$, $\ga_{ij}$, $\vareps_{ij}$, and $\pi_i$ (where $1\le i,j\le g$, $i\neq j$).
Let us also denote by $R[C]_1\sub R[C]$ the similar subalgebra where we do not allow to use $\pi_i$, and by 
$R[C]_0\sub R[C]_1$ the subalgebra generated by $\a_{ij}$, $\b_{ij}$ and $\ga_{ij}$.
Then we have 
$$\de_{ij}, \de_{ii}, \eta_{ij}, \vartheta_{ij}, \zeta_{ij}\in R[C]_0,$$ 
$$\pi_{ij}, a_{ij}\in R[C]_1,$$
$$b_{ij}, s_i, u_{ij}\in R[C],$$
where $i\neq j$.
If $g\ge 3$ then $u_{ij}\in R[C]_1$.
\end{lem}

\Pf . (i) With the given normalization of $f_i$ and $h_i$, the equations \eqref{ff-prod-eq},
\eqref{f-i-h-j-eq}, \eqref{h-i-2-f-i-3-eq}, \eqref{h-i-h-j-eq}, together with the expression \eqref{fi-2-eq}
for $f_i[4]$, give the following equations:
\begin{equation}\label{ff-prod-eq2}
f_if_j=\a_{ji}h_i+\a_{ij}h_j+\ga_{ji}f_i+\ga_{ij}f_j+\sum_{k\neq i,j}\a_{ik}\a_{jk}f_k+a_{ij},
\end{equation}
\begin{equation}\label{f-i-h-j-eq2}
f_ih_j=\a_{ij}f_j^2
+\ga_{ij}h_j+\b_{ji}h_i+\de_{ij}f_j
+(\vareps_{ji}-\a_{ij}\a_{ji}^2)f_i+\sum_{k\neq i,j}(\a_{ik}\b_{jk}-\a_{ij}\a_{jk}^2)f_k+b_{ij},
\end{equation}
\begin{equation}\label{h-i-2-f-i-3-eq2}
h_i^2=f_i^3+(2\vartheta_{ii}-3\pi_{ii})f_i
-\sum_{j\neq i}\a_{ij}^3h_j+\sum_{j\neq i}(\b_{ij}^2-3\a_{ij}^2\ga_{ij})f_j+s_i,
\end{equation}
\begin{equation}\label{h-i-h-j-eq2}
\begin{array}{l}
h_ih_j=\b_{ij}f_j^2+\b_{ji}f_i^2+\vareps_{ij}h_j+\vareps_{ji}h_i+
(\vartheta_{ij}-\b_{ji}\a_{ij}^2)f_j+(\vartheta_{ji}-\b_{ij}\a_{ji}^2)f_i+\\
\sum_{k\neq i,j}(\b_{ik}\b_{jk}-\b_{ij}\a_{jk}^2-\b_{ji}\a_{ik}^2)f_k+u_{ij}.
\end{array}
\end{equation}
Note that these equations have the same form as in Lemma \ref{marked-alg-main-lem}(iii).
Hence, by Lemma \ref{marked-alg-main-lem}(iv),
the elements $(f_i,h_i)$ are the canonical generators of the marked algebra
$\OO(C\setminus D)$.
Thus, all the coefficients in these equations should match.
From this we get most of the required equations. For $\psi_{ij}$ we get
$$\psi_{ij}=\vartheta_{ij}-\b_{ji}\a_{ij}^2$$
which we have to transform further using \eqref{pr-rel} and using the formula \eqref{p-ij-4-eq} for $\eta_{ij}$.

\noindent (ii)
First, relations \eqref{bc-rel}, \eqref{l-ii-eq} and \eqref{p-ij-4-eq} give the inclusions
$$\de_{ij}, \de_{ii}, \eta_{ij}\in R[C]_0.$$
It follows that $\vareps_{ii}=3\de_{ii}/2$ is also in $R[C]_0$.
Now relations \eqref{q-ij-2-eq}, \eqref{pr-rel} and \eqref{c-ij-4-eq} give
$$\pi_{ij}\in R[C]_1 \ \text {and } \vartheta_{ij},\zeta_{ij}\in R[C]_0.$$
Thus, \eqref{pi-a-rel} and \eqref{s-i-eq} (as well as \eqref{pi-i-th-pi-ii-eq}) give
$$a_{ij}\in R[C]_1 \ \text{ and } s_i\in R[C],$$

Next, we want to show that $b_{ij}$ is in $R[C]$. 
First, we note that \eqref{q-ii-u-eq} implies that
$$\zeta_{ii}-2\pi_{ii}\in R[C]_0,$$
while \eqref{rho-eq} gives
$$\rho_{ij}+\a_{ij}(\pi_{jj}-\vartheta_{jj})\in R[C]_0.$$
Hence, we have the following congruence modulo $R[C]_0$:
$$\rho_{ij}+\a_{ij}(\vartheta_{jj}-\zeta_{jj})\equiv \a_{ij}(2\vartheta_{jj}-\pi_{jj}-\zeta_{jj})\equiv \a_{ij}(2\vartheta_{jj}-3\pi_{jj})=
\a_{ij}\pi_j
\mod R[C]_0,$$
and so, $\rho_{ij}+\a_{ij}(\vartheta_{jj}-\zeta_{jj})$ is in $R[C]$. Now \eqref{b-ij-main-eq} implies that 
$b_{ij}$ is in $R[C]$.

It remains to express $u_{ij}$. In the case $g\ge 3$ using \eqref{u-i-j-eq} we get $u_{ij}\in R[C]_1$.
In the case $g=2$ we have to take a longer route. First, \eqref{l[4]-eq} shows that
$l_{ij}[4]\in R[C]_1$. Next, \eqref{nu-l[4]-eq} shows 
$\nu_{ij}+\b_{ij}(\pi_{jj}-\vartheta_{jj})\in R[C]$, and finally,
\eqref{nu-u-eq} gives
$$u_{ij}\equiv \nu_{ij}+\b_{ij}(\vartheta_{jj}-\zeta_{jj})\equiv \b_{ij}\pi_j \mod R[C]_1,$$
and so, $u_{ij}$ is in $R[C]$.
\ed

Recall that we have constructed an isomorphism of the moduli space of non-special curves (with tangent vectors
at the marked points)
$\wt{\UU}^{ns,a}_{g,g}\times\Spec(\Z[1/6])$ with the affine scheme $\SS_g$ over $\Z[1/6]$
parametrizing algebras with defining relations \eqref{main-eq-ff}--\eqref{main-eq-h2f3}
(see Theorem \ref{moduli-thm}).

\begin{prop}\label{constants-prop} Assume $g\ge 2$.

\noindent
(i) The algebra of functions on the affine scheme
$\wt{\UU}^{ns,a}_{g,g}\times \Spec(\Z[1/6])\simeq \SS_g$  
is generated over $\Z[1/6]$ by the elements
\begin{equation}\label{main-constants-eq}
\a_{ij}, \b_{ij}, \ga_{ij}, \vareps_{ij}, \pi_i
\end{equation}
(where $1\le i,j\le g$, $i\neq j$), defined for the universal curve.

\noindent
(ii) Let $k$ be a field of characteristic $\neq 2,3$. Let $[C^{\cusp}_g]\in \SS_g(k)$ 
be the point at which all the coordinates \eqref{big-generators-eq} vanish. Let
$\mg$ be the corresponding maximal ideal in the ring of functions on $\SS_g\times\Spec(k)$.
Then the elements \eqref{main-constants-eq} project to a basis in
$\mg/\mg^2$.
\end{prop}

\Pf . (i) Note that this algebra is generated by the coefficients
\eqref{big-generators-eq} defining equations of the universal affine curve $C\setminus D$. 
Hence, the assertion follows from Lemma \ref{two-coef-lem}. 

\noindent
(ii) We have to calculate the fiber of $\SS_g\bigl(k[u]/(u^2)\bigr)\to\SS_g(k)$ over $[C^{\cusp}_g]\in\SS_g(k)$.
In other words, we have to look at the condition of associativity for the relations 
\eqref{main-eq-ff}--\eqref{main-eq-h2f3}, where all
the coordinates \eqref{big-generators-eq} belong to the ideal $(u)\sub k[u]/(u^2)$.
As we have seen in Lemma \ref{two-coef-lem}, these relations have form
$$f_if_j=\a_{ji}h_i+\a_{ij}h_j+\ga_{ji}f_i+\ga_{ij}f_j,$$
$$f_ih_j=\a_{ij}f_j^2+\ga_{ij}h_j+\b_{ji}h_i+\vareps_{ji}f_i,$$
$$h_i^2=f_i^3+\pi_if_i,$$
$$h_ih_j=\b_{ij}f_j^2+\b_{ji}f_i^2+\vareps_{ij}h_j+\vareps_{ji}h_i.$$
It is easy to check that in fact these relations define an associative algebra over $k[u]/(u^2)$
with the basis $(f_i^n, h_if_i^n)$
for arbitrary $\a_{ij},\b_{ij},\ga_{ij},\vareps_{ij},\pi_i\in (u)\sub k[u]/(u^2)$,
which implies the result.
\ed

The above Proposition shows that \eqref{main-constants-eq} are minimal generators of the ring
of functions on $\SS_g$. 

We have the following geometric interpretation of the vanishing of the functions $\a_{ij}$ and $\b_{ij}$ on
$\SS_g$.

\begin{lem}\label{constants-vanishing-lem} Assume that $g\ge 2$.
Let $(C,p_1,\ldots,p_g)\in \wt{\UU}^{ns,a}_{g,g}(k)$,
where $k$ is a field of characteristic $\neq 2,3$. Let us consider the functions $\a_{ij},\b_{ij}$
on $\wt{\UU}^{ns,a}_{g,g}(k)\simeq \SS_g(k)$.

\noindent
(i) One has $\a_{ij}(C,p_1,\ldots,p_g)=0$ if and only if $h^0(2p_i+D_i-p_j)=2$, where $D_i=D-p_i$.

\noindent
(ii) One has $\a_{ij}(C,p_1,\ldots,p_g)=\b_{ij}(C,p_1,\ldots,p_g)=0$ if and only if 
$h^0(3p_i+D_i-p_j)=3$.

\noindent
(iii) If $C$ is smooth then for every $i$ either there exists $j\neq i$ such that $\a_{ij}\neq 0$ or there exists
$j\neq i$ such that $\b_{ij}\neq 0$. 
\end{lem}

\Pf . (i) We use the interpretation of $\a_{ij}$ as the coefficient of $t_j^{-1}$ in the Laurent series of
$f_i$ at $p_j$ (see Section \ref{relations-sec}). 
Thus, $\a_{ij}=0$ means that $f_i$ in fact is regular at $p_j$, i.e.,
$$H^0(C,\OO(2p_i+D_i))=H^0(C,\OO(2p_i+D_i-p_j))=0.$$
Since $h^0(2p_i+D_i)$=2, the assertion follows.

\noindent
(ii) We use in addition the interpretation of $\b_{ij}$ as the coefficient of $t_j^{-1}$ in the Laurent series
of $h_i$ at $p_j$. Thus, $\a_{ij}=\b_{ij}=0$ if and only if both $f_i$ and $h_i$ are regular at $p_j$.
Since $1,f_i,h_i$ is a basis of $H^0(C,\OO(3p_i+D_i))$, this is equivalent to having
$$H^0(C,\OO(3p_i+D_i))=H^0(C,\OO(3p_i+D_i-p_j)),$$
and the assertion follows.

\noindent
(iii) By part (i), if for some $i$ we have $\a_{ij}=0$ for all $j\neq i$ then $f_i$ is a nonconstant section
of $H^0(C,\OO(2p_i))$, so $C$ is hyperelliptic and $p_i$ is a Weierstrass point. Let $f:C\to \P^1$ be the double covering. Then $f_*\OO_C(p_i)=\OO_{\P^1}\oplus \OO_{\P^1}(-g)$, so
$$H^0\bigl(C,\OO(3p_i)\bigr)=H^0\bigl(\P^1,\bigl(f_*\OO(p_i)\bigr)(1)\bigr)\simeq H^0\bigl(\P^1,\OO(1)\oplus \OO(1-g)\bigr)$$
has dimension $2$. Hence, the result follows from part (ii).
\ed

\begin{rem} The interpretations of $\a_{ij}$ and $\b_{ij}$ from Section \ref{relations-sec}
imply that $\a_{ij}$ are exactly the functions
defined in \cite[Sec.\ 2.4]{FP} as triple Massey products, while the values of $\b_{ij}$, in the situation
when $\a_{ij}=0$ for all $j\neq i$, 
are given by the quadruple Massey products considered in \cite[Sec.\ 2.5]{FP}.
\end{rem}

\subsection{GIT quotients and weakly modular compactifications of $M_{g,g}$}\label{GIT-sec}

In this section we work over $\Z[1/6]$. In particular, we denote by $\A^n$ the affine space over $\Z[1/6]$.
%an algebraically closed field $k$ of characteristic zero,
%and we denote by $\wt{\UU}^{ns,a}_{g,g}$, $\SS_g$, etc., the schemes over $\Z[1/6]$ obtained
%from the schemes considered above by the base change. 
We use the version of GIT over general base schemes
as developed in \cite{Alper}.

We would like to consider GIT quotients with respect to the $\G_m^g$-action on the
affine scheme $\SS_g$ (see Theorem \ref{moduli-thm} for the definition of this action). 
Such a quotient is determined by a nontrivial character
$\chi:\G_m^g\to \G_m$ (as a linearized ample line bundle on $\SS_g$ we take the trivial line bundle
with the $\G_m^g$-action twisted by $\chi$) and is given by the scheme
$$\SS_g\sslash_\chi \G_m^g:=\Proj \left(\bigoplus_{n\ge 0} H^0(\SS_g,\OO)_{\chi^n}\right),$$
where the subscript $\chi^n$ denotes the subset of functions $f$ such that $(\la^{-1})^*f=\chi(\la)^nf$
for $\la\in\G_m^g$.

We have a closed embedding 
\begin{equation}\label{affine-emb-eq}
\SS_g\hra \A^{4g^2-3g},
\end{equation} 
where
$\A^{4g^2-3g}$ is the affine space with coordinates corresponding to the minimal generators
$\a_{ij}$, $\b_{ij}$, $\ga_{ij}$, $\vareps_{ij}$ and $\pi_i$ of $\OO(\SS_g)$ (see
Proposition \ref{constants-prop}).

Note that the minimal generators \eqref{main-constants-eq} 
are semihomogeneous for the $\G_m^g$-action on $\SS_g$. Namely, 
if we define the action of $\la\in\G_m^g$ on $\OO(\SS_g)$ by $f\mapsto (\la^{-1})^*f$
then the elements \eqref{main-constants-eq} transform under the action of
$\la=(\la_1,\ldots,\la_g)$ as follows:
$$\a_{ij}\mapsto \la_i^2\la_j^{-1}\a_{ij}, \ \ \b_{ij}\mapsto \la_i^3\la_j^{-1}\a_{ij}, \ \ 
\ga_{ij}\mapsto \la_i^2\ga_{ij}, \ \ \vareps_{ij}\mapsto \la_i^3\vareps_{ij}, \ \ \pi_i\mapsto \la_i^4\pi_i.$$
Thus, the embedding \eqref{affine-emb-eq} is $\G_m^g$-equivariant, and 
for every character $\chi$ we can check $\chi$-semistability (resp., stability) of a point in $\SS_g$ by
viewing it as a point of $\A^{4g^2-3g}$. In particular, we have a closed embedding of the corresponding GIT quotients
$$\SS_g\sslash_\chi \G_m^g \hra P_{g,\chi}:=\A^{4g^2-3g}\sslash_\chi \G_m^g.$$

Let us consider for a moment the general situation when an algebraic torus $T$ acts on an affine space $\A^N$
via a homomorphism $\rho:T\to \G_m^N$, so that the coordinates on $\A^N$ are semihomogeneous for the $T$-action
of weights $\om_1,\ldots,\om_N\in X(T)$, where $X(T)$ is the character lattice of $T$.
Assuming in addition that $\rho$ has finite kernel, let us recall a well known description of the walls in $X(T)\ot\R$, so that
the GIT quotients change only when crossing a wall (see e.g., \cite[Sec.\ 14.3, 14.4]{CLS}; for the reader's convenience
we will give a direct proof in Lemma \ref{GIT-walls-lem} below).

Let us set $\Om:=\{\om_1,\ldots,\om_N\}\sub X(T)$. The assumption that $\rho$ has finite kernel
is equivalent to the fact that the set $\Om$ spans $X(T)\ot\R$.
For every finite subset $S\sub X(T)$ we denote by
$\bC_S\sub X(T)\ot\R$ the closed cone generated by $S$, i.e., the set of linear combinations of elements of $S$ with
nonnegative real coefficients. In the case when $\bC_S$ is of full rank we denote by $\int(\bC_S)$ the interior of this
cone. Note that an element of $X(T)\ot\Q$ belongs to $\bC_S$ if and only if it is a linear combination of elements of $S$
with coefficients in $\Q_{\ge 0}$.

Let us define the subset $\Sigma\sub \bC_\Om$ by
\begin{equation}\label{Sigma-eq}
\Sigma=\cup_{\Om'\sub\Om: \dim \span(\Om')<\dim X(T)\ot\R} \bC_{\Om'}.
\end{equation}
Note that $\Sigma$ is the union of closed subcones of codimension $1$ in $\bC_\Om$ (the {\it walls}). 
The characters in $\bC_\Om\setminus\Sigma$ are called {\it generic}. The connected components of
$\bC_\Om\setminus\Sigma$ are called {\it chambers}. 
For a character $\chi\in X(T)$ let $(\A^N)^s_\chi$ and $(\A^N)^{ss}_\chi$ denote the subsets of $\chi$-stable
and $\chi$-semistable points in $\A^N$, respectively.

\begin{lem}\label{GIT-walls-lem}
(i) For $x\in \A^N$ let $\Om(x)\sub \Om$ be the set of $T$-weights of all the coordinates that do not vanish at $x$.
Then $x\in (\A^N)^{ss}_\chi$ if and only if $\chi\in \bC_{\Om(x)}$. In particular,
if $\chi\not\in\bC_\Om$ then $(\A^N)^{ss}_\chi=\emptyset$.

\noindent
(ii) If for a point $x\in \A^N$ the set $\Om(x)$ spans $X(T)\ot\Q$ over $\Q$, then for any $\chi\in \int(\bC_{\Om(x)})\cap X(T)$
the point $x$ is $\chi$-stable.

\noindent
(iii) For any character $\chi$ such that $\chi\not\in\Sigma$ one has $(\A^N)^s_\chi= (\A^N)^{ss}_\chi$.

\noindent
(iv) If $\chi$ and $\chi'$ belong to the same chamber of $\bC\setminus\Sigma$ then
$(\A^N)^{ss}_\chi=(\A^N)^{ss}_{\chi'}$, hence, the corresponding GIT quotients are the same.
\end{lem}

\Pf . (i) By definition, a point $x$ is $\chi$-semistable if and only if there exists a monomial $M$ in coordinates of a weight
proportional to $\chi$, such that $M$ does not vanish at $x$. Such a monomial exists if and only if $\chi\in \bC_{\Om(x)}$. 

\noindent
(ii) For $x\in (\A^N)^{ss}_\chi$ let $S\sub\{1,\ldots,N\}$ be the set of all the coordinates that do not vanish at $x$, 
so that $\Om(x)$ is the set of weights of coordinates in $S$. Consider the open subset
$$U_x=\{(y_1,\ldots,y_n)\in \A^N \ |\ \prod_{i\in S} y_i\neq 0\}.$$
Note that by (i), for every $y\in U_x$ we have $\chi\in\bC_{\Om(x)}\sub \bC_{\Om(y)}$, so all the points of
$U_x$ are $\chi$-semistable. On the other hand, since $\chi$ is an interior point of $\bC_{\Om(x)}$,
we can write $\chi$ as a linear combination of elements of $\Om(x)$ with positive rational coefficients.
Hence, there exists a monomial $M=\prod_{i\in S}x_i^{n_i}$ with $n_i>0$, of a weight proportional to $\chi$, so that
$U_x$ coincides with the open subset $M\neq 0$.
Furthermore, the assumption that $\Om(x)$ generates the sublattice of full rank
in $X(T)$ implies that all the points of $U_x$ have finite stabilizers. Hence, $T$ acts on $U_x$ with closed orbits, so
all the points of $U_x$ are $\chi$-stable.

\noindent 
(iii) Since $\chi$ is generic, any subset $S\sub \Om$, such that $\chi\in\bC_S$, spans $X(T)\ot\R$.
In particular, for any such $S$ one has $\chi\in\int(\bC_S)$. Now if $x$ is $\chi$-semistable then
$\chi\in\bC_{\Om(x)}$, so using (ii) we derive that $x$ is $\chi$-stable.

\noindent
(iv) Assume there exists a point $x\in (\A^N)^{ss}_\chi$ such that $x\not\in (\A^N)^{ss}_{\chi'}$.
Then we have $\chi\in \bC_{\Om(x)}$ and $\chi'\not\in\bC_{\Om(x)}$. But this implies that
the segment connecting $\chi$ and $\chi'$ in $X(T)\ot\R$ intersects the boundary of $\bC_{\Om(x)}$,
which is contained in $\Sigma$. This contradicts the assumption that this segment is contained entirely
in $\bC\setminus\Sigma$. Thus, the semistable loci for $\chi$ and $\chi'$ are the same. This implies that
the GIT quotients are  also the same, as the categorical quotients of the semistable loci.
\ed
 
Now let us return to our situation with the action of $T=\G_m^g$ on $\SS_g\sub\A^{4g^2-3g}$.
We identify the lattice of characters of $\G_m^g$ with the standard lattice $\Z^g\sub \R^g$, where
the basis vector $e_i$ corresponds to the projection $\G_m^g\to\G_m$ to the $i$th factor.
The set $\Om\sub\Z^g$ of $\G_m^g$-weights of the coordinates on $\A^{4g^2-3g}$ consists of
the vectors $2e_i-e_j$, $3e_i-e_j$, for $i\neq j$ and vectors $2e_i$, $3e_i$ and $4e_i$.
For $g\ge 2$ we denote by $\bC\sub\R^g$ the cone generated by the vectors $2e_i-e_j$, $i\neq j$. For $g=1$
we simply set $\bC=\R_{\ge 0}$.
Note that $\bC$ contains all the vectors $e_i=\frac{2}{3}(2e_i-e_j)+\frac{1}{3}(2e_j-e_i)$ and $3e_i-e_j$, 
hence, we have $\bC_\Om=\bC$.
As before, we denote by $\Sigma$ the union of the walls \eqref{Sigma-eq}.

\begin{prop}\label{GIT-prop}
(i) For every character $\chi\in\Z^g\cap\bC$
the quotient $P_{g,\chi}$ is a projective toric scheme over $\Z[1/6]$.
% of dimension $4(g^2-g)$.
Assume further that $\chi$ is in the interior of $\bC$.
Then for any smooth curve $C$ of genus $g$ and for generic points $p_1,\ldots,p_g\in C$,
with any nonzero tangent vectors $v_i$ at $p_i$'s, 
the point $(C,p_1,\ldots,p_g,v_1,\ldots,v_g)\in \wt{\UU}^{ns,a}_{g,g}\times\Spec(\Z[1/6])\simeq \SS_g$
is $\chi$-stable.

\noindent
(ii) Let $\bC_0\sub \bC$ be the subcone generated by all the vectors $e_i$. Assume that the character 
$\chi\in\Z^g$ is in $\int(\bC_0)$.
Then any point $(C,p_1,\ldots,p_g,v_1,\ldots,v_g)\in \SS_g$ with
$C$ smooth is $\chi$-stable.

\noindent
(iii) The interior of $\bC_0$ lies inside $\bC\setminus\Sigma$. Hence, for $\chi\in\int(\bC_0)$
the subset $(\SS_g)^{ss}_{\chi}=(\SS_g)^s_{\chi}$ does not depend on $\chi$.
\end{prop}
 
\Pf . (i) Since all  the vectors $2e_i-e_j$, $e_i$ and $3e_i-e_j$ have positive scalar product with 
$e_1+\ldots+e_g$, we see that any nontrivial combination of them with non-negative
coefficients is nonzero. Therefore, $H^0(\A^{4g^2-3g},\OO)^{\G_m^g}=\Z[1/6]$, which immediately implies
that $P_{g,\chi}$ is projective over $\Z[1/6]$ for any $\chi$. Hence, $\SS_g\sslash_\chi \G_m^g$ is also projective.
The fact that $P_{g,\chi}$ is a toric scheme is well known (see \cite[Thm.\ 14.2.13]{CLS} or \cite[Thm.\ 2.2]{Proudfoot}).

By Lemma \ref{GIT-walls-lem}(ii), for $\chi\in\int(\bC)$ all the points of the open subset
$U_\a\sub\A^{4g^2-3g}$ given by $\prod_{i\neq j}\a_{ij}\neq 0$, 
are $\chi$-stable.
It remains to observe that for any smooth curve $C$ and generic points $p_1,\ldots,p_g$
one has $\a_{ij}(C,p_1,\ldots,p_g,v_1,\ldots,v_g)\neq 0$ for every pair $i\neq j$, as follows
from Lemma \ref{constants-vanishing-lem}(i), so the corresponding point of $\SS_g$ lies in $U_\a$.

\noindent
(ii) Let $x=(C,p_1,\ldots,p_g,v_1,\ldots,v_g)\in \SS_g$ with $C$ smooth.
By Lemma \ref{constants-vanishing-lem}(iii), 
for any $i$ either there exists $j\neq i$ such that $\a_{ij}\neq 0$ or there exists $j\neq i$ such that $\b_{ij}\neq 0$.
Hence, for any $i$ the set $\Om(x)\sub \Z^g$ contains either a vector $2e_i-e_j$ or $3e_i-e_j$ for some $j\neq i$. 
By Lemma \ref{cone-lem} below (applied to appropriate rescalings of these vectors in $\Om(x)$), this implies the inclusion 
$\bC_0\sub\bC_{\Om(x)}$. Hence, $\chi\in\int(\Om(x))$ and by Lemma \ref{GIT-walls-lem}(ii), $x$ is $\chi$-stable.

\noindent
(iii) Given $\chi\in\int(\bC_0)$, assume that $\chi\in \Sigma$. Then 
$\chi=\sum_{p=1}^r t_p\om_p$, where $t_i>0$ and $\om_1,\ldots,\om_r\in\Om$
lie in a hyperplane in $X(T)\ot\R$.
Then since $\chi$ has positive coordinates, for every $i$ there exists $p$ such that $\om_p$ has
positive $i$th coordinate, i.e., $\om_p$ is one of the vectors $e_i$, $2e_i-e_j$ or $3e_i-e_j$ for some $j\neq i$.
Thus, we can again apply Lemma \ref{cone-lem} to conclude that $\span(\om_1,\ldots,\om_r)=X(T)\ot \R$,
which is a contradicton.
\ed

\begin{lem}\label{cone-lem}
The subcone $\R_{\ge 0}v_1+\ldots+ \R_{\ge 0}v_n\sub \R^n$ spanned by vectors of the form
$$v_i=e_i-a_ie_{\si(i)}, \ \ i=1,\ldots,n,$$
where $0\le a_i<1$, and $\si$ is a map from $\{1,\ldots,n\}$ to itself,
%, such that $\si(i)=i$ implies $a_i=0$,
contains every vector $e_i$.
\end{lem}

\Pf . First, assume that all $a_i$ are positive, and $\si(i)\neq i$ for every $i$.
Renumbering the indices we can assume that $\si(1)=2, \si(2)=3, \ldots, \si(m)=1$ for some
$2\le m\le n$. Then
$$v_1+a_1v_2+\ldots+a_1\ldots a_{m-1}v_m=(1-a_1\ldots a_m)e_1,$$
etc., which shows that $e_1,\ldots,e_m$ belong to $\bC(v_1,\ldots,v_n):=\R_{\ge 0}v_1+\ldots+ \R_{\ge 0}v_n$.
Similarly, we have $e_i\in \bC(v_1,\ldots,v_n)$ whenever $i$ is contained in a subset on which $\si$ acts
as a cyclic permutation (a {\it $\si$-cycle}). Now for any $i$ there exists $r\ge 1$ such that
$j=\si^r(i)$ is contained in a $\si$-cycle, so that $e_j\in \bC(v_1,\ldots,v_n)$. This easily implies that
$e_i\in \bC(v_1,\ldots,v_n)$ (by induction on $r$).

Let $I$ be the set of $i$ such that either $a_i=0$ or $\si(i)=i$. For $i\in I$ we have
$e_i=(1-a_i)^{-1}v_i$, so as above we deduce that $e_i\in \bC(v_1,\ldots,v_n)$ for any $i$
such that $\si^r(i)\in I$ for some $r\ge 1$. Now we can apply the argument from the first part of the proof
to the remaining set of indices.
\ed

\begin{cor}\label{GIT-cor} 
For every character $\chi\in\bC\setminus\Sigma$ the quotient stack
$(\SS_g)^s_\chi/\G_m^g$, which is an open substack of $\UU^{ns,a}_{g,g}$, is proper over $\Z[1/6]$.
\end{cor}

\Pf .
By Lemma \ref{GIT-walls-lem}(iii), for such $\chi$ there are no strictly $\chi$-semistable points.
Hence, the map 
$$\pi: (\SS_g)^s_\chi/\G_m^g\to \SS_g\sslash_\chi \G_m^g$$
is a tame moduli space (see \cite[Prop.\ 7.7]{Alper}). The map $\pi$ is proper (see \cite[Thm.\ 1.1]{Conrad})
and $\SS_g\sslash_\chi \G_m^g$ is projective by Proposition \ref{GIT-prop}(i). This implies that
$(\SS_g)^s_\chi/\G_m^g$ is proper.
\ed

\begin{rems}
1. Recall that Smyth considered in \cite{Smyth} {\it modular compactifications} of $\MM_{g,n}$, which correspond
to open substacks $\XX$ in $\NN_{g,n}$, the stack of smoothable curves, such that $\XX$ is proper.
Thus, by Corollary \ref{GIT-cor}, for each $\chi\in\bC\setminus\Si$
intersecting $(\SS_g)^s_\chi/\G_m^g$ with the locus of smoothable curves
we get such a modular compactification (with the only difference that we work over $\Z[1/6]$, rather than over $\Z$).
Furthermore, by Lemma \ref{GIT-walls-lem}, these compactifications are indexed by chambers in $\bC\setminus\Si$.
Note that for the chamber $\int(\bC_0)$ we get a modular compactification equipped with an action of $S_n$,
permuting the marked points.

\noindent
2. Given that for $\chi\in\int(\bC_0)$ every point $(C,p_1,\ldots,p_g,v_1,\ldots,v_g)\in\SS_g$ with $C$ smooth is
$\chi$-stable, it would be interesting to find out which
singular curves $(C,p_1,\ldots,p_g)$ correspond to $\chi$-stable points.
\end{rems}

\subsection{Some other consequences of relations}

Let us keep the setup of Section \ref{relations-sec}. We
assume that $6$ is invertible, and choose parameters $t_i$ of order $4$ as in 
Lemma \ref{parameter-lem}. In addition we assume that 
$f_i=f_i[2]$ and $h_i=f_i[3]$ are chosen (by adding a constant) in such a way that
$\ga_{ii}=0$ and $2\vareps_{ii}=3\de_{ii}$.
Then we can use the identities of Sections \ref{relations-sec} and \ref{min-gen-sec} for the corresponding functions
on the affine scheme $\SS_g$.
%polynomial relations between the minimal generators \eqref{main-constants-eq} of the algebra of functions on
%$\SS_g\times \Spec(\Z[1/6])$. where we express $\pi_{ii}$ in terms of $\pi_i$ using \eqref{pi-pi-eq}.

The relation \eqref{q-ij-2-eq} expresses $\pi_{ij}$ in terms of $(\a_*,\b_*,\vareps_*)$.
Plugging this expression into \eqref{pi-a-rel} we get
\begin{equation}\label{a-i-j-main-eq}
%\begin{array}{l}
a_{ij}=\a_{ij}(\vareps_{ji}-2\de_{jj})+
\b_{ji}\b_{ij}-\a_{ji}\vareps_{ij}-\ga_{ji}\ga_{ij}
+\sum_{k\neq i,j}\a_{ik}\a_{kj}\b_{jk}
-\sum_{k\neq i,j}\a_{ik}\a_{jk}\ga_{kj}.
%\end{array}
\end{equation}
Plugging the expression \eqref{a-i-j-main-eq} for $a_{ij}$ into \eqref{abcd-rel} we get after simplifying
\begin{equation}\label{M5-eq1}
\begin{array}{l}
-\a_{ik}\de_{jk}-\a_{jk}\de_{ik}+\b_{ij}\b_{ji}+\a_{ij}(\vareps_{jk}+\vareps_{ji}-2\de_{jj})+
\a_{ji}(\vareps_{ik}-\vareps_{ij})+\sum_{l\neq i,j}\a_{il}\a_{lj}\b_{jl}=\\
(\ga_{ij}-\ga_{ik})(\ga_{ji}-\ga_{jk})+\sum_{l\neq i,j}\a_{il}\a_{jl}(\ga_{lj}-\ga_{lk}).
\end{array}
\end{equation}

On the other hand, using \eqref{p-ij-4-eq} we can express $\eta_{ij}$ in terms of $(\a_*,\b_*,\ga_*)$.
Plugging the result into \eqref{r-rel} and \eqref{pi-theta-eq} we get
\begin{equation}\label{M5-eq2}
\begin{array}{l}
\a_{ik}(\varepsilon_{jk}-\varepsilon_{ji})=(\ga_{ij}-\ga_{ik})\b_{jk}+2\a_{ij}\a_{jk}\ga_{jk}
+\b_{ji}\b_{ik}+\de_{ij}\a_{jk}+\\
\sum_{l\neq i,j,k}\a_{il}\b_{jl}\a_{lk}-\sum_{l\neq j,k}\a_{ij}\a_{jl}^2\a_{lk},
\end{array}
\end{equation}
\begin{equation*}
%\label{}
%\begin{array}{l}
2\b_{ij}\vareps_{ij}=3\a_{ij}^2\de_{ij}+3\a_{ij}\ga_{ij}^2
+(2\vartheta_{ii}-3\pi_{ii})\a_{ij}-
\sum_{k\neq i,j}\a_{ik}^3\b_{kj}+
\sum_{k\neq i,j}(\b_{ik}^2-3\a_{ik}^2\ga_{ik})\a_{kj}.
%\end{array}
\end{equation*}
Using \eqref{pi-i-th-pi-ii-eq} we can rewrite the last equation as
\begin{equation}\label{M6-eq}
%\begin{array}{l}
-\a_{ij}\pi_i=
-2\b_{ij}\vareps_{ij}+3\a_{ij}^2\de_{ij}+3\a_{ij}\ga_{ij}^2-
\sum_{k\neq i,j}\a_{ik}^3\b_{kj}+
\sum_{k\neq i,j}(\b_{ik}^2-3\a_{ik}^2\ga_{ik})\a_{kj}.
%\end{array}
\end{equation}
%Now let us plug the expressions for $\pi_{ij}$ and $\vartheta_{ij}$ given by \eqref{q-ij-2-eq} and \eqref{pr-rel} into
%\eqref{agdpt-rel}.
%\begin{equation}
%\begin{array}{l}
%\a_{jk}\b_{ki}\b_{ik}+\a_{ik}\b_{kj}\b_{jk}+
%(\vareps_{ki}+\vareps_{kj}-2\vareps_{kk})\a_{ik}\a_{jk}+\sum_{l\neq i,k}\a_{jk}\a_{il}\a_{lk}\b_{kl}+\\
%\sum_{l\neq i,j}\a_{ik}\a_{jl}\a_{lk}\b_{kl}+
%(\ga_{ik}-\ga_{ij})\de_{jk}+(\ga_{jk}-\ga_{ji})\de_{ik}=\sum_{l\neq i,j}\a_{il}\a_{jl}\de_{lk}+\a_{ji}\vartheta_{ik}+\a_{ij}\vartheta_{jk}.
%\end{array}
%\end{equation}

%\vartheta_{ji}=\sum_{k\neq i,j}\a_{ik}\b_{jk}\a_{ki}+(\ga_{ij}-\ga_{ii})\b_{ji}+\de_{ij}\a_{ji}
%+\a_{ij}\eta_{ji},

%\subsection{Analyzing the relations}

Note also that expressing $\rho_{ij}$ from \eqref{rho-eq} and substituting this expression into
\eqref{b-ij-main-eq} we get
\begin{equation}
\begin{array}{l}
2\a_{ij}\vartheta_{jj}+\sum_{k\neq i,j}\a_{ik}\a_{jk}\de_{kj}+\a_{ji}\vartheta_{ij}+
\ga_{ji}\de_{ij}
=\a_{ij}\pi_{jj}+\a_{ij}\zeta_{jj}+\b_{ji}\vareps_{ij}+\ga_{ij}\vareps_{ji}\\
+\sum_{k\neq i,j}\a_{ik}\b_{jk}\ga_{kj}+b_{ij}.\nonumber
\end{array}
\end{equation}
Subtracting \eqref{pqe-rel} from this, we get the relation
\begin{equation*}
\begin{array}{l}
\a_{ij}(2\vartheta_{jj}-\pi_{jj}-\zeta_{jj})
+\a_{ji}\vartheta_{ij}-\a_{ik}\vartheta_{jk}+\sum_{k\neq i,j}\a_{ik}\a_{jk}\de_{kj}
-\de_{ik}\b_{jk}=\sum_{l\neq i,j}\a_{il}\b_{jl}(\ga_{lj}-\ga_{lk})\\
-\a_{ij}\zeta_{jk}-\de_{ij}(\ga_{ji}+\ga_{jk})+
\b_{ji}(\vareps_{ij}-\vareps_{ik})+(\ga_{ij}-\ga_{ik})(\vareps_{ji}-\vareps_{jk}).
\end{array}
\end{equation*}
Further, using \eqref{q-ii-u-eq}, \eqref{pi-i-th-pi-ii-eq} and \eqref{c-ij-4-eq} 
we can rewrite this as
\begin{equation}\label{compl-rel}
\begin{array}{l}
\a_{ij}\pi_j=-\a_{ji}\vartheta_{ij}+\a_{ik}\vartheta_{jk}-\sum_{k\neq i,j}\a_{ik}\a_{jk}\de_{kj}+\de_{ik}\b_{jk}\\
-2\a_{ij}\a_{jk}\de_{jk}-\a_{ij}(\ga_{jk})^2-\a_{ij}\sum_{l\neq j}\a_{jl}^2(\ga_{lj}-\ga_{lk})
-\de_{ij}(\ga_{ji}+\ga_{jk})+\\
\b_{ji}(\vareps_{ij}-\vareps_{ik})+(\ga_{ij}-\ga_{ik})(\vareps_{ji}-\vareps_{jk})+
\sum_{l\neq i,j}\a_{il}\b_{jl}(\ga_{lj}-\ga_{lk}).
\end{array}
\end{equation}

\begin{prop}\label{constants-prop2} Assume that $g\ge 3$.
Let $U_1\sub \SS_g$ be the open subset consisting of points $x$ such that
for every $i$ either there exists $j$ with $\a_{ij}(x)\neq 0$ or there exists $j$ with $\a_{ji}(x)\neq 0$.
Let also $U_0\sub U_1$ be the open subset of $x$ such that $\a_{ij}(x)\neq 0$ for all $i\neq j$. 
Then the natural maps to the affine spaces
$$U_1\to \Spec \Z[1/6][\a_{ij},\b_{ij},\ga_{ij},\vareps_{ij}] \ \text{ and }$$
$$U_0\to \Spec \Z[1/6][\a_{ij},\b_{ij},\ga_{ij}]$$ 
are locally closed embeddings (in the right-hand side we treat $\a_{ij},\b_{ij},\ga_{ij},\vareps_{ij}$
as independent variables).
\end{prop}

\Pf . Recall that by Proposition \ref{constants-prop}, the functions on $\SS_g$ are
generated by $\a_{ij}$, $\b_{ij}$, $\ga_{ij}$, $\vareps_{ij}$ and $\pi_i$.
Thus, the first assertion follows from \eqref{M6-eq} and \eqref{compl-rel} (recall also that
$\de_{ij}\in R[C]_0$ by \eqref{bc-rel}).
To check the second we use the relations \eqref{M5-eq1} and \eqref{M5-eq2}
that have the form
$$\a_{ij}(\vareps_{ji}+\vareps_{jk})+\a_{ji}(\vareps_{ik}-\vareps_{ij})=A_{ijk},$$
$$\a_{jk}(\vareps_{ik}-\vareps_{ij})=B_{ijk}$$
with $A_{ijk},B_{ijk}$ depending only on $\a_{ij},\b_{ij},\ga_{ij}$
(we switched $i$ and $j$ in \eqref{M5-eq2}).
Expressing $\vareps_{ik}-\vareps_{ij}$ from the second equation and substituting into the first we get
$$\a_{ij}(\vareps_{ji}+\vareps_{jk})=A_{ijk}-\frac{\a_{ji}}{\a_{jk}}B_{ijk}.$$
Switching $i$ and $k$ we get
$$\vareps_{ik}+\vareps_{ij}=(A_{jik}-\frac{\a_{ij}}{\a_{ik}}B_{jik})/\a_{ji}.$$
Since on the other hand 
$$\vareps_{ik}-\vareps_{ij}=B_{ijk}/\a_{jk},$$
we get an expression for $\vareps_{ij}$ in terms of $\a_{ij},\b_{ij},\ga_{ij}$ .
\ed

\begin{rems}\label{grading-rem} 
1. Recall that the condition $\a_{ij}=0$ for a curve $(C,p_1,\ldots,p_g)\in\UU^{ns,a}_{g,g}(k)$
is equivalent to $h^0(2p_i+D_i-p_j)=2$ (see Lemma \ref{constants-vanishing-lem}).
In particular, for smooth $C$ the condition that $(C,p_1,\ldots,p_g)$ is not
in $U_1$ means that $C$ is hyperelliptic and the points $p_i$ are Weierstrass points.

\noindent
2. The action of the diagonal subgroup $\G_m\sub\G_m^g$ gives the algebra of function on $\SS_g$
a grading so that $\deg \a_{ij}=1$, $\deg \b_{ij}=\deg \ga_{ij}=2$, $\deg \vareps_{ij}=3$ and $\deg \pi_i=4$.
%(the degree is equal to the row number in \eqref{main-constants-eq})
All the relations between the coordinates \eqref{main-constants-eq} on 
$\SS_g$ are homogeneous with respect to this grading.
For example, \eqref{ct-rel} is of degree $3$, \eqref{M5-eq1} and \eqref{M5-eq2} are of degree $4$,
\eqref{agdpt-rel}, \eqref{ab-eps-th-eta-eq}, \eqref{M6-eq} and \eqref{compl-rel} are of degree $5$.
By varying $j$ in \eqref{s-i-eq} and eliminating $s_i$ we get relations of degree $6$. Similarly,
we get relations of degree $6$ by eliminating $u_{ij}$ from \eqref{u-i-j-eq}.
It is possible in principle to write out all the defining relations explicitly using the Gr\"obner basis
approach (see the proof of Lemma \ref{grobner-lem}) but this does not seem to be very illuminating. 
The relations become more manageable on the hyperelliptic locus (see Theorem \ref{hyperell-thm}
below).
\end{rems}

%$$(\a_{ij}\a_{jk}\a_{ki}+\a_{ik}\a_{kj}\a_{ji})(x)\neq 0$$
%for any distinct $(i,j,k)$.

\subsection{Hyperelliptic locus}\label{hyperell-sec}

%In this section we assume that $g\ge 2$.
There is a natural hyperelliptic version of the moduli space $\wt{\UU}^{ns,a}_{g,g}$.

\begin{defi}
Let us define the moduli stack $\wt{\HU}^{ns,a}_{g,g}$ by considering families \break
$(C,p_1,\ldots,p_g;v_1,\ldots,v_g)$ in $\wt{\UU}^{ns,a}_{g,g}$ (where $v_i$ is a nonzero
tangent vector at $p_i$) equipped with an involution
$\tau:C\to C$ such that $\tau\circ p_i=p_i$ and $\tau_*v_i=-v_i$.
We denote by $\HU^{ns,a}_{g,g}$ the similar stack with no fixed choice of tangent vectors
(but we still require that $\tau$ acts as $-1$ on the tangent spaces at $p_i$).
\end{defi}

We will show that the involution with this property is unique and that in the smooth case it
is automatically a hyperelliptic involution.

\begin{prop}\label{hyperell-prop} 
(i) The stack $\wt{\HU}^{ns,a}_{g,g}\times \Spec(\Z[1/6])$ is the
closed subscheme of $\wt{\UU}^{ns,a}_{g,g}\times \Spec(\Z[1/6]$ given as the locus of fixed points of 
the involution given by the action of the element $(-1,\ldots,-1)\in\G_m^g$.
The corresponding subscheme of $\SS_g$ is given 
in terms of the coefficients \eqref{big-generators-eq}  by the equations 
\begin{equation}\label{hyperell-van-eq}
\a_{ij}=d_{ij}=r_{ij}=\de_{ij}=e_{ij}^k=b_{ij}=\vareps_{ij}=g_i^j=0.
\end{equation}

\noindent
(ii) Assume that $(C,p_1,\ldots,p_g)$ is in $\HU^{ns,a}_{g,g}(k)$, where $k$ is a field, and $C$ is smooth.
Then $C$ is hyperelliptic and $p_1,\ldots,p_g$ are Weierstrass points.
\end{prop}

\Pf . (i) Recall that the $\G_m^g$-action resclales the tangent vectors at the marked points. Thus, the fixed
points of $(-1,\ldots,-1)\in\G_m^g$ correspond to curves $(C,p_1,\ldots,p_g)$ such that
there exists an automorphism $\tau$ of $(C,p_1,\ldots,p_g)$ such that $\tau_*v_i=-v_i$. 
Then $\tau^2$ stabilizes $v_i$, which implies that $\tau^2=\id$.

Let $f_i, h_i$ be the canonical generators of the corresponding marked algebra 
$A=H^0(C\setminus D,\OO)$ (see Lemma \ref{marked-alg-main-lem}). Note that the involution $\tau$ defines
an automorphism of $A$ as a marked algebra.
Furthermore, $\tau^*f_i$ and $-\tau^*h_i$ are also canonical generators, hence we have
$$\tau^*f_i=f_i, \ \ \tau^*h_i=-h_i.$$
Therefore, looking at the defining relations \eqref{main-eq-ff}--\eqref{main-eq-h2f3}
we deduce the vanishing
\eqref{hyperell-van-eq}. Conversely, if these equations hold then we can define an involution of $A$
which acts on the marking as $(-1,\ldots,-1)$.

\noindent
(ii) Since $H^0\bigl(C,\om_C(-p_1-\ldots-p_g)\bigr)=H^1\bigl(C,\OO(p_1+\ldots+p_g)\bigr)^*=0$, the projection
$$H^0(C,\om_C)\to \bigoplus_{i=1}^g \om_C|_{p_i}$$
is an isomorphism. Hence, $\tau$ acts as $-\id$ on $H^0(C,\om_C)$, so it acts trivially
on the projectivization of this space. Note that the canonical morphism
$C\to \P H^0(C,\om_C)^*$ is compatible with the action of $\tau$.
If $C$ is not hyperelliptic then we get that $\tau$ acts trivially on $C$,
which is a contradiction. Hence, $C$ is hyperelliptic and $\tau$ is the hyperelliptic involution.
\ed

Using Proposition \ref{constants-prop} we get that over $\Spec(\Z[1/6])$ the functions on 
the hyperelliptic locus are generated by $\b_{ij}$, $\ga_{ij}$ and $\pi_i$. Furthermore, using
the associativity equations (see Lemma \ref{grobner-lem}) 
we can get a complete set of relations between these generators.

\begin{thm}\label{hyperell-thm} Assume $g\ge 2$.
The algebra of functions on the affine scheme $\wt{\HU}^{ns,a}_{g,g}\times \Spec(\Z[1/6])$ is
generated by $\b_{ij}$, $\ga_{ij}$ and $\pi_i$ with the defining relations
\begin{equation}\label{hyperell-main-relations}
\begin{array}{l}
\b_{ij}\b_{ji}=(\ga_{ij}-\ga_{ik})(\ga_{ji}-\ga_{jk}),\\
(\ga_{ij}-\ga_{ik})\b_{jk}+\b_{ji}\b_{ik}=0,\\
-(\pi_i+\ga_{ij}^2)\ga_{ij}+2\b_{ij}^2\ga_{ji}-\sum_{l\neq i,j}\b_{il}^2\ga_{lj}=
-(\pi_i+\ga_{ik}^2)\ga_{ik}+2\b_{ik}^2\ga_{ki}-\sum_{l\neq i,k}\b_{il}^2\ga_{lk},\\
(\pi_i+\ga_{ij}^2)\b_{ji}
+\b_{ij}\ga_{jk}(\ga_{ji}+\ga_{jk})+
\b_{ji}\ga_{ik}(\ga_{ij}+\ga_{ik})+\sum_{l\neq i,j,k}\b_{il}\b_{jl}\ga_{lk}=\\
=\b_{ik}\b_{jk}\ga_{kj}+\b_{jk}\b_{ik}\ga_{ki}+
2\b_{ij}\ga_{ji}^2+\sum_{k\neq i,j}\b_{ik}\b_{jk}\ga_{ki},\\
(\pi_i+3\ga_{ij}^2)\b_{ji}=(\pi_j+3\ga_{ji}^2)\b_{ij},
\end{array}
\end{equation}
where different indices are assumed to be distinct. 
The open part $C\setminus \{p_1,\ldots,p_g\}$
of the universal curve is given by the equations
\begin{equation}\label{hyperell-curve-eq}
\begin{array}{l}
f_if_j=\ga_{ji}f_i+\ga_{ij}f_j+a_{ij},\\
f_ih_j=\ga_{ij}h_j+\b_{ji}h_i,\\
h_i^2=f_i^3+\pi_if_i+\sum_{j\neq i}\b_{ij}^2f_j+s_i, \\
h_ih_j=\b_{ij}f_j^2+\b_{ji}f_i^2+\ga_{ji}\b_{ij}f_j+\ga_{ij}\b_{ji}f_i+\sum_{k\neq i,j}\b_{ik}\b_{jk}f_k+u_{ij},
\end{array}
\end{equation}
where 
\begin{equation}
\begin{array}{l}\nonumber
a_{ij}=\b_{ij}\b_{ji}-\ga_{ij}\ga_{ji}, \\
s_i=-(\pi_i+\ga_{ij}^2)\ga_{ij}+2\b_{ij}^2\ga_{ji}-\sum_{k\neq i,j}\b_{ik}^2\ga_{kj}, \\
u_{ij}=(\pi_i+\ga_{ij}^2)\b_{ji}-2\b_{ij}\ga_{ji}^2-\sum_{k\neq i,j}\b_{ik}\b_{jk}\ga_{ki}.
\end{array}
\end{equation}
\end{thm}

\Pf . By Proposition \ref{hyperell-prop}, the defining relations are obtained by looking at the condition
that the relations \eqref{hyperell-curve-eq} define an algebra with the basis
$(f_i^n, h_if_i^n)$. Now the relations \eqref{hyperell-main-relations} are obtained by the
standard Gr\"obner basis technique (see Lemma \ref{grobner-lem}).
\ed

\begin{ex} In the case of genus $2$ the relations of the above Theorem between the generators 
$\b_{12},\b_{21},\ga_{12},\ga_{21},\pi_1,\pi_2$ reduce to the single relation
$$(\pi_2+3\ga_{21}^2)\b_{12}=(\pi_1+3\ga_{12}^2)\b_{21}.$$
\end{ex}

Assume that $g\ge 4$. Then 
using the second equation in \eqref{hyperell-main-relations} 
we deduce the following set of quartic equations for $(\b_{ij})$ on the hyperelliptic locus:
\begin{equation}\label{quartic-eq}
\b_{ji}\b_{ik}\b_{kl}\b_{lj}+\b_{ki}\b_{il}\b_{lj}\b_{jk}+\b_{li}\b_{ik}\b_{jk}\b_{kl}=0,
\end{equation}
where $i,j,k,l$ are distinct.

\begin{rem}
The locus in $\wt{\HU}^{ns,a}_{g,g}$ corresponding to smooth curves can be identified with a $\G_m^g$-torsor
over the configuration space of $2g+2$ distinct points in $\P^1$, $g$ of which are ordered.
By Lemma \ref{constants-vanishing-lem}, the coordinates $\b_{ij}$ are all nonzero on this locus (and can be computed explicitly in terms of the positions of the points---see \cite[Prop.\ 2.5.5]{FP}). 
Using the second equation in \eqref{hyperell-main-relations}
we can express $\ga_{ij}-\ga_{ik}$ in terms of $(\b_{ij})$
(then the first equation in \eqref{hyperell-main-relations} becomes superfluous).
Also, for $g\ge 3$, using the fourth equation in \eqref{hyperell-main-relations} we can express $\pi_i$ in terms
of the other coordinates. Hence, for $g\ge 3$, on the smooth locus our relations can be viewed as equations for the coordinates $\b_{ij}$ only.
\end{rem}

%Make a choice $\ga_{ii}=0$

%$$f_if_j=\sum_{k\neq i,j}\a_{ik}\a_{jk}f_k+\a_{ji}h_i+\a_{ij}h_j+\ga_{ji}f_i+\ga_{ij}f_j+a_{ij},$$
%$$f_ih_j=\a_{ij}f_j[4]+\ga_{ij}h_j+\de_{ij}f_j+\b_{ji}h_i+\vareps_{ji}f_i+\sum_{k\neq i,j} \a_{ik}\b_{jk}f_k+b_{ij}$$
%$$f_i[4]=f_i^2-\sum_{j\neq i}\a_{ij}^2f_j.$$
%$$h_i^2=f_i^3-3\de_{ii}h_i+(2\vartheta_{ii}-3\pi_{ii})f_i-\sum_{j\neq i}\a_{ij}^3h_j+\sum_{j\neq i}(\b_{ij}^2-3\a_{ij}^2\ga_{ij})f_j+s_i$$
%$$h_ih_j=\b_{ij}f_j[4]+\b_{ji}f_i[4]+\vareps_{ij}h_j+\vareps_{ji}h_i+\vartheta_{ij}f_j+\vartheta_{ji}f_i+\sum_{k\neq i,j}\b_{ik}\b_{jk}f_k+u_{ij}$$

\subsection{Differentials and the canonical embedding}\label{Petri-sec}

Throughout this section we work with curves over an algebraically closed field $k$.

We are going to describe the relation of our picture with Petri's analysis of the defining ideal of the
canonical embedding of a smooth non-hyperelliptic curve $C$ (see \cite{Petri}, \cite[Ch.\ III.3]{ACGH}, \cite[Lec.\ I]{Mumford},
\cite{Saint-Donat}, \cite{Sch}). 
Recall that the starting point of this analysis is to consider $g$ distinct points $p_1,\ldots,p_g\in C$
and a basis $\om_1,\ldots,\om_g$ of $H^0(C,\om_C)$,
where $\om_C$ is the canonical line bundle, such that $\om_i$ vanishes at $p_j$ for $i\neq j$.
Such a basis exists precisely when $h^0(p_1+\ldots+p_g)=1$ and is uniquely determined by the
marked points up to rescaling (since by Serre duality $h^0\bigl(\om_C(-p_1-\ldots-p_g)\bigr)=h^1(p_1+\ldots+p_g)$).

First, we observe an important relation between the basis $(\om_i)$ and
the canonical formal parameters $t_i$ at $p_i$ obtained from Lemma \ref{parameter-lem}. 
Assume for a moment that the characteristic is zero.
Then these formal parameters are uniquely determined by a collection of
nonzero tangent vectors $v_i$ at $p_i$ by the condition that $\lan v_i, t_i\ran=1$ for each $i$ and
that for every $n\ge 2$ and every $i=1,\ldots,g$ there exists a rational function $f_i[n]\in H^0\bigl(C,\OO(np_i+D_i)\bigr)$,
where $D_i=\sum_{j\neq i}p_j$. We claim that in fact for each $i=1,\ldots,g$, one can rescale $\om_i$ so that
$\om_i=dt_i$ in the formal neighborhood of $p_i$. If the characterstic is positive we consider
the parameters $t_i$ up to $N$th order (where $N<\cha(k)$) and this statement has to be modified accordingly.

\begin{prop}\label{differentials-prop} 
Assume that $\cha(k)>N$ for some $N>1$.
Let $C$ be a smooth projective curve over $k$, 
$p_1,\ldots,p_g\in C$ distinct points such that
$h^0(p_1+\ldots+p_g)=1$, and let $v_1,\ldots,v_g$ be nonzero tangent vectors at these points.

\noindent
(i) The canonical formal parameters $t_i$ of order $N$ at $p_i$'s from Lemma
\ref{parameter-lem} are characterized by the property that
$$\om_i\equiv dt_i \mod \om_C(-Np_i)$$ 
in $\om_C/\om_C(-Np_i)$, for $i=1,\ldots,g$, where
$(\om_i)$ is a basis of $H^0(C,\om_C)$ such that $\om_i$ vanishes at $p_j$ for $j\neq i$,
and $\lan \om_i,v_i\ran=1$.

\noindent
(ii) For $i\neq j$ and $n\le N$,
let $p_{ij}[n]$ be the coefficient of $\frac{1}{t_j}$ in the Laurent series of $f_i[n]$ at $p_j$. Then
the expansion of $\om_i$ near $p_j$ (where $i\neq j$) has form
$$\om_i=-\sum_{n=2}^N p_{ji}[n]t_j^{n-1}dt_j \mod (t_j^N dt_j).$$
\end{prop}

\Pf . (i)  For every $n\le N$ and $i=1,\ldots,g$,the rational differential $f_i[n]\om_i$ can have a pole only at $p_i$,
so by the residue theorem, we get $\Res_{p_i}(f_i(n)\om_i)=0$. Thus, if we write $\om_i=\phi_i(t_i)dt_i$
at the formal neighborhood of $p_i$ then we deduce that 
$\phi_i=1 \mod (t_i^N)$ as claimed.

\noindent
(ii) This follows immediately from the residue theorem applied to the rational differentials $f_j[n]\om_i$
for $i\neq j$, since $\Res_{p_i}(f_j[n]\om_i)=p_{ji}[n]$ while $\Res_{p_j}(f_j[n]\om_i)$ is equal to the
coefficient of $t_j^{n-1}dt_j$ in the expansion of $\om_i$ at $p_j$. 
\ed

Recall that Petri proceeds by assuming that the differentials $\om_1$ and $\om_2$ have no common zeroes
(which is true if the marked points are sufficiently generic) and shows that the following relations
hold in $H^0(C,\om^2)$:
\begin{equation}\label{Petri-rel-1}
\om_i\om_j=\sum_{k>2}(\la_{ijk}\om_1+\mu_{ijk}\om_2)\om_k+\nu_{ij}\om_1\om_2,
\end{equation}
where $i\neq j$, $i,j>2$, for some constants $\la_{ijk}$, $\mu_{ijk}$ and $\nu_{ij}$. 
In addition, if for $i>2$, $\eta_i\in H^0(\om_C)$ is a linear combination of $\om_1$
and $\om_2$ with double zero at $p_i$ then for $i,j>2$, one has cubic relations in $H^0(C,\om_C^3)$ of the
form
\begin{equation}\label{Petri-rel-2}
\eta_i\om_i^2-\eta_j\om_j^2=\sum_{k>2}(a_{ijk}\om_1^2+b_{ijk}\om_1\om_2+c_{ijk}\om_2^2)\om_k+
d_{ij}\om_1^2\om_2+e_{ij}\om_1\om_2^2,
\end{equation}
for some constants $a_{ijk}$, $b_{ijk}$, $c_{ijk}$, $d_{ij}$ and $e_{ij}$.

Furthermore, Petri shows that the ideal of $C$ in the canonical embedding is generated by the relations
\eqref{Petri-rel-1}, \eqref{Petri-rel-2} for $g\ge 4$ (for generic points $p_1,\ldots,p_g$).

Now let us assume that $\cha(k)\neq 2, 3$
and let us choose canonical parameters $t_i$ of order $4$ as in Lemma \ref{parameter-lem}.
Assume also that
$H^0\bigl(C,\om_C^2(-3D)\bigr)=0$, where $D=p_1+\ldots+p_g$ (this is true for $p_1,\ldots,p_g$
generic). Then any quadratic differential is uniquely
determined by its expansions up to $t_i^2 dt_i^{\ot 2}$ at $p_i$, for
$i=1,\ldots,g$. Hence, the coefficients $\la_{ijk}$, $\mu_{ijk}$ and $\nu_{ij}$ 
from \eqref{Petri-rel-1} are determined from linear equations
obtained by looking at such expansions. Proposition \ref{differentials-prop}(ii) implies that 
these linear equations depend only on $\a_{ij}=p_{ij}[2]$ and $\b_{ij}=p_{ij}[3]$.

For generic curve of genus $g\ge 5$ the quadratic relations generate the ideal of $C$, so let us now show
how to determine the coefficients of the cubic relation \eqref{Petri-rel-2} in the case $g=4$ (assuming
that the marked points are sufficiently generic).
Note first that we can take 
$$\eta_i=\a_{i2}\om_1-\a_{i1}\om_2.$$
Next, we observe that for generic $p_1,\ldots,p_4$ we have $H^0\bigl(C,\om_C^3(-4D)\bigr)=0$.
Indeed, this is equivalent to the vanishing $H^1\bigl(\om_C^{-2}(4D)\bigr)=0$ which follows from the fact that the map
$D\mapsto \om_C^{-2}(4D)$ from $S^4C$ to $\Pic^4(D)$ is dominant.
Thus, any element of $H^0(C,\om_C^3)$ is determined uniqely by expansions up to $t_i^3 dt_i^{\ot 3}$ at $p_i$
(these expansions are well defined by our parameters of order $4$).
Using Proposition \ref{differentials-prop}(ii) we see that such expansions of the terms of \eqref{Petri-rel-2}
give equations on the coefficients that depend only on $\a_{ij}$, $\b_{ij}$ and $\eta_{ij}=p_{ij}[4]$.

This leads to the following result, where we use the notation of Section \ref{relations-sec}.

\begin{prop}\label{genus-4-5-prop} Assume $\cha(k)\neq 2,3$.
Then for $g\ge 4$ the isomorphism class of a generic pointed curve 
$(C,p_1,\ldots,p_g)$ is determined by
the constants $(\a_{ij},\b_{ij},\ga_{ij})$ (viewed up to an action of $\G_m^g$). For $g\ge 5$
this isomorphism class is determined by the constants $(\a_{ij},\b_{ij})$. 
\end{prop}

\Pf . We use Petri's theorem that $C$ is cut out in the canonical embedding by the equations \eqref{Petri-rel-1}
and \eqref{Petri-rel-2} (only by \eqref{Petri-rel-1} for $g\ge 5$)
together with the above considerations on expressing the coefficients of these
equations in terms of $\a_{ij}$, $\b_{ij}$ and $\eta_{ij}$. It remains to observe that
that $\eta_{ij}$ is given by some universal polynomial expression in $(\a_{ij},\b_{ij},\ga_{ij})$ (see Lemma 
\ref{two-coef-lem}(ii))
and that the point $p_i$ is recovered as intersection of zero divisors of $\om_j$ for $j\neq i$.
\ed

\begin{rem} By Theorem \cite[Thm.\ 3.2.1]{FP}, if the characteristic is zero, then
for $g\ge 6$ the isomorphism class of $(C,p_1,\ldots,p_g)$ is determined 
by the constants $(\a_{ij})$ alone (viewed up to $\G_m^g$-action). On the other hand, for $g=5$
the rational map $\MM_{5,5}\to \A^{15}$
given by $\a_{ij}$ modulo the $\G_m^5$-action is dominant and so its generic fiber has dimension $2$ 
(see \cite[Thm.\ 5.2.2]{FP}). 
\end{rem}

Now we are going to present a different way to determine quadratic relations for the canonical embedding of 
$C$, independent from Petri's analysis (see Proposition \ref{quadratic-rel-prop} below).

We need the following result based on the Serre duality.

\begin{lem}\label{Serre-inj-lem} 
Assume that for some line bundle $L$ on $C$ and an effective divisor $E\sub C$ one has 
$H^0\bigl(C,\om_C\ot L^{-1}(-E)\bigr)=0$. Then the map
$$H^0(\om_C\ot L^{-1})\to H^0\bigl(L(E)/L\bigr)^*: \eta\mapsto \left(s\mapsto\sum_{p\in E} \Res_{p}(\eta s)\right)$$
is injective.
\end{lem}

\Pf . Recall that for a line bundle $L$ the Serre duality isomorphism
$$H^0(\om_C\ot L^{-1})\rTo{\sim} H^1(L)^*$$
is given by $\eta\mapsto \bigl(\a\mapsto \tr(\eta\a)\bigr)$, where $\tr:H^1(\om_C)\to k$ is the canonical trace map.
Let us consider the commutative diagram
\begin{diagram}
H^0\bigl(L(E)/L\bigr) &\rTo{\de_L}& H^1(L)\\
\dTo{\eta}&&\dTo{\eta}\\
H^0\bigl(\om_C(E)/\om_C\bigr) &\rTo{\de_{\om}}& H^1(\om_C)
\end{diagram}
where the vertical arrows are given by the multiplication with $\eta\in H^0(\om_C\ot L^{-1})$,
while the horizontal arrows are the connecting homomorphisms in the natural long exact sequences.
It is well known that for $\xi\in H^0\bigl(\om_C(E)/\om_C\bigr)$ one has
$$(\tr\circ \de_{\om})(\xi)=\sum_{p\in E} \Res_{p}(\xi).$$
Hence, by the commutativity of the above diagram we get that
the composed map
$$H^0(\om_C\ot L^{-1})\rTo{\sim} H^1(L)^*\rTo{\de_L^*} H^0\bigl(L(E)/L\bigr)^*$$
is given by $\eta\mapsto \bigl(s\mapsto \sum_p \Res_{p}(\eta s)\bigr)$.
It remains to observe that the condition $H^0\bigl(\om_C\ot L^{-1}(-E)\bigr)=0$ is equivalent to the vanishing
of $H^1\bigl(L(E)\bigr)$ which implies the surjectivity of the map $\de_L$, and hence, the injectivity of
the dual map $\de_L^*$.
\ed

\begin{prop}\label{quadratic-rel-prop}
Assume $\cha(k)\neq 2,3$ and 
$H^0\bigl(C,\om_C^2(-3D)\bigr)=0$. Let $I_2\sub S^2 H^0(C,\om_C)$ be the kernel of the multiplication map
\begin{equation}\label{canonical-mult-map}
S^2 H^0(C,\om_C)\to H^2(C,\om_C^2).
\end{equation}
Let us denote by $\La_2\sub S^2 H^0(C,\om_C)$ the subspace spanned by $\om_i\om_j$ where $i\neq j$.
Let also denote by $V$ the $2g$-dimensional vector space with the basis $(v_i,w_i)_{i=1,\ldots,g}$.
Then $I_2$ is contained in $\La_2$ and coincides with the kernel of the map
$$\La_2\to V: \om_i\om_j\mapsto \sum_{k\neq i,j}\a_{ki}\a_{kj}w_k-\b_{ji}w_j-\b_{ij}w_i-\a_{ji}v_j-\a_{ij}v_i.$$
\end{prop}

\Pf . Applying Lemma \ref{Serre-inj-lem} to the line bundle $L=\om_C^{-1}$ and the divisor $E=3D$, we
derive that the map
$$\kappa:
H^0(\om_C^2)\to H^0\bigl(\om_C^{-1}(3D)/\om_C^{-1}\bigr)^*: \eta\mapsto 
\left(v\mapsto\sum_{i=1}^g \Res_{p_i}(\eta v)\right)$$
is injective.
Thus, the kernel of of the multiplication map \eqref{canonical-mult-map} coincides with the kernel of the map
$$S^2 H^0(\om_C)\to H^0\bigl(\om_C^{-1}(3D)/\om_C^{-1}\bigr)^*: \om\om' \mapsto 
\left(v\mapsto\sum_{i=1}^g \Res_{p_i}(\om\om' v)\right).$$
Let $(u_i,v_i,w_i)_{i=1,\ldots,g}$ be the basis in $H^0\bigl(\om_C^{-1}(3D)/\om_C^{-1}\bigr)^*$,
dual to the basis 
$$(t_i^{-1}\frac{\del}{\del t_i}, t_i^{-2}\frac{\del}{\del t_i}, t_i^{-3}\frac{\del}{\del t_i})_{i=1,\ldots,g}$$ 
of
$H^0\bigl(\om_C^{-1}(3D)/\om_C^{-1}\bigr)$ (where $t_i$ are the canonical parameters of order $4$ at $p_i$).
Then using Proposition \ref{differentials-prop}(ii) we find
\begin{align*}
&\kappa(\om_i^2)=u_i+\sum_{j\neq i}\a_{ji}^2 w_j,\\
&\kappa(\om_i\om_j)=\sum_{k\neq i,j}\a_{ki}\a_{kj}w_k-\b_{ji}w_j-\b_{ij}w_i-\a_{ji}v_j-\a_{ij}v_i, \text{ for }i\neq j,
\end{align*}
where $\a_{ij}=p_{ij}[2]$, $\b_{ij}=p_{ij}[3]$.
Thus, we have a commutative diagram
\begin{diagram}
S^2 H^0(\om_C) &\rTo{\kappa}& H^0\bigl(\om_C^{-1}(3D)/\om_C^{-1}\bigr)^*\\
\dTo{}&&\dTo{}\\
\lan \om_i^2 \ |\ i=1,\ldots,g\ran &\rTo{\sim}& \lan u_i \ |\ i=1,\ldots,g\ran
\end{diagram}
where the vertical arrows are the coordinate projections on the subspaces generated by 
$(\om_i^2)_{i=1,\ldots,g}$ and $(u_i)_{i=1,\ldots,g}$, respectively.
Hence, the kernel of $\kappa$ coincides with the kernel of its restriction to the subspace spanned by
$\om_i\om_j$ with $i\neq j$. Now our statement follows from the explicit formula for $\kappa$.
\ed

\section{Computation of higher products on a curve}\label{Cech-ainf-sec}

\subsection{Cech resolutions}\label{resolution-sec}

Let $\pi:C\to \Spec(R)$ be a flat proper family of curves equipped with a relative effective Cartier
divisor $D\sub C$,
such that $\pi$ is smooth near $D$ and $U=C\setminus D$ is affine.
Then for every quasicoherent sheaf $\FF$ on $C$ we can consider the two-term complex 
$K^\bullet(\FF)=K^\bullet_D(\FF)$ with 
\begin{equation}\label{K-D-def}
\begin{array}{l}
K^0(\FF)=\underset{n}{\varprojlim} H^0\bigl(C,\FF/\FF(-nD)\bigr) \oplus H^0(U,\FF),\nonumber \\
K^1(\FF)=\underset{m}{\varinjlim}\ \underset{n}{\varprojlim}  H^0\bigl(C,\FF(mD)/\FF(-nD)\bigr)
\end{array}
\end{equation}
and the differential
$$d(s_0, s)=\kappa(s)-\iota(s_0),$$
where we use natural maps $\kappa:H^0\bigl(C,\FF/\FF(-nD)\bigr)\to K^1(\FF)$ and 
$\iota:H^0(U,\FF)\to K^1(\FF)$.

The construction of $K^\bullet(\FF)$ immediately generalizes to the case when $\FF$ is a bounded complex
of vector bundles (by taking the total complex of the corresponding bicomplex).
Furthermore, if $\AA$ is a complex of quasicoherent sheaves 
equipped with the structure of an $\OO$-dg-algebra then
we can equip the complex $K^\bullet(\AA)$ with a structure of dg-algebra by 
using the natural componentwise multiplication on $K^0(\AA)$ and using
the multiplications 
\begin{equation}\label{gen-dg-products-eq}
\begin{array}{l}
K^0(\AA)\ot K^1(\AA)\to K^1(\AA): (s_0,s)\cdot u=\iota(s_0)\cdot u, \nonumber \\
K^1(\AA)\ot K^0(\AA)\to K^1(\AA): u\cdot (s_0;s)=u\cdot \kappa(s),
\end{array}
\end{equation}
where in the right-hand side we use the natural product on $K^1(\AA)$.

\begin{lem}\label{Cech-lem} 
(i) For a quasicoherent sheaf $\FF$ on $C$
there is a natural isomorphism in the derived category of $R$-modules
$R\Ga(C,\FF)\simeq K^\bullet(\FF)$.

\noindent (ii)
Assume that $R$ is a finitely generated $k$-algebra, where $k$ is a field.
Let $\AA(V)=\und{\End}_{\OO}(V)$ be the endomorphism dg-algebra over $\OO$ of a 
bounded complex of vector bundles $V$ on $C$. Then the dg-algebra
$K^\bullet\bigl(\AA(V)\bigr)$ is quasi-isomorphic to the dg-algebra of endomorphisms of $V$ computed
using any dg-enhancement of $D\bigl(\Qcoh(C)\bigr)$.
\end{lem}

\Pf . (i) Let $j:C-p \to C$ be the natural open embedding. Then $j_*j^*\FF=\underset{m}{\varinjlim} \FF(mD)$.
We have the sheaf version $\KK^\bullet$ of the complex
$K^\bullet=K^\bullet(\FF)$ with
$$\KK^0=\underset{n}{\varprojlim} \FF/\FF(-nD)\oplus j_*j^*\FF,$$
$$\KK^1=\underset{m}{\varinjlim}\ \underset{n}{\varprojlim} \FF(mD)/\FF(-nD).$$
with the differentials induced by the natural maps $\underset{n}{\varprojlim} \FF/\FF(-nD)\to\KK^1$ and
$j_*j^*\FF=\underset{m}{\varinjlim} \FF(mD)\to \KK^1$, so that $K^\bullet=H^0(C,\KK^\bullet)$. 
Since $D$ and $U$ are affine, the sheaves $\KK^0$ and $\KK^1$ have no higher cohomology.
It remains to check that $\KK^0\to \KK^1$ is a resolution of $\FF$.
Indeed, for each $m,n>0$ we have an exact sequence
$$0\to \FF\to \FF/\FF(-nD) \oplus \FF(mD) \to \FF(mD)/\FF(-nD)\to 0.$$
By passing to the inverse limit in $n$ we get an exact sequence
$$0\to \FF\to \underset{n}{\varprojlim} \FF/\FF(-nD) \oplus \FF(mD) \to 
\underset{n}{\varprojlim} \FF(mD)/\FF(-nD)\to 0$$
(note that the exactness is preserved since on the left we have a constant inverse system).
Then passing to the direct limit in $m$ we get the exactness of
$$0\to \FF\to \KK^0\to \KK^1\to 0.$$

\noindent
(ii)
Recall that $D\bigl(\Qcoh(C)\bigr)$ and the subcategory $\Per(C)$ of perfect complexes both
have unique dg-enhancements, by the results of Lunts-Orlov (see \cite[Cor.\ 2.11, Thm.\ 2.12]{LO}).
Now we observe that for a pair of bounded complexes $V$ and $W$ we can consider the complex
$K^\bullet(V^\vee\ot W)$ and get a dg-category structure on these using products similar to
\eqref{gen-dg-products-eq}. By part (i), this is a dg-enhancement of $\Per(C)$, and
$K^\bullet\bigl(\AA(V)\bigr)$ is precisely the corresponding dg-endomorphism algebra of $V$.
\ed

\subsection{dg-endomorphism algebra of $\OO_C\oplus\OO_{p_1}\oplus\ldots\oplus\OO_{p_g}$}

Let $\pi:C\to \Spec(R)$, $p_1,\ldots,p_g:\Spec(R)\to C$ be a family of curves in 
$\wt{\UU}^{ns,a}_{g,g}(R)$, where $R$ is a finitely generated $k$-algebra, and let
$v_i\in\OO_C(p_i)/\OO_C$ be the corresponding trivializations (for example, we can consider the universal
family over $\Z[1/6]$).
We are going to apply the construction of Section \ref{resolution-sec} for the divisor $D=p_1+\ldots+p_g$
to get an explicit presentation for dg-endomorphisms
of the object
\begin{equation}\label{G-eq}
G=\OO_C\oplus\OO_{p_1}\oplus\ldots\oplus\OO_{p_g}.
\end{equation}
Here we identify each $p_i:\Spec(R)\to C$ with its image, which is a divisor in $C$. 
The reason we chose this object is the fact that the corresponding $\Ext$-algebra is independent of the curve
(this was also the starting point of \cite{FP}).

\begin{lem}\label{Ext-lem}
One has a natural isomorphism of $R$-algebras $\Ext^*(G,G)\simeq E_{g,\Z}\ot R$, where
$E_{g,\Z}$ is the natural $\Z$-form of the algebra \eqref{Eg-eq}.
\end{lem}

\Pf . We have a natural identification $R\simeq\Hom(\OO,\OO_{p_i})$, and we denote by 
$A_i\in \Hom(\OO_C,\OO_{p_i})$ the corresponding generator.
The resolution $[\OO_C(-p_i)\to\OO_C]$ for $\OO_{p_i}$ gives isomorphisms
$$\Ext^1(\OO_{p_i},\OO)\simeq H^0\bigl(\OO_C(p_i)\bigr)/H^0(\OO_C)\simeq H^0\bigl(C,\OO_C(p_i)/\OO_C\bigr),$$
where the last isomorphism follows from the vanishing of $H^1\bigl(\OO_C(p_i)\bigr)$.
Let us denote by $B_i\in \Ext^1(\OO_{p_i},\OO)$ the generator corresponding to the chosen
trivializations $R\simeq H^0\bigl(C,\OO_C(p_i)/\OO_C\bigr)$. The same resolution also gives an isomorphism
$$\Ext^1(\OO_{p_i},\OO_{p_i})\simeq H^0\bigl(\OO_C(p_i)/\OO_C\bigr)\simeq R$$
with the generator given by the composition $Y_i=A_i B_i$.
Finally, the exact sequence 
$$0\to \OO_C\to \OO_C(p_1+\ldots+p_g)\to \bigoplus_{i=1}^g \OO_C(p_i)/\OO_C\to 0$$
together with the vanishing of $H^1(C,\OO)$ gives an isomorphism
$$\bigoplus_{i=1}^g \OO_C(p_i)/\OO_C\to H^1(C,\OO)=\Ext^1(\OO_C,\OO_C),$$
with the generators given by the compositions $X_i=B_i A_i$.
\ed

To calculate the dg-endomorphisms $G$, we first use the resolution
$P:=[\OO_C(-D)\to\OO_C]$ for the sheaf $\bigoplus_i\OO_{p_i}$.
Then we consider the bundle of dg-algebras 
$$\AA=\AA(C,p_1,\ldots,p_g)=\und{\End}(\OO_C\oplus P)$$
and the corresponding dg-algebra 
$$E^{dg}_{(C,p_1,\ldots,p_g)}=K^\bullet_D(\AA).$$ 

Let us fix some formal parameters $t_i$ at $p_i$, compatible with $v_i$, i.e., sections of \break
 $\underset{n}{\varprojlim} H^0\bigl(C,\OO_C(-p_i)/\OO_C(-np_i)\bigr)$ inducing isomorphisms
$$R[[t_i]]\rTo{\sim} \underset{n}{\varprojlim} H^0\bigl(C,\OO_C/\OO_C(-np_i)\bigr)$$
and such that $\lan v_i, t_i\mod \OO_C(-2p_i)\ran=1$.
Note that we also have the induced isomorphisms
$$R((t_i))\rTo{\sim} \underset{m}{\varinjlim}\ \underset{n}{\varprojlim} H^0\bigl(C,\OO_C(mp_i)/\OO_C(-np_i)\bigr).$$
Hence, for any integer $n$ we have an identification of
$K^1\bigl(\OO_C(nD)\bigr)$ with $\bigoplus_{i=1}^g R((t_i))$.
For an element $a(t_i)\in R((t_i))$ we denote by $[a(t)]$
the corresponding element of $K^1\bigl(\OO_C(nD)\bigr)$.

We have a direct sum decomposition
$$E^{dg}_{(C,p_1,\ldots,p_g)}=K_\OO\oplus K_{\OO,P}\oplus K_{P,\OO}\oplus K_{P,P},$$
where $K_\OO=K^\bullet_D(\OO)$ and $K_{P_1,P_2}=K^\bullet_D(P_2\otimes P_1^\vee)$.
We denote (local) sections of the $0$th term of $P$ by $\be\cdot f$, where $f\in\OO$,
and local sections of the $-1$th term of $P$ by $\bu\cdot f$, where $f\in\OO(-D)$.

We denote elements of $K_\OO$ as
$$v+f+[a],$$
where $v\in \bigoplus_{i=1}^g t_i R[[t_i]]$, $f\in \OO(U)$, $a\in \bigoplus_{i=1}^g R((t_i))$,
$v$ and $f$ have degree $0$ and $[a]$ has degree $1$.
The differential on $K_{\OO,\OO}$ is given by
$$d_\OO(v+f+[a])=[f-v],$$
where in the right-hand side we use the projection
$$\OO(U)\to K^1(\OO_C)\simeq\bigoplus_{i=1}^g R((t_i))$$
to view $f$ as an element of $\bigoplus_{i=1}^g R((t_i))$.

The summand $K_{\OO,P}$ decomposes as a graded space as
$$\bu\cdot \left(\bigoplus_{i=1}^g t_iR[[t_i]]\oplus \OO(U)\right)[1]\oplus
\bu\cdot \bigoplus_{i=1}^g R((t_i))\oplus 
\be\cdot \left(\bigoplus_{i=1}^g R[[t_i]]\oplus \OO(U)\right)\oplus 
\bigoplus_{i=1}^g \be\cdot R((t_i))[-1],$$
where $U=C-D$.
We will write elements of $K_{\OO,P}$ as formal sums
$$\bu\cdot v+\bu\cdot f+\bu\cdot [a]+\be\cdot w+\be\cdot h+\be\cdot [b],$$
where $v\in \bigoplus_{i=1}^g t_i R[[t_i]]$, $w\in \bigoplus_{i=1}^g R[[t_i]]$, $a,b\in\bigoplus_{i=1}^g R((t_i))$, 
$f,h\in \OO(U)$.
Here we treat $a,b,v,w,f,h$ as having degree $0$, and use the convention that $\deg(\bu)=-1$, $\deg(\be)=0$
and $\deg([x])=\deg(x)+1$.

Similarly, elements of $K_{P,\OO}$ are formal sums
$$v\cdot \bu^*+f\cdot \bu^*+[a]\cdot \bu^*+w\cdot \be^*+h\cdot \be^*+[b]\cdot \be^*,$$
where $\deg(\bu^*)=1$, $\deg(\be)=0$, $v \in \bigoplus_{i=1}^g t_i^{-1}R[[t_i]]$.

Elements of  $K_{P,P}$ are formal sums
\begin{align*}
&\bu\cdot(v_{uu}+f_{uu}+[a_{uu}])\cdot\bu^*+\be\cdot(v_{eu}+f_{eu}+[a_{eu}])\cdot\bu^*+
\bu\cdot(v_{ue}+f_{ue}+[a_{ue}])\cdot\be^*+\\
&\be\cdot(v_{ee}+f_{ee}+[a_{ee}])\cdot\be^*,
\end{align*}
where $v_{uu}\in \bigoplus_{i=1}^g R[[t_i]]$, $v_{eu}\in \bigoplus_{i=1}^g t_i^{-1}R[[t_i]]$ and
$v_{ue}\in \bigoplus_{i=1}^g t_iR[[t_i]]$.

The product on $K^0_{\OO,\OO}$ is simply that of the direct sum of rings. The remaining products
are determined by the rules 
$$\bu^*\bu=1,\ \  \be^*\be=1, \ \ \bu^*\be=\be^*\bu=0,$$
$$f\cdot[a]=0, \ \ v\cdot[a]=[va], \ \ [a]\cdot f=[af], \ \ [a]\cdot v=0,$$
where the product on $\bigoplus_i R((t_i))$ is componentwise.
%REVISE THE PRODUCT IN EARLIER FORMULAS!

The differentials on $K_{\OO,P}$, $K_{P,\OO}$ and $K_{P,P}$ are determined 
by the Leibnitz rule using the formulas 
$$d(\bu)=\be, \ d(\be)=0, \ d(\be^*)=-\bu^*, \  d(\bu^*)=0.$$
For example, the differentials on $K_{\OO,P}$ and $K_{P,\OO}$ are
\begin{align*}
& d(\bu\cdot x+\be\cdot y)=-\bu\cdot d_\OO(x)+\be\cdot x+\be\cdot d_\OO(y),\\
& d(x\cdot\bu^*+y\cdot\be^*)=d_\OO(x)\cdot\bu^*-(-1)^{\deg(y)}y\cdot\bu^*+d_\OO(y)\cdot \be^*
%& d(\bu\cdot x_{uu}\cdot\bu^*+\bu\cdot x_{ue}\be^*+\be\cdot x_{eu}\bu^*+\be\cdot x_{ee}\be^*)=
\end{align*}

By Lemmas \ref{Cech-lem}(i) and \ref{Ext-lem},
the cohomology algebra of $E^{dg}_{(C,p_1,\ldots,p_g)}$ is the algebra
$\Ext^*(G,G)=E_{g,\Z}\otimes R$.

Let us list convenient representatives for the cohomology of our complexes:
For $K_{\OO,P}$ these are
$$A_i=\be\cdot \unit_i+\bu\cdot[\unit_i]\in K^0_{\OO,P}, \ i=1,\ldots,g,$$
where $\unit_i$ is $1\in R[[t_i]]\sub\bigoplus_{i=1}^g R[[t_i]]$.
For $K_{P,\OO}$ we choose
$$ B_i=\frac{1}{t_i}\cdot \bu^*+[\frac{1}{t_i}]\be^*\in K^1_{P,\OO}.$$
We have induced classes 
$$ X_i= B_i A_i=[\frac{1}{t_i}]\in K^1_{\OO,\OO},$$
$$ Y_i=A_i B_i=\be\cdot\frac{1}{t_i}\cdot\bu^*+\be[\frac{1}{t_i}]\be^*\in K^1_{P,P}.$$
To complete the list of representatives we alse need the classes
$$\unit_\OO=\sum_{i=1}^g\unit_i+1\in K^0_{\OO,\OO} \ \text{ and}$$
$$e_{p_i}:=\be\cdot\unit_i\cdot\be^*+\bu\cdot\unit_i\cdot\bu^*+\bu[\unit_i]\be^*\in K^0_{P,P}.$$
Note that we have the following product formulas
$$ B_i Y_i= X_i B_i=0,$$
$$A_i X_i= Y_i A_i=\be[\frac{1}{t_i}].$$
The products involving $e_{p_i}$ and $\unit_{\OO}$ are the same as on the cohomology level.

\subsection{Homotopy operators and formulas for higher products}
\label{homotopy-sec}

Next, we are define the homotopy operator $Q$ on $E^{dg}=E^{dg}_{(C,p_1,\ldots,p_g)}$.
More precisely, it is easy to see that cohomology of $E^{dg}$ and the submodule $\im(d)\sub E^{dg}$
are free $R$-modules, so we can choose complements to $\im(d)$ in $\ker(d)$ and to
$\ker(d)$ in $E^{dg}$. These complements give a homotopy operator $Q$ that maps both these complements
to zero, has the image in the complement to $\ker(d)$ and is inverse to $d$ on $\im(d)$.
Note that in this case $\Pi=\id-Qd-dQ$ is the projector onto the free $R$-submodule of
cohomology representatives. Also, it is clear from the construction that we have 
$$Q^2=\Pi Q=Q\Pi=0.$$

To implement this construction we make a choice for each $n\ge 2$ of a
section $f_i[n]\in H^0\bigl(C,\OO(D+(n-1)p_i)\bigr)$ so that the polar part of $f_i[n]$ at $p_i$ 
has form $t_i^{-n}\mod t_i^{-n+1}R[[t_i]]$ with respect to
the formal parameter $t_i$ at $p_i$. Furthermore, we can choose such sections $f_i[n]$
uniquely up to adding a constant
(modifying them by $f_i[m]$ with $m<n$), so that for each $n\ge 2$ the polar part at $p_i$ is
\begin{equation}\label{f-i-n-def-eq}
f_i[n](t_i)=\frac{1}{t_i^n}+\frac{p_{ii}[n]}{t_i}+\ldots
\end{equation}
for some constants $p_{ii}[n]\in R$. Also, the polar part of $f_i[n]$ at $p_j$ has form
$$f_i[n](t_j)=\frac{p_{ij}[n]}{t_j}+\ldots$$
for some constants $p_{ij}[n]\in R$. Note that at this point we do not make any special choice of formal
parameters $t_i$ at $p_i$. Lemma \ref{parameter-lem} tells that by choosing these parameters in
a special way we may achieve the vanishing of some of the constants $p_{ii}[n]$ (in characteristic zero
of all of them).

Now to define $Q$ we use cohomology representatives defined above. Also, we use $f_i[n]$ to define
the complement to $\ker(d)$. For example, for $K_{\OO,\OO}^0$ this complement is spanned
by $(f_i[n])_{n\ge 2}$ and by $\bigoplus_{i=1}^g R[[t_i]]$.

Below for a series $a(t)=\sum_{i\ge -N} a_it^i\in R((t))$ we denote $a(t)_{\ge n}=\sum_{i\ge n}a_it^i$. 

We define the homotopy operator on $K_{\OO,\OO}$ by 
$$Q([\frac{1}{t_i}])=0, \ Q([v])=-v$$
for $v\in\bigoplus_i R[[t_i]]$ and
$$Q([\frac{1}{t_i^n}])=f_i[n](t_i)_{\ge 0}\unit_i+\sum_{j\neq i}f_i[n](t_j)_{\ge 0}\unit_j+f_i[n].$$
%(f_i[n]-\frac{1}{t_i^n}\frac{p_{ii}[n]}{t_i})\unit_i+\sum_{j\neq i}(f_i[n]-\frac{p_{ij}[n]}{t_j})\unit_j+f_i[n].
On $K_{\OO,P}$ we define $Q$ by
$$Q(\be[b])=\bu[b],\ Q(\bu[a]+\be\cdot v+\be\cdot f)=\bu\cdot v_{>0}+\bu\cdot f.$$
%where $v_{>0}$ is obtained by subtracting constant terms from all the components of $v$. 
On $K_{P,\OO}$ we set
$$Q([a]\bu^*)=[a]\be^*,\ Q([a]\be^*+v\cdot\bu^*+f\cdot\bu^*)=-v_{\ge 0}\be^*-f\cdot\be^*.$$
%where $v_{\ge 0}$ is obtained from $v$ by removing the polar parts.
Finally, on $K_{P,P}$ the homotopy operator is defined by
\begin{align*}
&Q(\be[a_{eu}]\bu^*)=\bu[a_{eu}]\bu^*,\\
&Q(\bu[a_{uu}]\bu^*+\be[a_{ee}]\be^*+\be\cdot v_{eu}\cdot\bu^*+\be\cdot f_{eu}\cdot\bu^*)=
\bu[a_{ee}-(v_{eu})_{<0}]\be^*+\bu\cdot (v_{eu})_{\ge 0}\cdot\bu^*+\bu\cdot f_{eu}\cdot \bu^*,\\
&Q(\bu[a_{ue}]\be^*+\bu\cdot v_{uu}\cdot\bu^*+\bu\cdot f_{uu}\cdot \bu^*+\be\cdot v_{ee}\cdot\be^*+
\be\cdot f_{ee}\cdot\be^*)=\bu\cdot(v_{ee})_{>0}\cdot\be^*+\bu\cdot f_{ee}\cdot\be^*.
\end{align*}

The corresponding projector $\Pi$ onto the cohomology representatives is given by
\begin{align*}
&\Pi(v)=\Pi([v])=\Pi(f_i[n])=0, \ \Pi(1)=\unit_\OO, \ \Pi([\frac{1}{t_i}])= X_i, \\ 
&\Pi([\frac{1}{t_i^n}])=-\sum_j p_{ij}[n] X_j,
 \text{ where }n\ge 2,\\
&\Pi(\bu[a]+\be\cdot v+\be\cdot f)=\sum_i v_i(0)A_i,\\
&\Pi([a]\cdot\be^*+v\cdot\bu^*+f\cdot\bu^*)=\sum_i \Res_i(v) B_i,\\
&\Pi(\bu[a_{uu}]\bu^*+\be[a_{ee}]\be^*+\be\cdot v_{eu}\cdot\bu^*+\be\cdot f_{eu}\cdot\bu^*)=
\sum_i\Res_i(v_{eu}) Y_i,\\
&\Pi(\bu[a_{ue}]\be^*+\bu\cdot v_{uu}\cdot\bu^*+\bu\cdot f_{uu}\cdot\bu^*+\be\cdot v_{ee}\cdot\be^*+
\be\cdot f_{ee}\cdot\be^*)=\sum_i (v_{ee})_i(0)e_{p_i},
\end{align*}
where $\Res_i(v)$ for $v=(v_i)$ is the coefficient of
$\frac{1}{t_i}$ in $v_i$.

Recall (see \cite{Merk}, \cite{KS})
that the above choice of the homotopy operator $Q$ gives homological perturbation formulas for the
$A_\infty$-structure on the cohomology of $E^{dg}$, which we identified with $E_{g,\Z}\ot R$. 
Namely, for cohomology representatives $b_1,\ldots,b_n$, we have
$$m_n(b_1,\ldots,b_n)=\sum_{T}\pm m_T(b_1,\ldots,b_n),$$
where $T$ runs over all oriented planar rooted $3$-valent trees with $n$ leaves (different from the root)
marked by $b_1,\ldots,b_n$
left to right, and the root marked by $\Pi$. The expression
$m_T(b_1,\ldots,b_n)$ is obtained by going down from leaves to the root, applying the
multiplication in $E^{dg}$ at
every vertex and applying the operator $Q$ at every inner edge (see \cite[Sec.\ 6.4]{KS},
\cite[Sec.\ 2.1]{P-ell} for details). 

\begin{lem}\label{ainf-unit-lem} 
All products $m_i$ with $i\ge 3$ involving $e_{p_i}$ or $\unit_\OO$ vanish. Thus,
the obtained $A_\infty$-structure on $E_{g,\Z}\ot R$, viewed as a category with $n+1$ objects,
is strictly unital.
\end{lem}

\Pf . Due to the vanishing $Q^2=\Pi Q=0$ it is enough to check
that the mutliplication by $e_{p_i}$ and by $\unit_\OO$ on the left or the right preserves the image
of $Q$. With $\unit_\OO$ this is clear since it is a unit in $K^0_{\OO,\OO}$. With $e_{p_i}$ this is easy to
check.
\ed

Let us apply the above construction to the universal family $(C,p_1,\ldots,p_g)$
over 
$$\wt{\UU}^{ns,a}_{g,g}\times\Spec(\Z[1/6])=:\Spec(R^\univ).$$ 
Recall that we have a natural $\G_m^g$-action on this moduli space and on the universal curve.
Let us consider the induced $\G_m$-action for the diagonal $\G_m\sub \G_m^g$. Note that
we have the corresponding $\Z$-grading on $R^\univ$.
We can choose all formal parameters $t_i$ at $p_i$ to be $\G_m$-equivariant, i.e., such that
\begin{equation}\label{la-t-i-eq}
\la^*t_i=\la t_i
\end{equation}
(recall that we normalize our formal parameters by $\lan v_i, t_i\mod\OO(-2p_i)\ran=1$ 
and that $\la_*v_i=\la^{-1} v_i$).
Then the corresponding sections $f_i[n]$ satisfying \eqref{f-i-n-def-eq}, where $n\ge 2$, can be chosen to
be $\G_m$-equivariant as well, i.e., such that
$$(\la^{-1})^*f_i[n]=\la^n f_i[n].$$
Therefore, the corresponding homotopy operator $Q$ on $E^{dg}$ 
is also going to be equivariant with respect to the $\G_m$-action induced by the
$\G_m$-action on the universal curve.

\begin{prop}\label{weight-prop} For a $\G_m$-equivariant homotopy operator $Q$, 
the structure constants of the higher product $m_n$ with respect to the standard basis of 
$E_{g,\Z}\ot R^{univ}$ belong to $R^{\univ}_{n-2}$, the component of degree $n-2$ with respect
to the $\Z$-grading of $R^{\univ}$.
\end{prop}

\Pf . Let us denote by $\wgt$ the degree of an element of $E^{dg}$ with respect to the $\G_m$-action
given by the operators $(\la^{-1})^*$.
 %action of $\G_m$ on the universal family. 
Since the differential $d$,
the homotopy $Q$ and the projector $\Pi$ are all homogeneous with respect to this grading,
the operations $m_n$ are also homogeneous.
Now we observe that $\wgt A_i=0$ while $\wgt B_i=1$ (the latter follows from \eqref{la-t-i-eq}). 
Since these elements generate $E_{g,\Z}$, it follows that $(\la^{-1})^*=\la^{\deg}$
on $E_{g,\Z}$, where $\deg$ is the cohomological grading on $E_{g,\Z}$. 
Since $m_n$ lowers the cohomological grading by $n-2$, we derive that $(\la^{-1})^*$ rescales the
coefficients of $m_n$ by $\la^{n-2}$.
\ed

Using Proposition \ref{constants-prop} we deduce the following result.

\begin{cor}\label{weight-cor} 
Consider the $A_\infty$-structure on $E_{g,\Z}\ot R^\univ[1/6]$
associated with some $\G_m$-equivariant homotopy operator $Q$ on $E^{dg}$.
Then the structure constants of $m_n$ 
are polynomials in generators \eqref{main-constants-eq} of degree $n-2$, where $\deg \a_{ij}=1$,
$\deg\b_{ij}=\deg\ga_{ij}=2$, $\deg\vareps_{ij}=3$ and $\deg\pi_i=4$.
\end{cor}

In the hyperelliptic case we can deduce vanishing of $m_n$ with odd $n$.

\begin{cor}\label{hyperell-ainf-cor}
Let $(C,p_1,\ldots,p_g)$ be the universal family over the hyperelliptic locus 
$\wt{\HU}^{ns,a}_{g,g}\times\Spec(\Z[1/6])=\Spec(R^\he)$.
Then the $A_\infty$-structure on $E_{g,\Z}\ot R^{he}$, in the equivalence class given by the homological
perturbation of $E^{dg}$, can be chosen so that $m_n=0$ for all odd $n$.
\end{cor}

\Pf . Let $p:R^\univ\to R^\he$ be the natural homomorphism. Note that it is compatible
with the $\Z$-grading and since $-1\sub\G_m$ acts trivially on $\Spec(R^\he)$ (see Proposition
\ref{hyperell-prop}),
we obtain that the odd degree components in $R^\he$ vanish. It remains to use the $A_\infty$-structure
on $E_{g,\Z}\ot R^{he}$ obtained from some $\G_m$-equivariant choices of $(t_i)$, $(f_i[n])$.
\ed

\begin{rem} Corollary \ref{hyperell-ainf-cor} is
a generalization of the fact that for elliptic curves the corresponding $A_\infty$-structure 
can be chosen so that $m_n=0$ for odd $n$, which was a consequence of a direct computation
in \cite{P-ell}.
\end{rem}

Now let us make some explicit computations. We set 
$$P_i=A_i X_i= Y_i A_i=\be[\frac{1}{t_i}].$$
Note that this is the only double product of $A_i$, $ X_i$, $ Y_i$, $ B_i$ on which $Q$ is not zero.
Thus, when calculating higher products we can restrict the sum over trees to only those trees in which
any two leaves joined to the same internal vertex are either $A_i$ and $ X_i$, or $ Y_i$ and $A_i$.
Note also that
$$Q(P_i)=\bu[\frac{1}{t_i}],$$
$$Q(P_i) B_i=Q(P_i) B_j=Q(P_i) X_i=Q(P_i) X_j=0,$$
$$(Q\circ Y)^n(Q(P_i)) B=(Q\circ Y)^n(Q(P_i)) X=0,$$
$$Q(K^1_{\OO,P}) B_i=Q(K^1_{\OO,P}) X_i=0,$$
$$ X_iQ(K^2_{P,\OO})=Q(K^2_{P,\OO})A_i=0,$$
$$Q(K^2_{P,P})A_i=Q(K^2_{P,P}) Y_i=Q(K^1_{P,P}) Y_i=0,$$
$$ X Q\bigl(Q( B Q(P)) X\bigr)= X Q\bigl(Q( B Q(P)) B\bigr)=0.$$
These vanishings further restrict the types of trees contributing to the higher
products. 

If $2$ is invertible in $R$ then we can choose parameters $t_i$ in such a way that
$p_{ii}[2]=0$ (by Lemma \ref{parameter-lem}). Then from the above observations
one can easily derive the following formulas for $m_3$:
\begin{equation}\label{m3-formulas-eq}
m_3( B_i, Y_i, A_i)=m_3( B_i, A_i, X_i)=-\sum_{j\neq i}\a_{ij} X_j,
\end{equation}
where $\a_{ij}=p_{ij}[2]$, while all other $m_3$ products of the basis elements vanish.
Similarly, if $6$ is invertible in $R$ we can also assume the vanishing of $p_{ii}[3]$ and $p_{ii}[4]$.
The obtained explicit formulas for $m_4$ and $m_5$ can be found in the Appendix.

\section{From moduli of curves to moduli of $A_{\infty}$-structures}\label{moduli-sec}

\subsection{Generalities on $A_\infty$ and $A_n$-structures}\label{ainf-gen-sec}

For a graded associative $S$-algebra $A$ (where $S$ is a commutative ring and $A$ is flat as $S$-module),
we denote the terms of the Hochschild cochain complex of $A$ over $S$ as follows: 
$CH^{s+t}(A/S)_t$ denotes the space of
$S$-multilinear maps $A^{\ot s}\to A$ of degree $t$ (where tensoring is over $S$).
We have the induced bigrading $HH^{s+t}(A/S)_t$ of the Hochschild cohomology.
The corresponding grading by the upper index is compatible with the 
definition of the Hochschild cohomology for $A_\infty$-algebras.

Below we use the notion of $A_n$-structure which is a truncated version of an $A_\infty$-structure
defined by Stasheff (see \cite[Def.\ 2.1]{Stasheff}). 
For a moment let $A$ be a graded $S$-module.
Recall that an $S$-linear $A_n$-structure is given by a collection of $S$-multilinear maps
$$(m_1,\ldots,m_n)\in CH^2(A/S)_1\times\ldots\times CH^2(A/S)_{2-n}$$ 
satisfying the standard $A_\infty$-identities involving only $m_1,\ldots,m_n$ (see below). 
Following \cite[(2.4)]{Stasheff}, 
$A_n$-structures can be described conveniently in terms of truncated bar-construction
$$\Bar_{\le n}(A)=\bigoplus_{i=1}^n T^i_S(A[1]).$$
It has a natural structure of coalgebra over $S$ (without counit), such that it is a sub-coalgebra
of the full bar-construction $\Bar(A)=\oplus_{i\ge 1} T^i_S(A[1])$.
For each cochain $c\in CH^{s+t}(A/S)_t$, where $s\ge 1$, we denote by $D_c$
the corresponding coderivation of $\Bar(A)$ of degree $s+t-1$, 
preserving each sub-coalgebra $\Bar_{\le n}(A)$
(we recover $c$ from the component $\Bar_{\le s}(A)\to A[1]$ of $D_c$).
Then the condition for $m=(m_1,\ldots,m_n)$ to define an $A_n$-algebra structure on $A$ 
is that 
$$D_m^2|_{\Bar_{\le n}(A)}=0.$$ 
We can rewrite this as a collection of identities
\begin{equation}\label{a-inf-id-eq}
\sum_{i=1}^r D_{m_i}D_{m_{r+1-i}}|_{\Bar_{\le n}(A)}=0,
\end{equation}
where $r=1,\ldots,n$.

We denote by $[?,?]$ the supercommutator of coderivations. Recall that
$$[D_c,D_{c'}]=D_{[c,c']},$$
where $[c,c']$ is the Gerstenhaber bracket. 
Also, if $D_c$ has degree $1$ then
$D_c^2$ is still a coderivation, so it corresponds to some cochain ($[c,c]/2$ when $2$ is invertible).
Thus, we can view the identity \eqref{a-inf-id-eq} as the linear equation for coderivations associated 
with some Hochschild cochains in $CH^3(A/S)_{3-r}$. Since such cochains $c$ are uniquely
determined from the restriction $D_c|_{\Bar_{\le r}(A)}$, we deduce the following result.

\begin{lem}\label{An-identities}
The elements $(m_1,\ldots,m_n)\in CH^2(A/S)_1\times\ldots\times CH^2(A/S)_{2-n}$
define an $S$-linear $A_n$-structure on $A$ if and only if
$$\sum_{i=1}^r D_{m_i}D_{m_{r+1-i}}=0,$$
for $r=1,\ldots,n$.
If $2$ is invertible in $S$ then this is equivalent to
$$\sum_{i=1}^r [m_i,m_{r+1-i}]=0,$$
$r=1,\ldots,n$.
\end{lem}

We are interested in {\it minimal} $A_n$-structures, i.e., $A_n$-structures with $m_1=0$.
Then for $n\ge 3$ the product $m_2$ is automatically associative. 
%and we are going to fix it.
When we talk about minimal $S$-linear $A_n$-structures on
a graded associative $S$-algebra $A$, unless otherwise specified, we always assume that $m_2$ is
the given product on $A$.
Note that for a Hochschild cochain $c\in CH^{s+t}(A/S)_t$ we have
$$[D_{m_2},D_c]=D_{m_2}D_c+(-1)^{s+t} D_cD_{m_2}=D_{\de(c)},$$
where $\de(c)=[m_2,c]$ is the Hochschild differential. 

Any $A_{n+1}$-structure induces an $A_n$-structure by forgetting $m_{n+1}$.
The following well-known result states that an obstacle to extending an $A_n$-structure
to an $A_{n+1}$-structure lies in $HH^3(A/S)_{2-n}$
(it is stated without proof as \cite[Lem.\ 2.3]{AAEKO}).

\begin{lem}\label{An-extension-lem}
(i) Let $\AA$ be an associative algebra with generators $D_1,\ldots,D_n$ and defining relations
$$\sum_{i=1}^r D_iD_{r+1-i}=0,$$
for $r=1,\ldots,n$. Set $S=\sum_{i=2}^n D_iD_{n+2-i}$. Then
$$D_1S-SD_1=0.$$

\noindent
(ii) For a minimal $S$-linear $A_n$-structure $m=(m_2,\ldots,m_n)$ on $A/S$ there exists a
Hochschild cocycle $\phi_n(m)\in CH^3(A/S)_{1-n}$ such that
\begin{equation}\label{phi-n-eq}
D_{\phi_n(m)}=\sum_{i=3}^n D_{m_i}D_{m_{n+3-i}}.
\end{equation}
The $A_n$-structure $m$ is extendable to an $A_{n+1}$-structure $(m_2,\ldots,m_n,m_{n+1})$
if and only $\phi_n(m)$ is a coboundary. 
\end{lem}

\Pf . (i) Let us give $\AA$ the grading by $\deg D_i=1$ and use the corresponding supercommutator
$[?,?]$. Then we have
$$[D_1,S]=\sum_{i=2}^n [D_1,D_i]D_{n+2-i}-\sum_{i=2}^n D_i[D_1,D_{n+2-i}].$$
%dab-abd= (da+ad)b-a(db+bd)
Applying the relations we can rewrite the sums in the right-hand side as
$$\sum_{i=2}^n [D_1,D_i]D_{n+2-i}=-\sum_{i\ge 2,j\ge 2,k\ge 2, i+j\le n+1}D_iD_jD_{n+3-i-j},$$
$$\sum_{i=2}^n D_i[D_1,D_{n+2-i}]=\sum_{i\ge 2,j\ge 2,i+j\le n+1}D_{n+3-i-j}D_iD_j.$$
Thus, both sums are equal to
$$\sum_{i\ge 2,j\ge 2,k\ge 2,i+j+k=n+3} D_iD_jD_k,$$
so they cancel out.

\noindent
(ii) The existence of the Hochschild cochain $\phi_n(m)$ follows from the fact that the expression
in the right-hand side of \eqref{phi-n-eq} is a coderivation. The fact
that $\phi_n(m)$ is $\de$-closed follows from (i).
By Lemma \ref{An-identities}, the condition on $m_{n+1}$ to extend $m=(m_2,\ldots,m_n)$ to
an $A_{n+1}$-structure is 
$$[D_{m_2},D_{m_{n+1}}]=-\sum_{i=3}^n D_{m_i}D_{m_{n+3-i}},$$
i.e., $\de(m_{n+1})=-\phi_n(m)$, which implies the assertion.
\ed

\begin{defi}\label{gauge-def} 
The {\it group of gauge transformations} $\fG$ is the group of degree-preserving
coalgebra automorphisms $\a:\Bar(A)\to \Bar(A)$ such that the component $\Bar(A)\to A[1]$
is given by a collection 
$$(f_1=\id,f_2,\ldots)\in CH^1(A/S)_{-1}\times CH^1(A/S)_{-2}\times\ldots.$$
The group of {\it extended gauge transformations} is defined similarly by requiring
$f_1$ just to be invertible.
Note that any such automorphism automatically preserves any sub-coalgebra
$\Bar_{\le n}(A)$ and the condition $f_1=\id$ is equivalent to the condition that
$\a$ acts as identity on every quotient $\Bar_{\le i}(A)/\Bar_{\le i-1}(A)$.
\end{defi}

We usually identify elements of $\fG$ with the corresponding collections $f=(f_1=\id,f_2\ldots)$ and
denote by $\a_f$ the corresponding automorphism of $\Bar(A)$.
Note that the group $\fG$ acts on the set of $A_n$-structures for every $n$: for $f\in\fG$ and an $A_n$-structure
$m$, the new $A_n$-structure $f*m$ is determined by
$$D_{f*m}=\a_f D_m \a_f^{-1},$$
where in the right-hand side we restrict $\a_f$ to $\Bar_{\le n}(A)$.

\begin{defi}
For each $n$ let us denote by $\fG_{\ge n}\sub \fG$ the subgroup of $f=(f_1=\id,f_2,\ldots)$
with $f_i=0$ for $2\le i<n$. To see that this is a subgroup we observe that this vanishing condition is
equivalent to the condition that $\a_f$ acts as identity on all the quotients $\Bar_{\le i}(A)/\Bar_{\le i-n+1}$.
In particular, $\fG_{\ge 2}=\fG$.
\end{defi}

\begin{lem}\label{gauge-lem} 
The subgroup $\fG_{\ge n+1}$ acts trivially on the set of $A_n$-structures.
The subgroup $\fG_{\ge n}$ acts trivally on the set of minimal $A_n$-structures.
For 
$$f=(f_1=\id,0,\ldots,f_{n-1},\ldots)\in \fG_{\ge n-1}$$ and a minimal $A_n$-structure 
$m=(m_2,\ldots,m_n)$ one has $f*m=(m_2,\ldots,m_{n-1},m'_n)$
with 
$$m'_n=m_n\pm \de(f_{n-1}).$$
\end{lem}

The proof is straightforward (cf.\ \cite[Lem.\ 2.2(ii)]{P-ext} for the last assertion).

When we are dealing with an associative algebra $A$ that has an underlying quiver $Q$, so
that $A$ is a quotient of the path-algebra $k[Q]$ by a homogeneous ideal (such as $E_g$), we 
treat $A$ as a category with the objects corresponding to the vertices of the quiver.
Furthermore, we always consider strictly unital $A_\infty$ (and $A_n$)-structures, and
we always consider the reduced Hochschild cochains of $A$ in the above definitions.
Note that this does not change the Hochschild cohomology groups 
(see \cite[(20d)]{Seidel-book}).

We will also use the notion of a homotopy between gauge transformations (see \cite{Keller-intro} and
\cite[Sec.\ 2.1]{P-ext}, where these are called homotopies between strict $A_\infty$-isomorphisms).

\begin{lem}\label{important-homotopy-lem} 
Let $(E,m_2)$ be an associative flat $S$-algebra such that
$HH^1(E/S)_{-i}=0$ for $i=1,\ldots,d-2$, where $d\ge 2$. Suppose $m=(m_\bullet)$ and $m'=(m'_\bullet)$
is a pair of minimal $S$-linear $A_n$-structures on $E$, where $n\ge d$, 
such that $m_i=m'_i$ for $i\le d$, and such that
there exists a gauge transformation $f$ such that $f*m=m'$. Then there exists a 
gauge transformation $f'$, homotopic to $f$, such that $f'*m=m'$ and $f'_i=0$ for $i\le d-1$.
\end{lem}

\Pf . We are going to construct by induction on $j=1,\ldots,d-1$
a gauge transformation $f'$, homotopic to $f$, such that $f'*m=m'$ and
$f'_i=0$ for $i\le j$. The base case $j=1$ is clear. Assume that we have $f'$ homotopic to $f$
such that $f'*m=m'$ and $f'_i=0$ for $i<j$.
We have to show that $f'$ can be improved
to make $f'_i=0$ for $i\le j$ (provided $j\le d-1$). Indeed, by Lemma \ref{gauge-lem}, we have
$$0=m'_{j+1}-m_{j+1}=\pm \de(f'_j),$$ 
so 
$f'_j$ is a Hochschild cocycle giving a class in $HH^1(E/S)_{1-j}$. Since this class is zero, 
there exists $\phi'_{j-1}$ of degree $1-j$ such that
$f'_j=[m_2,\phi'_{j-1}]$. By \cite[Lem.\ 2.1]{P-ext}, this implies the existence of a gauge transformation $f''$, homotopic to $f'$, with $f''_i=0$ for $i\le j$.
\ed

\subsection{Moduli of $A_\infty$-structures}\label{ainf-moduli-sec}

%Discuss some kind of families of $A_\infty$-structures.

In this section we fix a finite-dimensional graded associative $k$-algebra $E$, where $k$ is a field.
We are going to study the moduli functors of minimal $A_n$-structures on $E$ extending the
product $m_2$ on $E$. We set $CH^i(E)_j=CH^i(E/k)_j$ and denote by
$\de^i_j:CH^i(E)_j\to CH^{i+1}(E)_j$ the Hochschild differential.

We will view all the spaces $CH^i(E)_j$ as affine spaces over $k$ and the group
$\fG$ of gauge transformations as an affine algebraic group over $k$, which is the projective
limit of the affine algebraic groups of finite type over $k$,
$$\fG[2,n]:=\fG/\fG_{\ge n+1}.$$
Note that each projection $\fG[2,n]\to \fG[2,n-1]$ has a natural section (not compatible with the group
structure). The kernel of this projection is isomorphic to the direct product of several copies of the
additive group over $k$. In particular, we see that the $k$-groups $\fG$ and $\fG[2,n]$ are smooth.

\begin{defi} (i) The affine $k$-scheme $\AA_n=\AA_{E,n}$
is defined as the closed subscheme of the affine space 
$CH^2(A)_{-1}\times\ldots\times CH^2(A)_{2-n}$, given by the 
$A_\infty$-constraints \eqref{a-inf-id-eq} (where $m_2$ is the given product on $E$).
For a $k$-algebra $R$, 
the set of points $\AA_n(R)$ is simply the set of minimal $R$-linear $A_n$-structures on $E\ot R$.

\noindent
(ii) For $n\ge 3$ we have an action of the algebraic $k$-group 
$\fG[2,n-1]$ of gauge transformations $(f_2,\ldots,f_{n-1})$
on $\AA_n$. We define $\MM_n=\MM_{E,n}$ to be the functor associating with a commutative $k$-algebra
$R$ the set of orbits for the corresponding action on the set of $R$-points,
$$\MM_n(R):=\AA_n(R)/\fG[2,n-1](R),$$
which is the set of gauge equivalence classes of minimal $R$-linear $A_n$-structures on $E\ot R$.
\end{defi}

For example, $\AA_3$ is the affine space $ZH^2(E)_{-1}$ of Hochschild cocycles and
the group $\fG[2,2]$ is just the product of additive groups corresponding to $CH^1(E)_{-1}$
acting on $\AA_3$ via the Hochschild differential $\de^1_{-1}:CH^1(E)_{-1}\to ZH^2(E)_{-1}$.
Thus, the functor
%if $HH^1_{-1}(E)=0$ then this action is free and 
$\MM_3$ is represented 
by the affine space $HH^2(E)_{-1}$.
This observation will be generalized to the moduli functors $\MM_n$ in Theorem
\ref{a-inf-moduli-thm} below (under some additional assumptions on $E$).

Below we work in the category of schemes over $k$.

\begin{defi} Let $G$ be a group scheme, $X$ a scheme with $G$-action.
We say that a closed subscheme $S\sub X$ is a {\it nice section for the action of $G$ on $X$}
if there exists a morphism $(\rho,\pi):X\to G\times S$, 
such that the morphism $\pi:X\to S$ is $G$-invariant, $\rho(S)=1$ and $\rho(x)\pi(x)=x$ (this is an equality of morphisms
$X\to X$).
\end{defi}

Note that for any scheme $T$ and any pair of points $x,x'\in X(T)$ one has $\pi(x)=\pi(x')$ if and only if
there exists a point $g\in G(T)$ such that $x'=gx$.
Thus, if $S$ is a nice section for $(G,X)$ then the functor $T\mapsto X(T)/G(T)$ is represented by $S$.

\begin{lem}\label{group-sections-lem}
Let $G$ be a group scheme acting on a scheme $X$. Assume that
$G$ fits into an exact sequence of groups
$$1\to H\to G\to G'\to 1$$
and that the projection $G\to G'$ admits a section $\si:G'\to G$ which is a morphism of schemes
(not necessarily compatible with the group structures).
Suppose we have a scheme $X'$ with an action of $G'$ and a morphism $f:X\to X'$
compatible with the $G$-action via the homomorphism $G\to G'$.
Assume that $S_H\sub X$ is a nice section for the $H$-action on $X$ and $S'\sub X'$ is a nice section
for the $G'$-action on $X'$. Finally assume that the following condition holds: for
any scheme $T$ and any points $x\in X(T)$, $g\in G(T)$ such that $f(gx)=f(x)$ there exists
a point $h\in H(T)$ such that $gx=hx$.
Then $S:=S_H\cap f^{-1}(S')$ is a nice section for the $G$-action on $X$.
\end{lem}

\Pf . Without loss of generality we can assume that $\si(1)=1$.
Let $(\rho_H,\pi_H):X\to H\times S_H$ and $(\rho',\pi'):X'\to G'\times S'$
be the morphisms that exist by the definition of a nice section.
Let us define a morphism $\rho_f:X\to G$ by
$$\rho_f(x)=\si(\rho'(f(x))).$$
Then $\rho_f(f^{-1}(S'))=1$, and we have
$$f(\rho_f(x)^{-1}x)=\rho'(f(x))^{-1}f(x)=\pi'(f(x)).$$
Thus, we obtain a morphism 
$$\pi_f:X\to f^{-1}(S'):x\mapsto \rho_f(x)^{-1}x,$$ 
satisfying $f\circ \pi_f=\pi'\circ f$ and $\pi_f|_{f^{-1}(S')}=\id$.
Since the morphism $f:X\to X'$ is $H$-invariant, the subscheme $f^{-1}(S')$ is
$H$-invariant. Hence, $\pi_H$ sends $f^{-1}(S')$ to $S=S_H\cap f^{-1}(S')$. Thus, we can define a morphism
$$\pi=\pi_H\circ\pi_f:X\to S.$$
By definition, we have
$$x=\rho_f(x)\pi_f(x)=\rho(x)\cdot \pi(x), \ \text{ where}$$
$$\rho(x)=\rho_f(x)\rho_H(\pi_f(x)).$$
Note that for $x\in S$ we have $\pi_f(x)=x$, so $\rho_H(\pi_f(x))=1$ and hence $\rho(x)=1$.
It remains to check that $\pi$ is $G$-invariant.
Let $x\in X(T)$, $g\in G(T)$. Then we have
$$f(\pi_f(gx))=\pi'(f(gx))=\pi'(f(x))=f(\pi_f(x)).$$
Also, $\pi_f(gx)$ and $\pi_f(x)$ are in the same $G(T)$-orbit:
$$\pi_f(gx)=\rho_f(gx)^{-1}gx=\rho_f(gx)^{-1}g\rho_f(x)x.$$
Hence, by our assumption, there exists $h\in H(T)$ such that
$\pi_f(gx)=h\pi_f(x)$, so that
$$\pi(gx)=\pi_H(\pi_f(gx))=\pi_H(h\pi_f(x))=\pi_H(\pi_f(x))=\pi(x).$$
\ed

\begin{thm}\label{a-inf-moduli-thm} 
Assume that $HH^1(E)_{-i}=0$ for $1\le i\le n-3$.
Then the action of $\fG[2,n-1]$ on $\AA_n$ admits a nice section
$S_n\sub\AA_n$. Hence,
the affine scheme $S_n$ represents the functor $\MM_n$ of gauge equivalence classes
of minimal $A_n$-structures on $E$.
\end{thm}

\Pf . We use induction on $n$. For $n=3$ we define the subscheme
$$S_3\sub ZH^2(E)_{-1}=\AA_3$$
to be any subspace complementary to $\de^1_{-1}\bigl(CH^1(E)_{-1}\bigr)$. This gives a nice section for the 
action of $\fG[2,2]=CH^1(E)_{-1}$ on $\AA_3$. Suppose we already have 
the subscheme $S_{n-1}\sub\AA_{n-1}$
which is a nice section for the action of $\fG[2,n-2]$ on $\AA_{n-1}$.
%Let $\wt{S}_n\sub \AA_n(E)$ be the preimage of $S_{n-1}\sub \AA_{n-1}(E)$ under
%the projection $\AA_n(E)\to \AA_{n-1}(E)$.
Let us consider the natural projection $\pi_n:\AA_n\to \AA_{n-1}$. It is compatible with
the homomorphism of groups
$\fG[2,n-1]\to \fG[2,n-2]$ that fits into the exact sequence of groups
\begin{equation}\label{gauge-groups-sequence-eq}
0\to CH^1(E)_{2-n}\to \fG[2,n-1]\to \fG[2,n-2] \to 0.
\end{equation}
Note that the action of $x\in CH^1(E)_{2-n}$ on $(m_2,\ldots,m_n)\in\AA_n$ changes $m_n$ to $m_n+\de^1(x)$ and
does not change $(m_2,\ldots,m_{n-1})$. Thus, if $K_{2-n}\sub CH^2_{2-n}$ is any complement to the subspace
$\im \de^1_{2-n}$, then the condition $m_n\in K_{2-n}$ defines a closed subscheme $\AA'_n\sub \AA_n$,
which is a section for the $CH^1(E)_{2-n}$-action on $\AA_n$. 
Furthermore, it is easy to see that it is a nice section: to define the corresponding map $\rho:\AA_n\to CH^1(E)_{2-n}$
we use any linear map $Q^2:CH^2(E)_{2-n}/K_{2-n}\to CH^1(E)_{2-n}$, such that $\de^1Q^2(x)-x\in K_{2-n}$, and
set $\rho(m_2,\ldots,m_n)=Q^2(m_n)$.

Now we would like to apply Lemma \ref{group-sections-lem}
to the group $G=\fG[2,n-1]$ acting on $X=\AA_n$, to the exact sequence
\eqref{gauge-groups-sequence-eq} and to the projection $X\to X'=\AA_{n-1}$.
Note that the assumptions of this Lemma are satisfied by Lemma \ref{important-homotopy-lem} (which we apply with $d=n-1$).
Thus, we conclude that
$$S_n:=\AA'_n\cap \pi_n^{-1}(S_{n-1})$$ 
is a nice section for the action of $\fG[2,n-1]$ on $\AA_n$.
\ed

\begin{cor}\label{representability-cor}
Assume that $HH^1(E)_{<0}=0$. Then the functor $\MM_\infty$ associating with a $k$-algebra $R$ 
the set of gauge equivalence classes
of minimal $A_\infty$-structures on $E\ot R$ is represented by the
affine scheme $S_\infty$, the inverse limit of the affine schemes $S_n$ representing $\MM_n$.
\end{cor}

\Pf . It is enough to show that for every $k$-algebra $R$, the natural map
$$\MM_\infty(R)\to \liminv_n \MM_n(R)$$
%=\AA_n(R)/\fG[2,n-1](R)$.
is an isomorphism.
%Since the set of minimal $A_\infty$-structures on $E\ot R$ is $\underset{n}{\varprojlim} \AA_n(R)$,
%and $\fG(R)=\underset{n}{\varprojlim} \fG[2,n](R)$, this follows from surjectivity of the projections
%$\fG[2,n](R)\to \fG[2,n-1](R)$.

To prove the surjectivity, suppose we have a collection $(\a_n)_{n\ge 3}$ 
of minimal $A_n$-structures on $E\ot R$, and a set of gauge equivalences
$(u_n\in \fG[2,n-1])$ such that $(\a_{n+1})_{\le n}=u_n\cdot \a_n$. Then we can recursively construct minimal $A_n$-structures
$(\a'_n)$, such that $(\a'_{n+1})_{\le n}=\a'_n$, and gauge equivalences $(v_n\in\fG[2,n-1])$ such that
$\a'_n=v_n\cdot \a_n$ and $(v_{n+1})_{\le n-1}u_n=v_n$. 
Namely, if $(\a'_i)$, $(v_i)$ for $i\le n$ are already constructed, then we pick a gauge equivalence $v_{n+1}\in \fG[2,n]$,
such that $(v_{n+1})_{\le n-1}=v_nu_n^{-1}$, and set $\a'_{n+1}:=v_{n+1}\cdot \a_{n+1}$. Then $(\a'_n)$ defines
the required minimal $A_\infty$-structure.

For the injectivity, suppose $\a$ and $\b$ are minimal $A_\infty$-structures such that $\a_{\le n}$ is gauge equivalent
to $\b_{\le n}$ for each $n$. We are going to construct recursively
a sequence of gauge equivalences $u_2=\id, u_3, u_4,\ldots$,
such that $(u_n)_{\le n-2}=\id$ and 
$$\a_{\le n}=(u_nu_{n-1}\ldots u_2)\b_{\le n}.$$
Indeed, suppose $(u_i)$ for $i\le n-1$ are already constructed and satisfy the above property.
The $A_{n+1}$-structures $\a_{\le n+1}$ and $(u_nu_{n-1}\ldots u_2)\b_{\le n+1}$ agree up to $n$,
and are gauge equivalent. Hence, by Lemma \ref{important-homotopy-lem}, there exists a gauge equivalence $u_{n+1}$
such that $(u_n)_{\le n-1}=\id$ and 
$$\a_{\le n+1}=(u_{n+1}u_nu_{n-1}\ldots u_2)\b_{\le n+1}.$$
It remains to note that the infinite product $\ldots u_4u_3u_2$ converges in $\fG$ to some element $u$, such that
$\a=u\b$. 

\ed

An important particular case is when we have some additional vanishing for components
of $HH^2(E)$ and $HH^3(E)$.

\begin{cor}\label{finite-a-inf-moduli-cor} 
Assume that $HH^1(E)_{<0}=HH^2(E)_{<2-n}=0$.
Then the natural morphism $\MM_\infty\to \MM_n$ is a closed embedding.
Hence, in this case $\MM_\infty$ is represented by an affine scheme of finite type over $k$.
If in addition $HH^3(E)_{<2-n}=0$ then $\MM_\infty\to\MM_n$ is an isomorphism.
\end{cor}

\Pf . We use the notation from the proof of Theorem \ref{a-inf-moduli-thm}. By Corollary \ref{representability-cor},
the vanishing of $HH^1(E)_{<0}$ guarantees that $\MM_\infty$ is the inverse limit of the affine schemes $S_n$.
If $HH^2(E)_{2-i}$=0 then the projection $\AA'_i\to \AA_{i-1,i}$ is an isomorphism, hence, $S_i\to S_{i-1}$ is a closed
embedding. Applying this to all $i>n$ we get the first assertion. If $HH^2(E)_{2-i}=HH^3(E)_{<2-i}=0$ then
$\AA_{i-1,i}=\AA_{i-1}$, the projection $\AA'_i\to \AA_{i-1}$ is an isomorphism, hence, $S_i\to S_{i-1}$ is an isomorphism.
Using this for all $i>n$ gives the last assertion.
\ed

\subsection{The map from the moduli space of non-special curves to the moduli space of $A_\infty$-structures}
\label{map-sec}

%We are going to study moduli of minimal $A_\infty$-structures for the algebra $E_g\ot k$ 
%which we denote from now on simply as $E_g$. We denote by $\MM_\infty$ the
Let us denote by $\MM_\infty=\MM_{\infty,\Z}$
the functor of minimal $A_\infty$-structures on $E_{g,\Z}$ (the natural $\Z$-form of the algebra $E_g$)
up to a gauge equivalence.
We are going to relate it to the moduli of non-special curves.

Let $S$ be a commutative algebra.
Given a family of curves in $\wt{\UU}^{ns,a}_{g,g}(S)$, the construction of Section \ref{Cech-ainf-sec} gives
a minimal $S$-linear $A_\infty$-structure on
$E_{g,\Z}\ot S$, defined uniquely up to gauge equivalence (recall that the construction depends
on a choice of formal parameters at marked points).
This gives a morphism of functors on commutative algebras
\begin{equation}\label{a-inf-map}
a_\infty:\wt{\UU}^{ns,a}_{g,g}\to \MM_\infty.
\end{equation}
%sending $(C,p_1,\ldots,p_g)$ to the equivalence class of $A_{\infty}$-structures on $E_g$,
%such that ???
Later we will show that the restriction of the functor $\MM_\infty$ to commutative algebras over a field $k$ is actually represented by a $k$-scheme, so from \eqref{a-inf-map}
we will get a morphism of $k$-schemes.

Recall that $\G_m^g$ acts on $\wt{\UU}^{ns,a}_{g,g}$ by rescaling the tangent vectors at the marked point
by $v_i\mapsto \la_i^{-1}v_i$.
The morphism \eqref{a-inf-map} is compatible with the $\G_m^g$-action, where the action on 
$\MM_\infty$ is induced by the $\G_m^g$-action on $E_g$ given on generators by
\begin{equation}\label{la-A-B-eq}
\la: A_i\mapsto A_i, \ B_i\mapsto \la_iB_i.
\end{equation}

We are going to show how to recover a curve (and marked points) from the corresponding
$A_\infty$-structure.
Namely, for a non-special curve $(C,p_1,\ldots,p_g)$ over $\Spec(S)$,
let $\Per(C,p_1,\ldots,p_g)\sub D^b(C)$ denote the thick triangulated subcategory generated
by the objects $\OO_C, \OO_{p_1}, \ldots, \OO_{p_g}$, or equivalently, by 
their direct sum, $G$.
Then the functor $R\Hom(G,?)$ induces an equivalence
\begin{equation}\label{perfect-equivalence-eq}
\Per(C,p_1,\ldots,p_g)\to \Per(E_{g,\infty}),
\end{equation}
where $E_{g,\infty}$ is the (minimal) $A_\infty$-algebra of endomorphisms of $G$ in an $A_\infty$-enhancement
of $D^b(C)$ (so $E_{g,\infty}=E_{g,\Z}\ot S$ as an associative $S$-algebra) and 
$\Per(E_{g,\infty})$ is the derived category of perfect right modules over $E_{g,\infty}$.

\begin{prop}\label{set-theoretic-inj-prop} The map $a_\infty$ is injective.
\end{prop}

\Pf . Suppose we have a gauge equivalence between the $A_{\infty}$-structures on $E_{g,\Z}\ot R$
associated with the curves $(C,p_1,\ldots,p_g)$ and $(C',p'_1,\ldots,p'_g)$ over $\Spec(R)$.
Let us denote the corresponding $A_\infty$-algebras $E_{g,\infty}$ and $E'_{g,\infty}$. Note
that with each idempotent in $E_{g,\Z}$ we have the associated $A_\infty$-module over $E_{g,\infty}$ 
(resp., $E'_{g,\infty}$) that is the image of one of the objects $\OO_C, \OO_{p_1},\ldots,\OO_{p_g}$
under the equivalence \eqref{perfect-equivalence-eq}
(resp., $\OO_{C'},\OO_{p'_1},\ldots,\OO_{p'_g}$). The gauge equivalence gives a quasi-isomorphism 
$$E_{g,\infty}\rTo{\sim} E'_{g,\infty}$$ 
of $A_\infty$-algebras over $S$, which induces an
equivalence between the derived categories of perfect right modules 
$$\Per(E'_{g,\infty})\rTo{\sim} \Per(E_{g,\infty}).$$
Using \eqref{perfect-equivalence-eq}
we get an equivalence 
$$\Phi:\Per(C',p'_1,\ldots,p'_g)\rTo{\sim} \Per(C,p_1,\ldots,p_g)$$ 
sending $\OO_{p'_i}$ to $\OO_{p_i}$, $\OO_{C'}$ to
$\OO_C$, and inducing the identity maps on
$$\Hom(\OO_{C'},\OO_{p'_i})=\Hom(\OO_C,\OO_{p_i})=S \ \text{ and}$$
$$\Ext^1(\OO_{p'_i},\OO_{C'})=\Ext^1(\OO_{p_i},\OO_C)=S$$
(the latter identifications are given by the tangent vectors at the marked points).
Let $D=p_1+\ldots+p_g$ (resp., $D'=p'_1+\ldots+p'_g$).
Note that the line bundles $\OO(nD)$ belong to $\Per(C,p_1,\ldots,p_g)$. More precisely, we have
exact sequences
$$0\to \OO_C(nD)\to\OO_C((n+1)D)\to \bigoplus_{i=1}^g\OO_{p_i}\to 0$$
and similarly for $C'$, which give isomorphisms
$\Phi\bigl(\OO_{C'}(nD')\bigr)\simeq \OO_C(nD)$ compatible with these exact sequences.
Hence, $\Phi$ gives an isomorphism of graded algebras
$$\bigoplus_n H^0\bigl(C',\OO(nD')\bigr)\rTo{\sim} \bigoplus_n H^0\bigl(C,\OO(nD)\bigr),$$
inducing the identity on 
$$H^0\bigl(C',\OO(nD')\bigr)/H^0\bigl(C',\OO((n-1)D')\bigr)=S^g=H^0\bigl(C,\OO(nD)\bigr)/H^0\bigl(C,\OO((n-1)D)\bigr)$$
for $n\ge 2$.
Since the divisor $D$ on $C$ (resp., $D'$ on $C'$) is ample, passing to $\Proj$ we get
an isomorphism $C\rTo{\sim} C'$, sending $D$ to $D'$ and compatible with the given
trivializations of $\OO_C(D)/\OO_C$ (resp., $\OO_{C'}(D')/\OO_{C'}$).
\ed

\begin{rems}\label{ainf-map-rem}
1. When the curve $C$ is smooth (where $S=k$, a field, and $g\ge 2$),
it can be recovered just  from the {\it derived Morita equivalence class} 
of the $A_{\infty}$-structure on $E_g$ associated with $(C,p_1,\ldots,p_g)$
for any non-special divisor $p_1+\ldots+p_g$.
Namely, this follows from the fact that $C$ can be recovered from the equivalence class of its derived category
$D^b(C)$ (see \cite{BO}). There is an extension of this result to singular $C$ provided $C$ is
Gorenstein and $\om_C$ is ample (see \cite{Ballard}). But of course,
the gauge equivalence class of an $A_\infty$-structure used in Proposition \ref{set-theoretic-inj-prop} above
gives more information than the derived Morita equivalence class.

\noindent
2. The map \eqref{a-inf-map} extends naturally to the stack of curves with a non-special divisor
$p_1+\ldots+p_g$ which is not necessarily ample. The extended map is not injective: it does not distinguish some curves that become isomorphic after contracting all irreducible components not containing any marked points. For example, suppose $f:C\to\ov{C}$ a contraction of a tree $T$ of $\P^1$'s that is attached transversally
to one of the components of $C$. Then $Rf_*\OO_C=\OO_{\ov{C}}$.
Hence, the pull-back functor $f^*:\Per(\ov{C})\to \Per(C)$ is fully faithful.
Since this functor lifts to a dg-level, it induces an equivalence between the $A_\infty$-algebra of
endomorphisms of $G$ and of $f^*G$ for any $E\in\Per(\ov{C})$. 
This implies that if $T$ does not contain any marked points then the $A_\infty$-structure
on $E_g$ associated with $(C,p_1,\ldots,p_n)$ is gauge equivalent to the one associated with 
$(\ov{C},\ov{p}_1,\ldots,\ov{p}_n)$ (where $\ov{p}_i=f(p_i)$).
\end{rems}

\subsection{Hochschild cohomology of the cuspidal curve of genus $g$}\label{cusp-sec}

Starting from this section we work over a field $k$.

Recall that our moduli space $\wt{\UU}^{ns,a}_{g,g}$ has a special point $[C^{\cusp}_g]$ corresponding to the cuspidal
curve of genus $g$ (see Definition \ref{cusp-curve-def}).

The cuspidal curve $C^{\cusp}_g$ is equipped with the natural action of $\G_m^g$ such that
for generators $f_i, h_i\in H^0(C^{\cusp}_g\setminus D,\OO)$ we have
$$(\la^{-1})^*f_i=\la_i^2f_i, \ \ (\la^{-1})^*h_i=\la_i^3h_i$$
for $\la=(\la_1,\ldots,\la_g)\in\G_m^g$. 
%(this action is the opposite of the one induced by the action on
%the universal curve considered in Sections \ref{marked-moduli-sec} and \ref{homotopy-sec}???).
The marked points $p_1,\ldots,p_g\in C^{\cusp}_g$ are invariant under this action and $\la$ acts on the
tangent space at $p_i$ as rescaling by $\la_i^{-1}$.

The induced $\G_m$-action coming from the diagonal $\G_m\sub\G_m^g$
will be especially important for us. 
Recall that the $\G_m$-orbit of any point of $\wt{\UU}^{ns,a}_{g,g}$ 
contains the cuspidal curve in its closure, so it is important to understand deformations of
$C^{\cusp}_g$ and of its derived category. 

%Recall that we have a $\G_m^g$-action on $C^{\cusp}_g$. Consider the diagonal copy of $\G_m$.
%So the action sends $f_i$ to $\la^2f_i$ and $h_i$ to $\la^3f_i$ (on affine part; to extend it to
%projective curve, consider the same action of $F_i$, $F_j$, fixing $T$). 
%Hence, the induced $\G_m$-action on the algebra $\Ext^*(G,G)$, where $G=???$
%acts trivially on $\Ext^0$ and rescales by $\la^{\pm 1}$ the $\Ext^1$.

%For our conventions on the grading on Hochschild cohomology, see ???

\begin{prop}\label{cusp-equiv-prop} 
The functor
$$\Per(C^{\cusp}_g)\to \Per(E_g)$$
associated with the generator $\OO\oplus \OO_{p_1}\oplus\ldots\oplus\OO_{p_g}$
is an equivalence of categories.
Hence, it induces a $\G_m^g$-equivariant isomorphism $HH^*(C^{\cusp}_g)\simeq HH^*(E_g)$.
The second grading on $HH^*(E_g)$
corresponds to the weights of the $\G_m$-action on it.
\end{prop}

\Pf . We claim that the gauge equivalence class of minimal $A_\infty$-structures
 on $E_g$ coming from $C^{\cusp}_g$ contains the structure with $m_i=0$ for $i>2$.
Indeed, we can choose
$\G_m$-equivariant formal parameters at $p_i$'s similarly to the construction of Section \ref{Cech-ainf-sec}.
Then the vanishing of the corresponding products $m_i$ for $i>2$ follows as in the proof of Proposition \ref{weight-prop}. Namely, the $\G_m$-weights on $E_g$ coincide with the cohomological grading
(recall that the $\G_m$-action on $C^{\cusp}_g$ is the opposite of the one used
in Proposition \ref{weight-prop}). The operation $m_i$ lowers this grading by $i-2$. On the other hand,
by the construction, all higher products have weight $0$. This implies that $m_i=0$ for $i>2$. 
The identification of the $\G_m$-weights on $E_g$ with the cohomological grading 
also implies that the second grading
on $HH^*(E_g)$ is given by the $\G_m$-weights. 
%hence it corresponds to the $\G_m$-weights
%on $HH^*(C^{\cusp}_g)$.
\ed

We also can look at the induced $\G_m$-action on other invariants of $C^{\cusp}_g$.

\begin{lem}\label{tangent-action-lem} Assume that the $\cha(k)\neq 2$, and either $g\ge 2$ or $\cha(k)\neq 3$.

\noindent
(i) The $g$-dimensional space $H^1(C^{\cusp}_g,\OO)$ has weight $1$ with respect to the $\G_m$-action.

\noindent
(ii) The space $H^0(C^{\cusp}_g,\TT)$, where $\TT$ is the tangent sheaf, decomposes as a direct sum
$$H^0(C^{\cusp}_g,\TT)=H^0\bigl(C^{\cusp}_g,\TT(-D)\bigr)\oplus V,$$
where $V$ is a $g$-dimensional subspace of weight $1$ with respect to the $\G_m$-action such that
the composition
$$V\to H^0(C^{\cusp}_g,\TT)\to H^0(D,\TT|_D)$$
is an isomorphism. Furthermore,
$$H^0\bigl(C^{\cusp}_g,\TT(-D)\bigr)=H^0(C^{\cusp}_g,\TT)^{\G_m}$$ 
and this space is spanned by $g$ linearly independent derivations coming from the $\G_m^g$-action on
$C^{\cusp}_g$. 

\noindent
(iii) One has $H^0\bigl(C^{\cusp}_g,\TT(-2D)\bigr)=0$.

\noindent
(iv) One has $H^1(C^{\cusp}_g,\TT)=0$.

\noindent
(v) The natural map $H^0\bigl(C^{\cusp}_g,\TT(nD)\bigr)\to H^0\bigl(C^{\cusp}_g,\TT(nD)|_D\bigr)$ is surjective for $n\ge -1$ and is an isomorphism
for $n=-1$.
\end{lem}

\Pf . (i) Let us write for brevity $C=C^{\cusp}_g$. 
Set $U=C\setminus D$, $V=C\setminus q$, where $q$ is the singular point of $C$.
Then $(U,V)$ is an affine covering of $C$. Furthermore, $V=\sqcup_{i=1}^g V_i$, where
$V_i\simeq\A^1$ and $U\cap V_i=V_i\setminus p_i$. Let $(C_i)$ be irreducible
components of $C$, where $C_i=\{q\}\cup V_i$. We can view $C_i$ as $\P^1$ with the point
$0\in\A^1\sub\P^1$ pinched: if $x_i$ is the coordinate on $\A^1$ then the algebra of
functions on $C_i\setminus p_i$
is the subring $k[x_i^2,x_i^3]\sub k[x_i]$. Hence, we can consider $x_i^{-1}$ as a coordinate on
$V_i=C_i\setminus q$. Thus, functions on $U\cap V$ are collections of Laurent polynomials
$(P_i(x_i,x_i^{-1}))_{i=1,\ldots,g}$. Such a function extends to $V$ if and only if $P_i\in k[x_i^{-1}]$.
On the other hand, for $P_i\in x_ik[x_i]$ the collection $(P_i)$ extends to $U$ if and only if
$P_i$ has no linear term in $x_i$. Thus, $H^1(C,\OO)$, that can be identified with the cokernel of the map
$$H^0(U,\OO)\oplus H^0(V,\OO)\to H^0(U\cap V,\OO),$$
is represented by the functions $(a_i x_i)$ for $a_i\in k$. It remains to check that $(\la^{-1})^*x_i=\la x_i$ for
$\la\in \G_m$. Indeed, we have $x_i=t_i^{-1}$ where $t_i$ is the canonical parameter near $p_i$, 
so the assertion follows from the fact that $\la^*t_i=\la t_i$.

\noindent
(ii) Since $C$ is smooth near $D$, the restriction map $H^0(C,\TT)\to H^0(U,\TT)$ is injective,
and its image consists of vector fields that extend regularly to $D$. We have further embedding
$H^0(U,\TT)\to \prod_{i=1}^g H^0(U_i,\TT)$, where $U_i=C_i\setminus p_i$. Thus, 
using the coordinate $x_i$ on the normalization of $U_i$ as before, each derivation $v$ on $C$ gives
rise to derivations $v_i$ of $k[x_i^2,x_i^3]$ of the form $v_i=P_i(x_i,x_i^{-1})\partial_{x_i}$. 
Assume first that $g\ge 2$. 
Applying our vector field
$v$ to the elements $x_i^2$, complemented by zeros on the other components,
 and using the assumption $\cha(k)\neq 2$, we see that $P_i\in k[x_i]$.
Furthermore, if $P_i\equiv a_i\mod x_ik[x_i]$ for $i=1,\ldots,n$, 
then applying $v$ to the same elements $x_i^2$ 
twice we derive that $a_i=0$ for each $i$, i.e., $P_i\in x_ik[x_i]$.
Since $x_i^n\frac{\partial}{\partial x_i}$ has a pole at infinity of order $n-2$ for $n>2$, we see that
in order for $v$ to define a global section of $\TT_C$
each $P_i$ should be a linear combination of $x_i$ and $x_i^2$, 
%which implies that the weights of $\G_m$ on $H^0(C,\TT_C)$
%are $0$ and $1$.
Conversely, every such collection $(v_i)$
defines a derivation of $\OO(U)$, regular at infinity. Now we note that
$v\in H^0\bigl(C,\TT(-D)\bigr)$ if and only if $b_i=0$ and let $V$ be a subspace consisting of 
$(v_i=b_ix_i^2\frac{\partial}{\partial x_i})$.

In the case $g=1$ the argument is similar but we have to use the assumption that $\cha(k)\neq 3$ to 
derive that $P_1\in x_1k[x_1]$ (e.g., using the identity $3x_1v(x_1^2)=2v_1(x_1^3)$).

\noindent
(iii) This follows immediately from (ii), since the derivations spanning $H^0\bigl(C,\TT(-D)\bigr)$ project to
a basis in $H^0\bigl(C,\TT(-D)|_D\bigr)$.

\noindent
(iv) The proof is similar to that of \cite[Lem.\ 4.16]{LP2}. Consider the affine covering 
$(U,V)$ of $C$ defined above. Then $H^0(U\cap V,\TT)$ consists of vector fields
$(v_i)$ of the form $v_i=P_i(x_i,x_i^{-1})\frac{\partial}{\partial x_i}$. Such a vector field extends
to a derivation over $V$ provided $P_i$ is a linear combination of $x_i^n$ with $n\le 2$. On the other hand,
if all $P_i\in x_ik[x_i]$ then $(v_i)$ extends to a derivation over $U$. 

\noindent
(v) As we have seen in the proof of (ii), sections of $\TT$ on $U$ with poles of order at most $n$
at $p_1,\ldots,p_g$ correspond to derivations $v_i$ of the form $v_i=x_iP_i\frac{\partial}{\partial t_i}$,
where $P_i$ is a polynomial of degree $n+1$, and the assertion follows from the
fact such sections generate $\TT(np_1)|_{p_1}\oplus\ldots\oplus \TT(np_g)|_{p_g}$ for $n\ge -1$.
\ed

\begin{lem}\label{HH-C-lem} 
(i) Let $C$ be a quasiprojective curve, $q\in C$ a point such that $C\setminus q$ is smooth, $U\sub C$
an open affine neighborhood of $q$. Then one has natural exact sequences
$$0\to H^1(C,\OO)\to HH^1(C)\to H^0(C,\TT)\to 0,$$
$$0\to H^1(C,\TT)\to HH^2(C)\to HH^2(U)\to 0.$$

\noindent
(ii) Assume the characteristic of $k$ is not $2$ or $3$.
Let $U=C^{\cusp}_g\setminus D$. Then the natural map
$$HH^2(C^{\cusp}_g)\to HH^2(U)$$
is an isomorphism.
Also, $HH^1(C^{\cusp}_g)_{<0}=0$.
\end{lem}

\Pf . (i) This is proved in \cite[Sec.\ 4.1.3]{LP2}.

\noindent
(ii) By Lemma \ref{tangent-action-lem}(iv), we have $H^1(C^{\cusp}_g,\TT)=0$, so the first assertion follows from the
second exact sequence in (i). On the other hand, 
the first exact sequence in (i) together with Lemma \ref{tangent-action-lem}(i)(ii) imply that
the only weights of $\G_m$ that occur on $HH^1(C^{\cusp}_g)$ are $0$ and $1$.
\ed

\begin{cor}\label{HH1-cor} Assume that $\cha(k)\neq 2$ or $3$.
Then one has $HH^1(E_g)_{<0}=0$.
\end{cor}

\Pf . Combine the vanishing of $HH^1(C^{\cusp}_g)_{<0}$ with Proposition \ref{cusp-equiv-prop}.
\ed

For a $k$-scheme $X$ we denote by $\bL_X$ the cotangent complex of $X$ over $k$
(see \cite{Illusie}).

\begin{lem}\label{cotangent-complex-lem} Assume that $\cha(k)\neq 2$ or $3$.
Let $C=C^{\cusp}_g$. 

\noindent
(i) The natural maps
$$\Ext^i\bigl(\bL_C,\OO_C(-2D)\bigr)\to \Ext^i\bigl(\bL_C,\OO_C(-D)\bigr)\to \Ext^i(\bL_C,\OO_C)$$
are isomorphisms for $i=1,2$, while the natural maps
\begin{equation}\label{Hom-L-C-D-map}
\Hom\bigl(\bL_C,\OO_C(-D)\bigr)\to \Hom\bigl(\bL_C,\OO_D(-D)\bigr) \ \text{ and}
\end{equation}
\begin{equation}\label{Hom-L-C-D-map2}
\Hom(\bL_C,\OO_C)\to \Hom(\bL_C,\OO_D)
\end{equation}
are surjective. Also, one has 
$$\dim\Hom(\bL_C,\OO_C)=2g, \ \ \dim\Hom\bigl(\bL_C,\OO_C(-D)\bigr)=g,$$
$$\Hom\bigl(\bL_C,\OO_C(-2D)\bigr)=0.$$

\noindent
(ii) Let $U=C\setminus D$. Then the natural map
$$\Ext^i\bigl(\bL_C,\OO_C(-2D)\bigr)\to \Ext^i_U(\bL_U,\OO_U)$$
is an isomorphisms for $i=1,2$.
\end{lem}

\Pf . (i) First, we observe that since $\bL$ is a locally free sheaf on the smooth part of $C$, it follows that
$$\Ext^{>0}(\bL_C,\OO_D)=0.$$
Thus, applying the functor $\Ext^\bullet(\bL_C,?)$ to the exact sequences 
$$0\to \OO_C(-2D)\to \OO_C(-D)\to \OO_D\to 0,$$
$$0\to \OO_C(-D)\to \OO_C\to \OO_D\to 0,$$
we get the required isomorphisms for $i=2$. 
Since $\bL_C\in D^{\le 0}(C)$ and $\und{H}^0(\bL_C)=\Om_C$, the maps
\eqref{Hom-L-C-D-map} and \eqref{Hom-L-C-D-map2}
can be identified with the maps
$$H^0\bigl(C,\TT(-D)\bigr)\to H^0\bigl(D,\TT(-D)|_D\bigr) \ \text{ and}$$
$$H^0(C,\TT)\to H^0(D,\TT|_D)$$
which are surjective by
Lemma \ref{tangent-action-lem}(v) (and the first is an isomorphism).
This implies the required isomorphisms for $i=1$. The computation of dimensions also follows from
Lemma \ref{tangent-action-lem}.

\noindent
(ii) Recall that the natural map $H^0\bigl(C,\TT(nD)\bigr)\to H^0\bigl(C,\TT(nD)|_D\bigr)$ is surjective for $n\ge -1$
by Lemma \ref{tangent-action-lem}(v).
Using the exact sequences
$$0\to \OO_C(nD)\to \OO_C\bigl((n+1)D\bigr)\to \OO_D\to 0$$
we deduce that the morphisms 
$\Ext^i\bigl(\bL_C,\OO_C(nD)\bigr)\to \Ext^i\bigl(\bL_C,\OO_C((n+1)D)\bigr)$
are isomorphisms for $n\ge -2$. 
Finally, since $\bL_C$ is a perfect complex, we have
$$\Ext^i_U(\bL_U,\OO_U)\simeq \Ext^i_C(\bL_C,\OO_U)\simeq \varinjlim_n \Ext^i\bigl(\bL_C,\OO_C(nD)\bigr),$$
and our assertion follows.
\ed

We need the following technical lemma only for the cuspidal curve $C^{\cusp}_g$ over $k$. 
However, we formulate it in a slightly bigger generality, since this could be of independent interest. 
Let $\bL_{X/R}$ denote the cotangent complex of a scheme $X$ over $\Spec(R)$.

\begin{lem}\label{curve-Hochschild-lem}  
Let $R$ be a Noetherian ring, $C$ a proper flat $R$-scheme of relative dimension $1$ 
such that the natural map $R\to H^0(C,\OO)$ is an isomorphism, $D\sub C$ a relative effective divisor such that
the projection $\pi:C\to \Spec(R)$ is smooth along $D$, and such that the scheme $U=C\setminus D$
is affine.
Assume also that the fibers of $\pi:C\to \Spec(R)$ are generically smooth.
Then the natural map
\begin{equation}\label{cot-HH-map}
\Ext^i(\bL_{U/R},\OO_U)\to HH^{i+1}(U/R)
\end{equation}
is an isomorphism for $i=1$, and is an injection for $i=2$.
\end{lem}

\Pf . By \cite[Thm.\ 8.1]{Quillen}, there is a spectral sequence
$$E_2^{pq}=\Ext^p({\bigwedge}^q\bL_{U/R},\OO_U)\implies HH^{p+q}(U/R)$$
where $\bigwedge^\bullet(?)$ denotes the exterior power functor on perfect complexes.
Since $\bL_{U/R}\in D^{\le 0}$, it follows that $\bigwedge^i\bL_{U/R}\in D^{\le 0}$, so $E_2^{pq}\neq 0$
only for $p\ge 0$ (and $q\ge 0$). Since $U$ is affine, we also have $E_2^{p0}=0$ for $p>0$.
We claim also that $E_2^{0q}=0$ for $q>1$. Indeed, we have
$$\Hom({\bigwedge}^q\bL_{U/R},\OO_U)=\Hom\bigl(\und{H}^0({\bigwedge}^q\bL_{U/R}),\OO_U\bigr).$$
Let $Z$ be the support of the coherent sheaf $\und{H}^0(\bigwedge^q\bL_{C/R})$ on $C$.
Then $Z$ is contained in the singular locus of the morphism $\pi$ (since $q>1$), so it is finite over $\Spec(R)$
and is contained in $U$. Therefore, we have 
$$\Hom\bigl(\und{H}^0({\bigwedge}^q\bL_{U/R}),\OO_U\bigr)\simeq \Hom\bigl(\und{H}^0({\bigwedge}^q\bL_{C/R}),\OO_C\bigr).$$
The image of a nonzero morphism $\und{H}^0(\bigwedge^q\bL_{C/R})\to \OO_C$
would be a nonzero subsheaf $\FF\hra\OO$ supported on $Z$.
But then we would have an injection $H^0(C,\FF)\to H^0(C,\OO)$. Since all global functions on $C$ are
constant on the fibers of $\pi$, it follows that $H^0(C,\FF)=0$. Since $Z$ is affine, we get that $\FF=0$ and our claim follows.
%Since $\Ext^2(\OO,\OO)=H^2(C,\OO)=0$, we get the required formula for $HH^2(C)$.

Thus, the spectral sequence implies that the map \eqref{cot-HH-map} is an isomorphism for $i=1$, while
for $i=2$ it fits into an exact sequence
$$0\to \Ext^2(\bL_{U/R},\OO_U)\to HH^3(U/R)\to \Ext^1({\bigwedge}^2\bL_{U/R},\OO_U)\to 0.$$
\ed

\begin{rem} Using the main result of \cite{BF}, one can similarly deduce that in the situation of Lemma
\ref{curve-Hochschild-lem}, assuming in addition that we work in characteristic zero, one has
$$HH^3(C/R)\simeq\Ext^2(\bL_{C/R},\OO)\oplus \Ext^1({\bigwedge}^2\bL_{C/R},\OO)$$
$$HH^2(C/R)\simeq\Ext^1(\bL_{C/R},\OO),$$
$$HH^1(C/R)\simeq H^1(C,\OO)\oplus \Hom(\bL_{C/R},\OO).$$
\end{rem}

%We deduce the following numerical result.
%
%\begin{lem} For $C=C^{\cusp}_g$ one has 
%$\dim HH^2(C)=\dim 

\subsection{Deformation theory}\label{def-theory-sec}

We refer to \cite{Manetti} for the basic notions of deformation theory used below.

Let us define two deformation functors
$$F_{g,g}, F_\infty: \Art_k\to \Sets$$
from the category of local Artinian $k$-algebras with the residue field $k$ to the category of sets.
The functor $F_{g,g}$ describes deformations of the cuspidal curve $C^{\cusp}_g$
viewed as an object of $\wt{U}^{ns,a}_{g,g}(k)$. More precisely,
$F_{g,g}(R)$ consists of isomorphism classes of flat proper families $\pi:C\to \Spec(R)$ equipped
with sections $p_1,\ldots,p_g$ and trivializations $v_i$
of the relative tangent bundle along each $p_i$, such that the base change of 
$(C,p_1,\ldots,p_g,v_1,\ldots,v_g)$ with respect to the reduction $R\to k$ gives
the cuspidal curve $(C^{\cusp}_g,p_1,\ldots,p_g)$ 
with the standard trivializations of the relative tangent bundle at $p_i$. 
Note that for such a family we automatically get that $\pi$ is smooth near the sections $p_i$ and
the divisor $p_1+\ldots+p_g$ is relatively ample 
(by \cite[(9.6.4)]{EGAIV}). Together with semicontinuity this implies that such a family defines an element of 
$\wt{U}^{ns,a}_{g,g}(R)$.

The functor $F_\infty$ is defined as follows.
We consider minimal $R$-linear $A_\infty$-structures on $E_g\ot R$, such that upon the
specialization $R\to k$ we get the trivial $A_\infty$-structure (i.e., the one with $m_i=0$ for $i>2$)
on $E_g$.  
We consider gauge
transformations between such structures that reduce to the identity under the
specialization $R\to k$. We define
$F_\infty(R)$ as the corresponding set of equivalence classes.

As in Section \ref{ainf-moduli-sec}, let $\MM_\infty$ denote the functor of minimal $A_\infty$-structures on $E_g$ up to gauge equivalence,
defined on all commutative $k$-algebras.

\begin{lem}\label{FS-homog-lem}
(i) We have a natural identification of $F_\infty(R)$ with 
the fiber of $\MM_\infty(R)\to \MM_\infty(k)$ over the
equivalence class of the  trivial $A_\infty$-structure on $E_g$.

\noindent
(ii) The functor $\MM_\infty$ is representable by an affine $k$-scheme.

\noindent
(iii) The functor $F_\infty$ is homogeneous.
\end{lem}

\Pf . (i) First, we have to check that if a minimal $R$-linear $A_\infty$-structure $m$ on $E_g\ot R$ 
reduces to an $A_\infty$-structure on $E_g$ that is gauge equivalent to the trivial one, then
there exists a gauge transformation $f$ over $R$ such that $f*m$ reduces to the trivial $A_\infty$-structure
under $R\to k$. This immediately follows from the fact that we can lift any gauge transformation defined over
$k$ to a gauge transformation defined over $R$.

It remains to show that if we have minimal $R$-linear $A_\infty$-structures
$m$ and $m'$ on $E_g\ot R$, reducing to the trivial one on $E_g$, and a gauge equivalence $f$ 
such that $f*m=m'$ then there exists a gauge equivalence $f'$ reducing to the identity on $E_g$
and such that we still have $f'*m=m'$. 
Let $\ov{f}$, $\ov{m}$, etc., denote the reduction with respect to $R\to k$.
Thus, $\ov{m}=\ov{m'}$ is the trivial $A_\infty$-structure on $E_g$.
Since $HH^1(E_g)_{<0}=0$ (see Corollary \ref{HH1-cor}), 
by Lemma \ref{important-homotopy-lem}, there exists a homotopy
$\ov{h}=(\ov{h}_n)$ over $k$ from the identity to $\ov{f}$. We can lift $\ov{h}$ to a homotopy $h$ over $R$
from the identity transformation of $m'$ to some gauge transformation $f_1$ with $f_1*m'=m'$
and $\ov{f_1}=\ov{f}$. Then setting $f'=f_1^{-1}\circ f$ gives the gauge transformation with the required
properties.

\noindent
(ii) By Corollary \ref{HH1-cor}, we have $HH^1(E_g)_{<0}=0$, so the representability of $\MM_\infty$ follows
from Corollary \ref{representability-cor}.

\noindent
(iii) Parts (i) and (ii) imply that $F_\infty$ is prorepresentable by the completion of
the algebra of functions on $\MM_\infty$ at the trivial point. Hence, $F_\infty$ is homogeneous.
\ed

The following result is well known (and in characteristic zero follows easily from the standard
description via the Maurer-Cartan equation) but we give a proof for completeness.
Note that for any $A_\infty$-algebra $A$ the Hochschild cochains
\begin{equation}\label{truncated-HH-eq}
CH^*(A)_{\le 0}=\prod_{t\le 0}CH^*(A)_t
\end{equation}
form a subcomplex of the Hochschild complex. We denote the corresponding cohomology
by $HH^*(A)_{\le 0}$.

\begin{lem}\label{ainf-obstruction-lem}
Let $E$ be a minimal $A_\infty$-algebra over $k$. 
Let us consider the functor on $\Art_k$ associating
with $R\in \Art_k$ the set of deformations of $E$ to a minimal $A_\infty$-algebra structure on $E\ot R$ 
up to extended gauge transformations (see Definition \ref{gauge-def}). Note that here we allow to deform $m_2$ as well.
Then the tangent space to this deformation functor is naturally identified with
$HH^2(E)_{\le 0}$. Furthermore, there is a complete obstruction theory for this functor
with values in $HH^3(E)_{\le 0}$.
The similar statements hold for deformations of a small minimal $A_\infty$-category.
\end{lem}

\Pf . It is convenient to think of deformations of an $A_\infty$-structure $m$ as that of the coderivation
$D_m$ on $\Bar(E)$ such that $D_m^2=0$ (see Section \ref{ainf-gen-sec}).
To compute the tangent space we have to look at $A_\infty$-structures on $E\ot k[t]/(t^2)$
reducing to the given one modulo $t$. The corresponding coderivation has
the form $D_m+tD_c$ with $c\in CH^2(E)_{\le 0}$, so the equation 
$(D_m+tD_c)^2=0$ is equivalent to $[D_m,D_c]=0$, i.e., to $c$ being a Hochschild cocycle.
An extended gauge equivalence corresponds to an automorphism of $\Bar(E)\ot k[t]/(t^2)$ of the form 
$\id+tD_f$, so it changes $c$ to $c-\de(f)$. This gives the required identification of the tangent space.

Given a small extension 
$$0\to I\to \wt{R}\to R\to 0$$
in $\Art_k$ (recall that this means that $I$ is annihilated by the maximal ideal $M$ of $\wt{R}$),
and an $A_\infty$-structure on $E\ot R$, we can lift each $m_i$ to some
Hochschild cochain $\wt{m}_i\in \wt{R}\ot CH^2(E)_{2-i}$ (where for $i=2$ we take the standard product
on $E\ot\wt{R}$).
Let $D_0$ be the coderivation of $\Bar(E\ot \wt{R})$ which corresponds to the $\wt{R}$-linear
extension of the original $A_\infty$-algebra structure on $E$.
Then the coderivation associated with $\wt{m}$ has form 
$D_0+D$, where $D$ takes values in $M\ot\Bar(E)$.
The $A_\infty$-equations hold modulo $I$, hence
$$(D_0+D)^2=[D_0,D]+D^2=D_{\phi}$$
for some $\phi\in I\ot CH^3(E)_{\le 0}$. 
We have 
$$[D_0,D_{\phi}]=[D_0,D^2]=[[D_0,D],D]=[D_\phi-D^2,D]=[D_\phi,D],$$
which vanishes since $MI=0$.
It follows that $\phi$ is a Hochschild cocycle. If we choose different liftings of $m_i$ then
$D$ would change to $D+D'$ where $D'$ takes values in $I\ot\Bar(E)$. Then
$$(D_0+D+D')^2=D_\phi+[D_0,D']+[D,D']+(D')^2.$$
Here $[D,D']=(D')^2=0$ since $MI=I^2=0$, so $\phi$ would change by a Hocshchild coboundary.
Thus, the class of $\phi$ in $I\ot HH^3(E)_{\le 0}$
is well defined. Conversely, if this class is zero then we can correct our choice of $\wt{m}_i$
to make $\phi=0$, so that $\wt{m}_i$ define an $A_\infty$-structure.
Thus, $\phi$ is a complete obstruction for the functor of $A_\infty$-structures.
The fact that extended gauge equivalences act trivially on this obstruction follows from the general theory
(see \cite[Lem.\ 2.20]{Manetti})
since this group is smooth (as a functor on $\Art_k$).
\ed 

%Need to check that it is a deformation functor - define it as the quotient of the functor of $A_\infty$-structure
%by the group functor of strict homotopies as in \cite{Manetti}??? 
%Note ???

\begin{lem}\label{comp-def-lem} 
Let $F\rTo{f} G\rTo{g} H$ be morphisms of deformation functors $\Art_k\to \Sets$.
Assume that $g\circ f:F\to H$ is smooth and the map of tangent spaces $t_F\to t_G$ is surjective.
Then $f:F\to G$ is smooth.
\end{lem}

\Pf . Let $O_F\to O_G\to O_H$
be the maps between the universal obstruction spaces induced by $f$ and $g$.
Since $F\to H$ is smooth, the composed map $O_F\to O_H$ is injective by
\cite[Prop.\ 2.18]{Manetti}. Therefore, the map $O_F\to O_G$ is injective.
Applying the cited result again we deduce that $F\to G$ is smooth.
\ed

The following result is a key step in the proof of Theorem A.

\begin{prop}\label{main-def-prop} 
Assume the characteristic of $k$ is not $2$ or $3$.
Then the morphism of deformation functors $a_\infty:F_{g,g}\to F_\infty$ is an isomorphism.
The tangent space to $F_\infty$ is naturally isomorphic to 
$$HH^2(E_g)_{<0}=HH^2(E_g).$$
\end{prop}

\Pf . Let $j$ denote the open embedding of $U=C^{\cusp}_g\setminus \{p_1,\ldots,p_g\}$ into $C^{\cusp}_g$,
and let $M=\Ext^*(G,j_*\OO_U)=\OO(U)$ be the corresponding $A_\infty$-module over $E_g$. 
Note that endomorphisms of $M$ in the derived category of $A_\infty$-modules are identified with the algebra $\OO(U)$.
Let us define the $A_\infty$-category $\CC$ as follows: it has $E_g$ as a full $A_\infty$-subcategory 
(where we think of $E_g$ as a category with objects $\OO_C$, $\OO_{p_1},\ldots,\OO_{p_g}$) and
one additional object $\OO_U$. The only additional morphisms in $\CC$ are elements of 
$M$, viewed as morphisms from $\OO_C$ to $\OO_U$, and the algebra $\OO(U)$ viewed as endomorphisms
of the object $\OO_U$. The $A_\infty$-structure on $\CC$ comes
from the products on $E_g$, the $A_\infty$-module structure on $M$, and the identification of
$\OO(U)$ with $A_\infty$-endomorphisms of $M$.
Let $F_{\CC}$ be the functor
of deformations of $\CC$ as a minimal $A_\infty$-category up to extended gauge transformations.
%, preserving the standard $m_2$ on the subcategory $\OO, \OO_{p_1},\ldots, \OO_{p_g}$. 
Let also 
$F_{U,nc}$ (resp., $F_U$) be the functor of deformations of $\OO(U)$ as an 
associative algebra (resp., as a commutative algebra). More precisely,
for $R\in \Art_k$ the set $F_{U,nc}$ (resp., $F_U$) consists of associative (resp., commutative)
$R$-algebra structures on $\OO(U)\ot R$ reducing to the given on under the homomorphism $R\to k$.
We have a commutative diagram of natural morphisms of functors
\begin{equation}\label{def-fun-diagram}
\begin{diagram}
F_U&\lTo{}& F_{g,g} &\rTo{} & F_\infty\\
\dTo{}&&\dTo{}&&\dTo{}\\
F_{U,nc}&\lTo{}& F_{\CC} &\rTo{}& \wt{F}_\infty 
\end{diagram}
\end{equation}
where $\wt{F}_\infty$ is the functor of deformations of $E_g$ as a minimal
$A_\infty$-category up to extended gauge transformations.

\noindent
{\bf Step 1}. The morphism $F_{\CC}\to \wt{F}_\infty$ is \'etale, and the tangent spaces to both
functors are naturally isomorphic to $HH^2(E_g)$.
To prove that this morphism is \'etale it is enough to check that it induces an isomorphism on tangent spaces
and an embedding on obstruction spaces (see \cite[Prop.\ 2.17]{Manetti}).
By Lemma \ref{ainf-obstruction-lem}, these spaces are given by $HH^2_{\le 0}$ and $HH^3_{\le 0}$
(applied to $\CC$ and $E_g$), respectively.  Note that our morphism corresponds to
the embedding of the full subcategory on the objects $\OO, \OO_{p_1},\ldots,\OO_{p_g}$ into 
$\CC$. Since $\OO(U)$ is the algebra of endomorphisms of the $A_\infty$-module $M$, it follows 
that this embedding induces an isomomorphism
$$HH^*(\CC)\rTo{\sim} HH^*(E_g)$$
(see \cite{Keller-di}, \cite[Thm.\ 4.1.1]{L-VdB}).
We have an exact sequence
\begin{equation}\label{Hochschild-ex-seq}
\begin{array}{l}
HH^1(\CC)\to 
H^1\bigl(CH(\CC)_{\ge 1}\bigr)\to HH^2(\CC)_{\le 0}\to HH^2(\CC)\to H^2\bigl(CH(\CC)_{\ge 1}\bigr)\to \\ 
HH^3(\CC)_{\le 0}\to HH^3(\CC),
\end{array}
\end{equation}
where $CH(?)_{\ge i}:=CH(?)/CH(?)_{\le i-1}$.
Since $E_g$ has trivial higher products, we have a canonical decomposition
$HH^i(E_g)=\prod_j HH^i(E_g)_j$, so the similar exact sequence for $E_g$
has trivial connecting homomorphisms. 
We claim that the map
$$HH^2(\CC)_{\le 0}\to HH^2(\CC)$$
is an isomorphism, while the map
$$HH^3(\CC)_{\le 0}\to HH^3(\CC)$$
is injective.
Indeed, first we observe that
$$H^2\bigl(CH(\CC)_{\ge 1}\bigr)=H^2\bigl(CH(E_g)_{\ge 1}\bigr)=0.$$
Indeed, since $\CC$ has morphisms only of degree $0$ and $1$, 
the only possible cochains in $CH^{s+t}(\CC)_t$ with $s+t=2$ and $t\ge 1$
correspond to $(s,t)=(1,1)$. But such cochains
should have form $\Hom^0(X,Y)\to\Hom^1(X,Y)$, so they all vanish
(since we exclude the identities in the Hochschild complex). The same argument
works for $E_g$.
Similar considerations with cochains show that 
$$H^1\bigl(CH(\CC)_{\ge 1}\bigr)=CH^1(\CC)_1=
\Ext^1(\OO,\OO)\oplus\bigoplus_{i=1}^g \Ext^1(\OO_{p_i},\OO_{p_i})\oplus\bigoplus_{i=1}^g \Ext^1(\OO_{p_i},\OO),$$
and the same formula holds for $E_g$.
Thus, in the commutative square
\begin{diagram}
HH^1(\CC)&\rTo{}&  H^1\bigl(CH(\CC)_{\ge 1}\bigr)\\
\dTo{}&&\dTo{}\\
HH^1(E_g)&\rTo{}&  H^1\bigl(CH(E_g)_{\ge 1}\bigr)
\end{diagram}
both vertical arrows are isomorphisms. Since the bottom horizontal arrow is surjective,
the top horizontal arrow is surjective too, and our claim follows from the exact sequence
\eqref{Hochschild-ex-seq}.
Our argument also shows that the map
$$HH^2(E_g)_{\le 0}\to HH^2(E_g)$$
is an isomorphism.
Now our assertion follows by considering the map from the exact sequence \eqref{Hochschild-ex-seq}
to the corresponding
exact sequence for $E_g$.

\noindent
{\bf Step 2}. The morphisms $F_{g,g}\to F_U$ and $F_U\to F_{U,nc}$ are both \'etale.
For the first morphism we can think geometrically, as
passing from a deformation of $(C,p_1,\ldots,p_g,v_1,\ldots,v_g)$ to that of $U=C\setminus D$. Thus, the maps of
tangent spaces and of obstruction spaces are the natural maps
$$\Ext^1_C\bigl(\bL_C,\OO_C(-2D)\bigr)\to \Ext^1_U(\bL_U,\OO_U) \ \text{ and} $$
$$\Ext^2_C\bigl(\bL_C,\OO_C(-2D)\bigr)\to \Ext^2_U(\bL_U,\OO_U).$$
By Lemma \ref{cotangent-complex-lem}(ii), these are isomorphisms.

Similar assertions for the morphism $F_U\to F_{U,nc}$ follow from Lemma 
\ref{curve-Hochschild-lem}.

\noindent
{\bf Step 3}. The morphisms $F_{g,g}\to F_{\CC}$ and
$F_{\CC}\to F_{U,nc}$ induce isomorphisms on tangent spaces. 
Indeed, by Step 2 and by the commutativity of diagram \eqref{def-fun-diagram},
the map on tangent spaces for the morphism $F_{\CC}\to F_{U,nc}$ is surjective. 
By Step 1, the tangent space to $F_{\CC}$ is identified with $HH^2(E_g)$.
On the other hand, the tangent space to $F_{U,nc}$ is
$$HH^2(U)\simeq HH^2(C^{\cusp}_g)\simeq HH^2(E_g)$$
(see Lemma \ref{HH-C-lem}(ii) and Proposition \ref{cusp-equiv-prop}), so it is a $k$-vector space
of the same dimension, which implies the assertion for $F_{\CC}\to F_{U,nc}$. 
The assertion for the other morphism
follows using Step 2 and the commutativity of diagram \eqref{def-fun-diagram}.
%Thus, $\G_m$ has negative weights on $HH^2(U)$???

\noindent
{\bf Step 4}. The morphisms $F_{g,g}\to F_\infty$ and $F_\infty\to\wt{F}_\infty$ 
induce isomorphisms on tangent spaces.
Indeed, Steps 1 and 3 plus the commutativity of \eqref{def-fun-diagram} imply that
the morphism of tangent spaces induced by $F_\infty\to\wt{F}_\infty$ is surjective.
Similarly to Lemma \ref{ainf-obstruction-lem} we can identify the tangent space to
$F_\infty$ with $HH^2(E_g)_{<0}$. Thus, the morphism of tangent spaces induced
by $F_\infty\to\wt{F}_\infty$ is the natural embedding
$$HH^2(E_g)_{<0}\hra HH^2(E_g),$$
hence, it is in fact an isomorphism. The assertion for the other morphism
follows using Steps 1 and 3 and the commutativity of diagram \eqref{def-fun-diagram}.

\noindent
{\bf Step 5}. By Step 2 and the commutativity of diagram \eqref{def-fun-diagram},
the composition $F_{g,g}\to F_{\CC}\to F_{U,nc}$ is \'etale.  
Applying Lemma \ref{comp-def-lem}, we deduce that the morphism 
$F_{g,g}\to F_{\CC}$ is smooth, hence \'etale (using Step 3).
Therefore, the composition 
$$F_{g,g}\to F_{\CC}\to \wt{F}_\infty$$
is also \'etale (using Step 1). Thus, we obtain that the composition 
$$F_{g,g}\to F_{\infty}\to \wt{F}_\infty$$
is \'etale. Applying Lemma \ref{comp-def-lem} again and using Step 4, we deduce that
the morphism $a_\infty: F_{g,g}\to F_\infty$ is \'etale.
But $F_\infty$ is homogeneous
(see Lemma \ref{FS-homog-lem}(iii)), so by \cite[Cor.\ 2.11]{Manetti}, $a_\infty$ is an isomorphism.
\ed

\subsection{Proof of Theorem A}\label{proof-sec}

The fact that $\wt{\UU}^{ns,a}_{g,g}\times \Spec(\Z[1/6])$ is an affine scheme was proved in Theorem \ref{moduli-thm}.

For the rest of the proof we change the base to a field $k$ of characteristic $\neq 2,3$, and 
denote the $k$-schemes obtained by the base change from 
$\wt{\UU}^{na,a}_{g,g}$ and $\SS_g$ by the same symbols. 
We have to prove that the morphism 
$$a_\infty:\wt{U}^{ns,a}_{g,g}\to \MM_\infty$$
(see Section \ref{map-sec}) is an isomorphism.
We will use the comparison of the deformation theories worked out above and in addition
will exploit the compatibility with $\G_m$-actions.

Recall that $\wt{U}^{ns,a}_{g,g}$ has a natural $\G_m$-action, coming from the diagonal $\G_m\sub \G_m^g$,
rescaling the trivializations of the tangent spaces at $p_1,\ldots,p_g$.
The morphism $a_\infty$ is compatible with this action and with the $\G_m$-action on the functor
$\MM_\infty$ of moduli of minimal $A_\infty$-structures on $E_g$, where
$\la\in\G_m$ sends $(m_n)\mapsto (\la^{-n+2}m_n)$. 

Now recall that by Corollaries \ref{HH1-cor} and \ref{representability-cor} the functor
$\MM_\infty$ of minimal $A_\infty$-deformations of $E_g$ is representable by an affine $k$-scheme.
Thus, we can view the morphism $a_\infty$
as a $\G_m$-equivariant morphism of affine schemes over $k$.
%(where we use an isomorphism $\wt{U}^{ns,a}_{g,g}\simeq \SS_g$ from Theorem \ref{moduli-thm}).
The cuspidal curve is sent by $a_\infty$ to the trivial $A_\infty$-structure and
Proposition \ref{main-def-prop} implies that $a_\infty$ induces an isomorphism between the completions of
our affine schemes at these points (which are both $\G_m$-invariant).
Recalling the explicit construction of the representable scheme in Theorem \ref{a-inf-moduli-thm} and
the description of the affine scheme $\SS_g$, we see that in both cases the maximal ideal of
the corresponding point is generated by elements of positive weight with respect to the $\G_m$-action
(note that the coefficients of $m_n$ viewed as functions on $\MM_\infty$ are sent by $\la$
to the coefficients of $(\la^{-1})^*m_n=\la^{n-2}m_n$, i.e., have degree $n-2$). 
Thus, in our situation we have a degree-preserving homomorphism of graded $k$-algebras
$$f:A=\bigoplus_{n\ge 0}A_n\to B=\bigoplus_{n\ge 0}B_n$$
with $A_0=B_0=k$ which induces an isomorphism between the completions with respect
to powers of $A_{>0}$ (resp., $B_{>0}$). Considering the induced
isomorphism
$$A/(A_{>0})^N\to B/(B_{>0})^N$$
for $N\gg 0$ we see that the map $A_n\to B_n$ is an isomorphism for every $n\ge 0$, so
in fact,  $f$ is an isomorphism.
\ed

%For deformations of $A$ this should mean that we are considering graded algebras $\wt{A}$ over

%For deformations of $C^{\cusp}$ this means that we are considering families $\CC\to \Spec(R[t]/t^N)$
%equipped with an action of $\G_m$, such that the special fiber is $C^{\cusp}\times\Spec(R)$.
%If we have such a family over $R[t]/t^N$ and looking at extensions to $R[t]/t^{N+1}$ then
%the spaces $\Ext^*(\bL,\OO_{C^{\cusp}})$ should get replaced by their $N$th components with respect to
%the grading associated with the $\G_m$-action on $C^{\cusp}$.

\subsection{Hochschild cohomology and normal forms of $A_n$-structures on $E_g$}\label{normal-sec}

We continue to work over a field $k$.
From our results we get a geometric way to calculate some Hochschild cohomology spaces
of the algebra $E_g$ (that were first determined in \cite{Fisette} via quite tedious algebraic calculations).

\begin{prop}\label{Hochschild-prop}
Assume that $\cha(k)$ is not $2$ or $3$.
Then the nonzero values of $\dim HH^2(E_g)_t$ are given by 
$$\dim HH^2(E_g)_{-1}=g^2-g, \ \ \dim HH^2(E_g)_{-2}=2g^2-2g,$$ 
$$\dim HH^2(E_g)_{-3}=g^2-g,\ \ \dim HH^2(E_g)_{-4}=g.$$
The nonzero values of $\dim HH^1(E_g)_t$ are
$$\dim HH^1(E_g)_0=g, \ \ \dim HH^1(E_g)_1=2g.$$
%Also $HH^3(E_g)$ is finite-dimensional???
\end{prop}

\Pf . By Proposition \ref{main-def-prop}, the space $HH^2(E_g)$ is the tangent space to 
$\MM_\infty\simeq \wt{\UU}^{ns,a}_{g,g}$ at the point corresponding to the cuspidal curve $C^{\cusp}_g$.
Furthermore, the component $HH^2(E_g)_t$ corresponds to the weight $t$ subspace with respect
to the $\G_m$-action (see Proposition \ref{cusp-equiv-prop}).
By Theorem \ref{moduli-thm}, we can identify $HH^2(E_g)$ with the tangent space to the affine scheme 
$\SS_g$ at the point where all coordinates are zero. 
By Proposition \ref{constants-prop}(ii),
the cotangent space at this point has the basis corresponding to the generators $\a_{ij}$, $\b_{ij}$, $\ga_{ij}$,
$\vareps_{ij}$ and $\pi_i$. Of these $\a_{ij}$ have weight $1$, $\b_{ij}$ and $\ga_{ij}$---weight $2$,
$\vareps_{ij}$---weight $3$, and $\pi_i$---weight $4$. Passing to the dual space we get the answers
for $HH^2(E_g)_t$.
To compute the components
$HH^1(E_g)_t$ we use the identification $HH^1(E_g)\simeq HH^1(C^{\cusp}_g)$
(see Proposition \ref{cusp-equiv-prop}), the exact sequence
$$0\to H^1(C,\OO)\to HH^1(C^{\cusp}_g)\to H^0(C,\TT)\to 0$$
from Lemma \ref{HH-C-lem}(i), and Lemma \ref{tangent-action-lem}.
\ed

Using the correspondence between non-special curves and $A_\infty$-structures on $E_g$
we can calculate normal forms of $A_n$-structures on $E_g$ for small $n$.
Let $\MM_n$ denote the moduli of minimal $A_n$-structures on $E_g$.
% (which is the same as $k[t]/(t^{n-1})$-linear $A_\infty$-structures on $E_g\ot k[t]/(t^{n-1})$). 
We know by Theorem \ref{a-inf-moduli-thm} and Corollary \ref{HH1-cor}
that these are affine $k$-schemes. 
Note that since $HH^3(E_g)\simeq HH^3(C^{\cusp}_g)$ is finite-dimensional,
by Corollary \ref{finite-a-inf-moduli-cor}, we have $\MM_\infty\simeq\MM_n$ for sufficiently
large $n$.

\begin{prop}\label{normal-prop} Assume that $\cha(k)$ is not $2$ or $3$.
Let 
$$\Theta=\{\a_{ij}, \b_{ij}, \ga_{ij}, \vareps_{ij}, \pi_i\}$$
be the set of minimal generators of the algebra of functions on the affine scheme $\SS_g$,
and let $I_r$ (resp., $I_{\le r}$) be the component of degree $r$ (resp., the sum of components of degree
$\le r$) in the ideal of relations between these generators
(see Section \ref{min-gen-sec} and Remark \ref{grading-rem}).

\noindent
(i) For $n\ge 3$ the algebra of functions on $\MM_n$ is isomorphic to
$$\Spec k[\Theta_{\le n-2}]/(I_{\le n-2}),$$
where $\Theta_{\le i}$ denotes the subset of elements in $\Theta$ of degree $\le i$.

\noindent
(ii) The schemes $\MM_3$ and $\MM_4$ are isomorphic to the affine spaces with coordinates
$\a_{ij}$ and $\a_{ij},\b_{ij},\ga_{ij}$, respectively
(equivalently, $I_{\le 2}=0$).
The univeral $A_3$-structure on $E_g$ is given by \eqref{m3-formulas-eq}.
The universal $A_4$-structure on $E_g$ is given by the formulas of the Appendix.
The algebra of functions on $\MM_5$ is isomorphic to $\Spec k[\a_{ij},\b_{ij},\ga_{ij},\vareps_{ij}]/(I_3).$
The universal $A_5$-structure on $E_g$ is given by the formulas of the Appendix.
\end{prop}

\Pf . (i) By Theorem A and Theorem \ref{moduli-thm}, we have an isomorphism of
affine schemes $\SS_g\simeq \MM_\infty$, compatible with $\G_m$-action.
Note that both these schemes have $\G_m$-equivariant closed embeddings into affine spaces
$$\SS_g\sub\A^{4g^2-3g}, \ \ \MM_\infty\sub \A^N,$$
where the former corresponds to the minimal generators 
\eqref{main-constants-eq}, while the latter corresponds to considering structure constants
of the normal forms $m_i$ (we use Theorem \ref{a-inf-moduli-thm} and
the fact that $\MM_\infty\simeq \MM_n$ for $n\gg 0$). 
Furthermore, by Corollary \ref{weight-cor}, there exists a closed $\G_m$-equivariant embedding
$\A^{4g^2-3g}\to \A^N$, which is a section for a $\G_m$-equivariant coordinate projection $\A^N\to \A^{4g^2-3g}$, such that we have a commutative diagram 
\begin{diagram}
\SS_g&\rTo{\sim}&\MM_\infty\\
\dTo{}&&\dTo{}\\
\A^{4g^2-3g}&\rTo{}&\A^N
\end{diagram}
By definition, the equations of $\MM_\infty$ in $\A^N$ are given by the $A_\infty$-constraints
$$[m_2,m_n]+[m_3,m_{n-1}]+\ldots=0,$$
$n\ge 3$, where the constraint corresponding to $n$ gives equations homogeneous of 
degree $n-2$. To get the presentation of $\MM_n$ from this we take coordinates of degree $\le n-2$
(corresponding to $m_i$ with $i\le 2$) and equations of degree $\le n-2$.
The assertion follows from this using Lemma \ref{graded-lem} below.

\noindent
(ii) The formula \eqref{m3-formulas-eq} gives $m_3$ as a linear combination of certain $g^2-g$ Hochschild
cohomology classes. To prove that this is a universal $A_3$-structure we need to check
that these Hochschild cohomology classes form a basis of $HH^2(E_g)_{-1}$. But this
follows from the description of the tangent space to $\SS_g$ at zero given in Proposition \ref{constants-prop}(ii)
(see the proof of Proposition \ref{Hochschild-prop}).

Next, we need to check that every product $m_3$ given by \eqref{m3-formulas-eq} extends
to an $A_4$-structure, i.e., that $[m_3,m_3]$ is a Hochschild coboundary. Indeed, let us set
$$m_4( X_j, B_i, Y_i, A_i)=\sum_{k\neq j}\a_{ij}\a_{jk} X_k,$$
$$m_4( X_j, B_i, A_i, X_i)=\sum_{k\neq j}\a_{ij}\a_{jk} X_k,$$
and let all other $m_4$ products of the basis elements be zero.
Then one can check by a direct computation that the $A_\infty$-constraint $[m_2,m_4]+[m_3,m_3]/2=0$ 
is satisfied. Note that the general formulas for $m_4$ in the Appendix differ from this one by
adding the linear combinations of certain $2(g^2-g)$ Hochschild cocycles with $\b_{ij}$ and $\ga_{ij}$
as coefficients. It remains to use Proposition \ref{constants-prop}(ii) again to see that the corresponding
Hochschild cohomology classes form a basis in $HH^2(E_g)_{-2}$.
The assertion for $A_5$-structures follows similarly using part (i).
\ed

The following Lemma was used in the proof of Proposition \ref{normal-prop}(i) (with $S$ being the algebra of functions on $\A^N$,
$\ov{S}$ the algebra of functions on $\A^{4g^2-3g}$ and $J$ the ideal of relations defining $\MM_\infty$
in $\A^N$).

\begin{lem}\label{graded-lem} Let $S$ be a $\Z$-graded ring, $\ov{S}\sub S$ a graded subring and
$K\sub S$ a homogeneous ideal such that $S=K\oplus\ov{S}$. Let also $J\sub S$ be a homogeneous ideal
containing $K$, so that $J=K\oplus\ov{J}$, where $\ov{J}$ is a homogeneous ideal in $\ov{S}$.
For an integer $n$ let $S'\sub S$ (resp., $\ov{S}'\sub\ov{S}$)
denote the subring generated by $S_{\le n}$ (resp., $\ov{S}_{\le n}$), and let $J'\sub S'$ (resp.,
$\ov{J}'\sub \ov{S}'$) be the ideal generated by $J_{\le n}$ (resp., $\ov{J}_{\le n}$).
Then the natural homomorphism 
$$S'/J'\to \ov{S}'/\ov{J}'$$
is an isomorphism.
\end{lem}

\Pf . Let $K'\sub S'$ denote the ideal generated by $K_{\le n}$. Since $\ov{S}$ is a subring of $S$,
from the decomposition
$$S_{\le n}=K_{\le n}\oplus \ov{S}_{\le n}$$
one can easily derive that $S'=K'+\ov{S}'$. Since $K'\sub K$ and $\ov{S}'\sub\ov{S}$, we get a direct sum decomposition
$$S'=K'\oplus\ov{S}'.$$
Now the equality $J_{\le n}=K_{\le n}\oplus \ov{J}_{\le n}$ implies the decomposition
$$J'=K'\oplus\ov{J}',$$
and the assertion follows.
\ed

\begin{rems} 1. One can check by a direct but tedious inspection of the equations defining $\SS_g$ 
(see the proof of Lemma \ref{grobner-lem})
that in fact the equations \eqref{ct-rel} span $I_3$, so they
form a defining set of relations for the scheme $\MM_5$. 

\noindent
2. By Corollary \ref{finite-a-inf-moduli-cor}, the vanishing of $HH^2(E_g)_{<-4}$ implies that 
the morphism $\MM_\infty\to \MM_6$ is a closed embedding.
%every minimal $A_\infty$-structure on $E_g$ is determined by the corresponding $A_6$-structure, i.e.,
On the other hand, using the isomorphism of Theorem A we can interpret Proposition \ref{constants-prop2} as
saying that the morphisms $\MM_\infty\to \MM_5$ and $\MM_\infty\to \MM_4$ become closed embeddings
when restricted to the preimages of explicit open subsets given in terms of nonvanishing of
some of $\a_{ij}$ (all of them, for the second morphism).
\end{rems}

\section{$\psi$-(pre)stable curves of genus $0$ and $A_\infty$-structures} 
\label{genus-0-sec}

In this section we will analyze what happens with the picture described above if we replace curves of genus $g$
with $g$ marked points by curves of genus $0$ with $n$ marked points (where $n\ge 3$). 

\subsection{The moduli stack of $\psi$-prestable curves of genus $0$}

First, similarly to Sec.\ \ref{curves-sec}
we consider certain normal forms for curves of genus $0$ with $n$ marked points, and the resulting
presentation of the moduli space by explicit equations. 

\begin{defi}\label{psi-stab-def} 
(i) Let $(C,p_1,\ldots,p_n)$ be a reduced connected complete curve 
of arithmetic genus $0$ with $n$ smooth marked points over an algebraically closed field.
We say that $(C,p_1,\ldots,p_n)$ is {\it $\psi$-prestable} if $\OO_C(p_1+\ldots+p_n)$ is ample,
or equivalently, if every component of $C$ contains at least one marked point.

\noindent
(ii) A $\psi$-prestable curve is called {\it $\psi$-stable} if every component of $C$ contains at least three
distinguished points (i.e., singular or marked points).
\end{defi}

It is known that the only singular points that can occur for reduced curves of arithmetic genus zero are
{\it rational $m$-fold points}, i.e., points for which the completion of the local ring has form
$$\hat{\OO}_{C,p}\simeq k[[x_1,\ldots,x_m]]/(x_ix_j \ |\  1\le i<j\le m)$$
(see \cite[Lem.\ 1.17]{Smyth} or \cite{Stevens}). It follows that the above definition of a $\psi$-stable curve
is equivalent to the one given in \cite[Def.\ 2.26]{FS} and in \cite{Boggi}.

For $n\ge 3$ let us denote by $\UU_{0,n}[\psi]$ the stack of $\psi$-prestable curves classifying families 
(flat, proper, finitely presented morphisms) with $n$ sections, whose geometric fibers are $\psi$-prestable
curves. Let also $\wt{\UU}_{0,n}[\psi]\to \UU_{0,n}[\psi]$ denote the $\G_m^n$-torsor corresponding to
choices of nonzero tangent vectors at the marked points.

For a $\psi$-prestable curve $(C,p_1,\ldots,p_g)$ with fixed nonzero tangent vectors $v_i$ at $p_i$ we
can construct a natural presentation for the algebra $H^0(C\setminus D,\OO)$, where
$D=p_1+\ldots+p_n$. Namely, for every $i=1,\ldots,g$ let us pick a rational function $f_i\in H^0\bigl(C,\OO(p_i)\bigr)$
such that $f_i\equiv v_i\mod \OO_C$ in $\OO_C(p_i)/\OO_C\simeq T_{p_i}C$. 
Note that such $f_i$ is defined uniquely up to an additive constant. To get rid of this ambiguity we
will impose the condition $f_i(p_{i+1})=0$ for $i=1,\ldots,n$ (where $p_{n+1}=p_1$).

A version of this construction works for a family $\pi:C\to\Spec(R)$ of $\psi$-prestable curves with 
$n$ marked points $p_i:\Spec(R)\to C$, where we choose $v_i\in H^0\bigl(C,\OO_C(p_i)/\OO_C\bigr)$ inducing
trivializations $R\simeq H^0\bigl(C,\OO_C(p_i)/\OO_C\bigr)$.
Then the natural map $R\to H^0(C,\OO_C)$ is an isomorphism and $H^1(C,\OO_C)=0$, so
we get an exact sequence
$$0\to R\to H^0\bigl(C,\OO_C(p_i)\bigr)\to H^0\bigl(C,\OO_C(p_i)/\OO_C\bigr)\to 0$$
so we can choose $f_i\in H^0\bigl(C,\OO_C(p_i)\bigr)$ projecting to $1\in R\simeq H^0\bigl(C,\OO_C(p_i)/\OO_C\bigr)\to 0$.

Let us set $D=p_1\bigl(\Spec(R)\bigr)+\ldots+p_n\bigl(\Spec(R)\bigr)$. Also, for $i\neq j$ we set
$$\a_{ij}=f_i(p_j).$$

\begin{lem}\label{genus-0-equation-lem} 
Assume $n\ge 3$. Let $\pi:C\to\Spec(R)$ be a family of $\psi$-prestable curves with $n$ marked
points $p_i:\Spec(R)\to C$ and with trivializations $R\simeq H^0\bigl(C,\OO_C(p_i)/\OO_C\bigr)$.
Then the algebra $H^0(C\setminus D,\OO)$ is generated over $R$ by the elements $f_1,\ldots,f_g$
with the defining relations
\begin{equation}\label{genus-0-curve-eq}
f_if_j=\a_{ij}f_j+\a_{ji}f_i+c_{ij}, \ \text{ for } i\neq j,
\end{equation}
where 
\begin{equation}\label{c-ij-a-eq}
c_{ij}=\a_{ik}\a_{jk}-\a_{ij}\a_{jk}-\a_{ji}\a_{ik}
\end{equation}
for any three distinct indices $i,j,k$.
Furthermore, for every $N\ge 1$ the elements $1$ and $f_i^m$, with $i=1,\ldots,n$, $1\le m\le N$, 
form a basis of $H^0\bigl(C,\OO_C(ND)\bigr)$ over $R$.
\end{lem}

\Pf . Since $H^1(C,\OO)=0$, for any $i\neq j$ we have an exact sequence
$$0\to R=H^0(C,\OO)\to H^0\bigl(C,\OO(p_i+p_j)\bigr)\to H^0\bigl(C,\OO(p_i+p_j)/\OO\bigr)\to 0,$$
so $(1, f_i, f_j)$ is an $R$-basis of $H^0\bigl(C,\OO(p_i+p_j)\bigr)$. Since $f_if_j\in H^0\bigl(C,\OO(p_i+p_j)\bigr)$,
we get the equation \eqref{genus-0-curve-eq} (the coefficients with $f_i$ and $f_j$ are determined
by looking at polar parts at $p_i$ and $p_j$).
Now the equation \eqref{c-ij-a-eq} is obtained by evaluating both parts of \eqref{genus-0-curve-eq} at $p_k$.
The last equation follows by induction on $N$ from the exact sequences
$$0\to H^0\bigl(C,\OO(jD)\bigr)\to H^0\bigl(C,\OO((j+1)D)\bigr)\to H^0\bigl(C,\OO((j+1)D)/\OO(jD)\bigr)\to 0,$$
where $j\ge 0$.
\ed

We also have a much simpler version of Lemma \ref{grobner-lem} in our situation.

\begin{lem}\label{genus-0-grobner-lem}
Let $R$ be a commutative ring.
An associative $R$-algebra generated by $f_1,\ldots,f_g$ with
defining relations \eqref{genus-0-curve-eq} has elements $1$, $(f_i^m)_{i=1,\ldots,n, m\ge 1}$ as 
an $R$-basis if and only if the relations \eqref{c-ij-a-eq} hold for any distinct $i,j,k$.
\end{lem}

\Pf . This is an immediate application of the Gr\"obner basis technique: expanding $f_if_jf_k$, for
distinct $i,j,k$, in two ways using the defining relations \eqref{genus-0-curve-eq} we get
the equation \eqref{c-ij-a-eq} and in addition, the equation
$$\a_{ij}c_{jk}+\a_{ji}c_{ik}=\a_{ik}c_{kj}+\a_{ki}c_{ij}.$$
It remains to observe that the latter equation follows from \eqref{c-ij-a-eq} (applied for appropriate permutations
of indices $i,j,k$).
\ed

Let us consider the subscheme $\TT_n$ in the affine space (over $\Z$)
with coordinates $\a_{ij}$, where $1\le i,j\le n$,
$i\neq j$, defined by the equations 
\begin{equation}\label{Boggi-eq}
\begin{array}{l}
\a_{i,i+1}=0, \\
\a_{ik}\a_{jk}-\a_{ij}\a_{jk}-\a_{ji}\a_{ik}=\a_{il}\a_{jl}-\a_{ij}\a_{jl}-\a_{ji}\a_{il},
\end{array}
\end{equation}
where $i,j,k,l$ are distinct (with the convention $\a_{n,n+1}=\a_{n,1}$ in the first equation).
Note that for $n=3$ we have $\TT_n\simeq\A^3$ (we can take $\a_{13}$, $\a_{21}$ and $\a_{32}$
as coordinates).

\begin{thm}\label{genus-0-moduli-thm} Assume $n\ge 3$.
There is an isomorphism of the stack $\wt{\UU}_{0,n}[\psi]$ with the affine scheme $\TT_n$, so that
the open part $C\setminus D$ of the universal curve over $\TT_n$ is given by
the equations \eqref{genus-0-curve-eq} in affine coordinates $f_i$, where $c_{ij}$ is given by
\eqref{c-ij-a-eq}. Furthermore, the isomorphism $\wt{\UU}_{0,n}[\psi]\simeq\TT_n$
is compatible with the $\G_m^n$-actions, where $\la=(\la_i)\in\G_m^n$ acts on $\wt{\UU}_{0,n}[\psi]$
by rescaling the tangent vectors $v_i\mapsto \la_i^{-1}v_i$, and on the scheme $\TT_n$ so that
$$(\la^{-1})^*\a_{ij}=\la_i\a_{ij}.$$
\end{thm}

\Pf . We mimic the proof of Theorem \ref{moduli-thm}.
Lemma \ref{genus-0-equation-lem} gives a natural map
\begin{equation}\label{genus-0-moduli-map-1}
\wt{\UU}_{0,n}[\psi](R)\to \TT_n(R)
\end{equation}
for every commutative ring $R$.
To construct the map in the other direction we note that by Lemma \ref{genus-0-grobner-lem},
an $R$-point of $\TT_n$ gives rise
to an associative $R$-algebra $A$ with generators $f_1,\ldots,f_n$ and defining relations
of the form \eqref{genus-0-curve-eq} with $c_{ij}$ given by \eqref{c-ij-a-eq} (where $\a_{ij}\in R$
are the coordinates of our point in $\TT_n(R)$), such that $1$ and $(f_i^m)$ form an $R$-basis of $A$.
For $N\ge 1$ let $F_NA\sub A$ denote the $R$-submodule generated by $f_i^m$ with $m\le N$,
and let $F_0A=R\sub A$. Then we define the projective family $\pi:C\to\Spec(R)$ by setting
$C=\Proj(\RR A)$, where $\RR A=\bigoplus_{j\ge 0} F_jA$ is the Rees algebra associated
with the filtration $(F_\bullet A)$. Let $T$ denotes the element $1\in F_1A\sub \RR A$. Then
$\RR A/(T)$ is isomorphic to the associated graded quotient 
$$\gr^\bullet_{F} A\simeq R[f_1,\ldots,f_n]/(f_if_j \ |\ 1\le i<j\le n).$$
Hence, $\Proj\bigl(\RR A/(T)\bigr)$ is the disjoint union of $n$ copies of $\Spec(R)$, which gives $n$ disjoint
sections $p_1,\ldots,p_n:\Spec(R)\to C$, so that $D=p_1\bigl(\Spec(R)\bigr)+\ldots+p_n\bigl(\Spec(R)\bigr)$ is the
divisor $(T=0)$. Let $F_i$ denote the element $f_i\in F_1A\sub \RR A$.
Then the algebra $\RR A$ is defined by the homogeneous versions of the relations
\eqref{genus-0-curve-eq}:
\begin{equation}\label{genus-0-homog-eq}
F_iF_j=\a_{ij}F_jT+\a_{ji}F_iT+c_{ij}T^2.
\end{equation}
Since $\deg(F_i)=\deg(T)=1$, we can think of $C=\Proj(\RR A)$ as a closed subscheme in the projective
space $\P^n_R$. In particular, $\OO_C(D)=\OO_C(1)$, so $D$ is ample.
Let $U_i\sub C$ denote the open subset $F_i\neq 0$.
Then $D\sub U_i$ is the vanishing locus of the function $T/F_i$ and the relations \eqref{genus-0-homog-eq}
show that $F_j/F_i$ 
is divisible by $T/F_i$ on $U_i$. This implies that $p_i\bigl(\Spec(R)\bigr)$ is locally given by one equation 
$T/F_i=0$. The same argument as in the proof of Theorem \ref{moduli-thm} shows that
the morphism $\pi:C\to \Spec(R)$ is flat of relative dimension $1$. Also, since $f_i=F_i/T$, we obtain
that $f_i$ projects to a generator of $\OO_C(p_i)/\OO_C$ and is regular near the images of $p_j$ for 
$j\neq i$. Thus, we deduce that for $m\ge 0$, $H^0\bigl(C,\OO(mD)\bigr)$ is a free $R$-module of rank $mn+1$ (with
the basis given by the appropriate powers of $f_i$).
This implies that the geometric fibers of $\pi$ are connected reduced curves of arithmetic genus $0$.
Together with ampleness of $D$ this shows that they are $\psi$-prestable, so we get a map
\begin{equation}\label{genus-0-moduli-map-2}
\TT_n(R)\to \wt{\UU}_{0,n}[\psi](R).
\end{equation}
As in Theorem \ref{moduli-thm} one immediately checks that the maps \eqref{genus-0-moduli-map-1}
and \eqref{genus-0-moduli-map-2} are inverse of each other. The compatibility with the $\G_m^g$-actions
is straightforward. 
\ed

\begin{ex} For $n=3$ we have four types of $\psi$-prestable curves, characterized by the number 
of vanishing coordinates $\a_{21},\a_{32},\a_{13}$.
For a smooth curve none of them vanish. If $\a_{32}=0$ but $\a_{21}\a_{13}\neq 0$ 
then $C$ is the union of two
(smooth) components $C_{12}$ and $C_3$, intersecting transversally at one point, with $p_1,p_2\in C_{12}$ and $p_3\in C_2$. If $\a_{21}=\a_{32}=0$ then $C$ is the union of three smooth components $C_1,C_2,C_3$
such that $p_i\in C_i$, where $C_2$ and $C_3$ are disjoint and both $C_2$ and $C_3$ intersect
$C_1$ transversally (at distinct points). Finally, when all the coordinates are zero we get the union of
three components intersecting at one point which is a rational threefold point.
\end{ex}

%Coordinates $\a_{ij}$.

\subsection{$A_\infty$-structures for curves of arithmetic genus $0$}

Given a curve of arithmetic genus $0$ with $g$ (smooth) marked points we can consider the corresponding
$A_\infty$-structure on the algebra $\Ext^*(G,G)$ for $G=\OO_C\oplus\OO_{p_1}\oplus\ldots\oplus\OO_{p_n}$.
As before, a choice of nonzero tangent vectors at all $p_i$ gives a way to fix generators 
$B_i\in\Ext^1(\OO_{p_i},\OO_C)$ so that the algebra $\Ext^*(G,G)$, as a graded associative algebra
gets identified with the algebra 
$$E_{0,n}=k[Q]/J_0,$$
where $Q=Q_n$ is the same quiver that we used before (see \eqref{Eg-eq}), while the ideal $J_0$ is 
generated by the elements
$$B_iA_i, A_iB_j, \ \text{ where } i\neq j.$$

More precisely, as in Section \ref{Cech-ainf-sec}, the construction of the $A_\infty$-structure on 
$E_{0,n}$ associated with a curve works well in families with affine bases (depending on
a choice of relative parameters at the marked points). We only have to make appropriate changes in
the choices of cohomology representatives and of the homotopy operator.
Namely, since we now have $H^1(C,\OO_C)=0$, we don't need the classes $X_i$, and the homotopy
operator on $K_{\OO,\OO}$ has to be defined by the formulas $Q([v])=-v$ and
$$Q([\frac{1}{t_i^n}])=f_i[n](t_i)_{\ge 0}\unit_i+\sum_{j\neq i}f_i[n](t_j)_{\ge 0}\unit_j+f_i[n]$$
for all $n\ge 1$, where $f_i[n]$ are now defined for all $n\ge 1$ (in fact, we can take $f_i[n]=f_i^n$,
where $f_i$ are defined as in Lemma \ref{genus-0-equation-lem}).

Let us fix a field $k$. We denote still by $\wt{\UU}_{0,n}[\psi]$ the moduli scheme of $\psi$-prestable curves
over $k$ (isomorphic to the affine scheme $\TT_n\times\Spec(k)$).
From the above construction of $A_\infty$-structures, as in section \ref{map-sec}, we get a morphism of functors on commutative algebras over $k$
\begin{equation}\label{genus-0-a-inf-map}
\wt{\UU}_{0,n}[\psi]\to \MM_\infty(E_{0,n}),
\end{equation}
where $\MM_\infty(E_{0,n})$ is the functor of minimal $A_\infty$-structures on 
the algebra $E_{0,n}$, up to a gauge equivalence. Furthermore, the map \eqref{a-inf-map} is compatible
with the natural $\G_m^n$-action, where the action on $\MM_\infty(E_{0,n})$ is induced by the rescalings
$B_i\mapsto \la_iB_i$.

Let us denote by $C_{0,n}$ the most singular curve appearing in the moduli space $\wt{\UU}_{0,n}$:
it  is the union of $n$ projective lines intersecting at one point which is a rational $n$-fold point (with one
marked point on each component). 
Under the isomorphism of Theorem \ref{genus-0-moduli-thm} this point in the moduli space corresponds
to the origin in the affine scheme $\TT_n$ (where all $\a_{ij}$ vanish).
Note that as in the case of the cuspidal curve, we have a natural $\G_m^n$-action on $C_{0,n}$,
and we will use the induced action of the diagonal $\G_m$.
The analog of Proposition \ref{cusp-equiv-prop} states in our case that we have an equivalence
of categories
$$\Per(C_{0,n})\simeq \Per(E_{0,n})$$
inducing a $\G_m^n$-equivariant isomorphism 
\begin{equation}\label{genus-0-HH-isom-eq}
HH^*(C_{0,n})\simeq HH^*(E_{0,n}),
\end{equation}
so that the second grading on these spaces is given by the weights of the $\G_m$-action.

Next, we observe that the assertions of Lemma \ref{tangent-action-lem}(ii)-(v) hold for the $n$-fold rational singularity curve
$C_{0,n}$, with the same proofs but without any assumption on the characteristic.
It follows that the analog of Lemma \ref{HH-C-lem}(ii) also holds for $C_{0,n}$, so 
$$HH^1(C_{0,n})_{<0}=0$$
and the natural map
$$HH^2(C_{0,n})\to HH^2(U)$$
is an isomorphism, where $U\sub C_{0,n}$ denotes the complement to the $g$ marked points.
Using the isomorphism \eqref{genus-0-HH-isom-eq} we deduce the vanishing
$$HH^1(E_{0,n})_{<0}=0,$$ 
which implies in particular, that the functor $\MM_\infty(E_{0,n})$ is represented by an affine scheme over 
$k$ (see Corollary \ref{finite-a-inf-moduli-cor}). 

Next, similarly to Lemma \ref{cotangent-complex-lem}(ii), we derive that for $C=C_{0,n}$ the natural maps
$$\Ext^1_C\bigl(\bL_C,\OO(-2D)\bigr)\to \Ext^1_U(\bL_U,\OO_U) \ \text{ and} $$
$$\Ext^2_C\bigl(\bL_C,\OO(-2D)\bigr)\to \Ext^2_U(\bL_U,\OO_U)$$
are isomorphisms.
This allows us to run the same argument as in Proposition \ref{main-def-prop} to show that
the deformation theory of $C_{0,n}$ matches with the deformation theory of the trivial $A_\infty$-structure
on $E_{0,n}$. Finally, using the $\G_m$-actions as in \ref{proof-sec} we deduce the following result
(this time we don't need any assumptions on the characteristic).

\begin{thm}\label{genus-0-a-inf-moduli-thm} 
For any field $k$ and any $n\ge 3$,
the morphism \eqref{genus-0-a-inf-map} induces an isomorphism of the moduli scheme
of $\psi$-prestable curves $\wt{\UU}_{0,n}[\psi]$ with the moduli scheme of minimal $A_\infty$-structures
on $E_{0,n}$, up to a gauge equivalence. This isomorphism is compatible with the natural $\G_m^n$-actions.
\end{thm}

Similarly to Proposition \ref{normal-prop} we can interpret the fact that the algebra of functions on
the scheme $\wt{\UU}_{0,n}[\psi]\simeq\TT_n$ is generated in degree $1$ with quadratic relations
as saying that the natural projection to the moduli space of $A_4$-structures,
$$\MM_\infty(E_{0,n})\to \MM_4(E_{0,n}),$$ 
is an isomorphism, while the map $\MM_4(E_{0,n})\to\MM_3(E_{0,n})$ is a closed embedding (where
$\MM_3(E_{0,n})$ is an affine space).
In particular, any minimal 
$A_\infty$-structure on $E_{0,n}$ is determined up to a gauge equivalence by $m_3$.
Using the formulas from Section \ref{Cech-ainf-sec} one can check that the nonzero values of $m_3$ on
the basis vectors are given by
$$m_3(A_j,B_i,A_i)=\a_{ij}A_j, \ \ \ m_3(B_j,A_j,B_i)=-\a_{ji}B_i,$$
where $i\neq j$.

\subsection{$\psi$-stability as GIT stability}

In this section we continue to work over a field $k$.
We are going to show that the notion of $\psi$-stability (see Definition \ref{psi-stab-def}(ii)) appears
naturally as a GIT stability for the action of $\G_m^n$ on the affine scheme $\wt{\UU}_{0,n}[\psi]\simeq\TT_n$. 

Namely, for a nontrivial character $\chi:\G_m\to \G_m$ we consider the GIT quotient
$\TT_n\sslash_\chi \G_m^n$ (see section \ref{GIT-sec} for our conventions on this).

Let us identify the character lattice of $\G_m$ with $\Z^n\sub\R^n$ in the standard way, and let 
$\bC_0\sub\R^n$ be the cone generated by all the basis vectors $e_i$.
Recall that $\TT_n$ is a closed subscheme of the affine space $\prod_{i=1}^n \A^{n-2}$,
where the coordinates $(\a_{ij})$ on the $i$th factor (where $j\neq i,i+1$)
satisfy $(\la^{-1})^*\a_{ij}=\la_i\a_{ij}$. 
Thus, we have a closed embedding of GIT quotients
$$\TT_n\sslash_\chi \G_m^n\sub \prod_{i=1}^n (\A^{n-2}\sslash_{a_i}\G_m),$$
where $\chi=\sum_i a_ie_i$. 
Thus, for $\chi\not\in\bC_0$ the GIT quotients will be empty, while for $\chi$ in the interior of $\bC_0$
the ambient GIT quotients will be just the product of $n$ copies of $\P^{n-3}$ (for $\chi$ on the boundary
it will be the product of the copies of $\P^{n-3}$ corresponding to positive coordinates of $\chi$).

\begin{prop}\label{psi-moduli-prop} 
Assume $n\ge 3$. For any $\chi$ in the interior of $\bC_0$ the $\chi$-semistable locus in
$\wt{\UU}_{0,n}[\psi]$ coincides with the $\chi$-stable locus and consists of $\psi$-stable curves.
Thus, the corresponding GIT quotient is exactly the moduli scheme of $\psi$-stable curves
$\ov{M}_{0,n}[\psi]$. The obtained embedding of $\ov{M}_{0,n}[\psi]$ into the product of $n$ copies of $\P^{n-3}$
corresponds to the line bundles $\psi_1,\ldots,\psi_n$ (the cotangent lines at the marked points).
\end{prop}

\Pf . Note that the conditions of stability and semistability for $\chi$ in the interior of $\bC_0$ are eqivalent, 
since this is true for the action on the ambient affine space. The condition of $\chi$-semistability is simply that 
for every $i$ there exists $j\neq i$ such that $\a_{ij}\neq 0$. We claim that this condition is equivalent
to the $\psi$-stability of the corresponding curve. Indeed, recall that $\a_{ij}=f_i(p_j)$, where 
$f_i\in H^0\bigl(C,\OO_C(p_i)\bigr)$. Let $C_i\sub C$ be the component containing $p_i$. Then $f_i$ has a unique
zero on $C_i$. If $(C,p_1,\ldots,p_n)$ is $\psi$-stable then $C_i$ has at least two distinguished points
$q$ and $q'$ other than $p_i$. Let $p_j$ (resp., $p_{j'}$) be either $q$ (resp., $q'$) if it is a marked point, 
or a point on the component attached to $C_i$ at $q$ (resp., at $q'$). Note that such marked points $p_j$,
$p_{j'}$ exist since $(C,p_1,\ldots,p_n)$ is $\psi$-prestable.
Then we have either $f_i(q)\neq 0$
or $f_i(q')\neq 0$ which implies that either $\a_{ij}=f_i(p_j)\neq 0$ or $\a_{ij'}=f_i(p_{j'})\neq 0$
(note that $f_i$ is constant on all components different from $C_i$), hence, we get $\chi$-semistability of
$(C,p_1,\ldots,p_n)$. Conversely, assume there exists a component $C_i\sub C$ with only two distinguished points and let us show that $(C,p_1,\ldots,p_n)$ is not $\chi$-semistable. Let $p_i\in C_i$ be a distinguished point and let $q\in C_i$ be another distinguished point (it exists by connectedness of $C$). Then
our normalization $f_i(p_{i+1})=\a_{i,i+1}=0$ implies that $f_i(q)=0$. Hence, $f_i$ is zero on all other components of $C$, and so $\a_{ij}=0$ for all $j\neq i$.

Let us take $\chi=e_1+\ldots+e_n$.
Then the identification of the pull-back of $\OO(1)$ from the $i$th copy of $\P^{n-3}$ with $\psi_i$ is 
immediate since both line bundles become trivial on $\wt{\UU}_{0,n}[\psi]$ and correspond to the same
action of $\G_m^n$ on the trivial line bundle (through the $i$th factor $\G_m$.
\ed

\begin{cor}\label{Boggi-cor} 
The moduli scheme $\ov{M}_{0,n}[\psi]$ is scheme-theoretically cut out in $(\P^{n-3})^n$ by the
equations \eqref{Boggi-eq}, where $(\a_{ij})_{j\neq i,i+1}$ are the homogeneous coordinates on the $i$th
copy of $\P^{n-3}$.
\end{cor}

\begin{rems} 1. Recall (see \cite[Rem.\ 2.30]{FS}) that there is a natural morphism 
$$f_\psi:\ov{M}_{0,n}\to \ov{M}_{0,n}[\psi],$$
where $\psi(C)$ is obtained from $C\in\ov{M}_{0,n}$ by contracting all components without marked points.
In fact, Proposition \ref{psi-moduli-prop} gives another proof of the fact that the line bundle 
$\psi=\psi_1+\ldots+\psi_n$ on $\ov{M}_{0,n}$ is semiample and that $f_\psi$ is the corresponding contraction
(see \cite[Sec.\ 4.2.1]{FS}). 
The composed morphisms $\ov{M}_{0,n}\to \ov{M}_{0,n}[\psi]\to \P^{n-3}$ 
given by the linear series $\psi_i$ were first considered by Kapranov \cite{Kapranov}. 
Note also that we have a commutative diagram
\begin{diagram}
\ov{M}_{0,n}&\rTo{f_\psi}&\ov{M}_{0,n}[\psi]\\
\dTo{}&&\dTo{}\\
\UU_{0,n}&\rTo{}&\UU_{0,n}[\psi]
\end{diagram}
where $\UU_{0,n}$ is the stack of (reduced connected complete) curves of arithmetic genus $0$ with
$n$ smooth marked points. Here the bottom horizontal arrow can be intepreted in view of Theorem
\ref{genus-0-a-inf-moduli-thm} as associating with a curve
$(C,p_1,\ldots,p_n)$ the corresponding $A_\infty$-structure on $E_{0,n}=\Ext^*(G,G)$
(where we identify $\UU_{0,n}[\psi]$ with the moduli stack of $A_\infty$-structures on $E_{0,n}$).

\noindent 2. One can similarly check that the 
GIT quotients $\TT_n\sslash_\chi \G_m^n$ for $\chi$ on the boundary of $\bC_0$ are given by
the images of the maps $\ov{M}_{0,n}\to (\P^{n-3})^r$ given by the line bundles $\psi_{i_1}+\ldots+\psi_{i_r}$
for $i_1<\ldots<i_r$.

\noindent 3.
If $(C,p_1,\ldots,p_n)$ is a $\psi$-prestable curve then one has also another natural generator
of the perfect derived category $\Per(C)$, namely
\begin{equation}\label{V-dec-eq}
V=\OO_C\oplus\bigoplus_{i=1}^n\OO_C(p_i).
\end{equation}
It is easy to check that $H^1\bigl(C,\OO_C(p_i-p_j)\bigr)=0$, so in fact, $V$ is a tilting bundle.
The coefficients $\a_{ij}$ appear in the structure constants of the algebra $\End(V)$ and
the equations among them correspond to the associativity equations. One can show that
in this way one gets an identification of  $\UU_{0,n}[\psi]$ with an open substack in the moduli stack of
$k^{n+1}$-algebra structures on $\End(V)$ (where the embedding $k^{n+1}\sub\End(V)$  comes
from the decomposition \eqref{V-dec-eq}). This open substack is characterized by the condition that
for all $i\neq j$ the compositions
\begin{align*}
&\Hom\bigl(\OO_C(p_i),\OO_C(p_j)\bigr)\rTo{x\mapsto x\ot 1}
\Hom\bigl(\OO_C(p_i),\OO_C(p_j)\bigr)\ot\Hom\bigl(\OO_C,\OO_C(p_i)\bigr)\rTo{\mu} \\
&\Hom\bigl(\OO_C,\OO_C(p_j)\bigr)\to 
H^0\bigl(\OO_C(p_j)/\OO_C\bigr),
\end{align*}
where $\mu$ comes from a $k^{n+1}$-algebra structure, are surjective.
\end{rems}

\bigskip

\noindent
{\sc APPENDIX: Formulas for $m_4$ and $m_5$.}

\bigskip

Here we assume that the characteristic is not $2$ or $3$.
We use the setup of Section \ref{homotopy-sec}, and in addition, use the constants introduced
in Section \ref{relations-sec}. By Proposition \ref{constants-prop}, 
all of these constants are some universal polynomials of the
generators \eqref{main-constants-eq}.
Recall that $P_i=A_i X_i= Y_i A_i$.

The computation of $m_4$ is based on the following formulas:
$$ B_iQ(P_j)=0, \  B_iQ(P_i)=[\frac{1}{t_i^2}], \  Y_iQ(P_j)=0, \  Y_iQ(P_i)=\be[\frac{1}{t_i^2}],$$
$$ X_iQ([\frac{1}{t_i^2}])=[\frac{f_i}{t_i}], \  X_jQ([\frac{1}{t_i^2}])=[\frac{f_i}{t_j}], \
Q([\frac{1}{t_i^2}]) X_i=[(f_i-\frac{1}{t_i^2})],\ Q([\frac{1}{t_i^2}]) X_j=[(f_i-\frac{\a_{ij}}{t_j})\frac{1}{t_j}],$$
$$A_iQ([\frac{1}{t_i^2}])=\be\cdot(f_i-\frac{1}{t_i^2})\unit_i+\bu[f_i\unit_i],\
A_jQ([\frac{1}{t_i^2}])=\be\cdot(f_i-\frac{\a_{ij}}{t_j})\unit_j+\bu[f_i\unit_j],$$
$$Q([\frac{1}{t_i^2}]) B_i=(f_i-\frac{1}{t_i^2})\frac{\unit_i}{t_i}\bu^*+[(f_i-\frac{1}{t_i^2})\frac{\unit_i}{t_i}]\be^*, \
Q([\frac{1}{t_i^2}]) B_j=(f_i-\frac{\a_{ij}}{t_j})\frac{\unit_j}{t_j}\bu^*+[(f_i-\frac{\a_{ij}}{t_j})\frac{\unit_j}{t_j}]\be^*, 
$$
$$ B_iQ(\be[\frac{1}{t_i^2}])=[\frac{1}{t_i^3}], \  Y_iQ(\be[\frac{1}{t_i^2}])=\be[\frac{1}{t_i^3}], \
 B_iQ(\be[\frac{1}{t_j^2}])= Y_iQ(\be[\frac{1}{t_j^2}])=0.$$
In the formulas below we omit commas between the arguments of $m_i$ for brevity. 
The indices denoted by different letters are assumed to be distinct.
Let $\wt{P_i}$ stand for either of the strings $A_i X_i$ or $Y_i A_i$.
Then the nonzero $m_4$ products of the basis elements are:
$$m_4( B_i Y_i Y_i A_i)=\sum_{j\neq_i}\b_{ij} X_j,$$
$$m_4( B_i Y_i A_i B_j)=-\ga_{ij} B_j,$$
$$m_4( B_i Y_i A_i X_i)=\sum_{j\neq_i}\b_{ij} X_j,$$
$$m_4( B_i Y_i A_i X_j)=-\ga_{ij} X_j,$$
$$m_4( B_i A_i X_i B_j)=-\ga_{ij} B_j,$$
$$m_4( B_i A_i X_i X_j)=-\ga_{ij} X_j,$$
$$m_4(A_j B_i Y_i A_i)=-\ga_{ij}A_j,$$
$$m_4(A_j B_i A_i X_i)=-\ga_{ij}A_j,$$
$$m_4( X_i B_i Y_i A_i)=\sum_{j\neq_i}\b_{ij} X_j,$$
$$m_4( X_j B_i Y_i A_i)=-\ga_{ij} X_j+\sum_{k\neq j}\a_{ij}\a_{jk} X_k,$$
$$m_4( X_i B_i A_i X_i)=\sum_{j\neq_i}\b_{ij} X_j,$$
$$m_4( X_j B_i A_i X_i)=-\ga_{ij} X_j+\sum_{k\neq j}\a_{ij}\a_{jk} X_k,$$

Let $Q_i$ stand for either $ B_i Y_i$ or $ X_i B_i$. Then the nonzero values of $m_5$ are:
$$m_5( B_i Y_i Y_i Y_i A_i)=\sum_{j\neq i}\eta_{ij} X_j,$$
$$m_5( B_i Y_i Y_i A_i X_i)=\sum_{j\neq i}\eta_{ij} X_j,$$
$$m_5( B_i\wt{P}_iQ_i)=\de_{ii} B_i,$$
$$m_5( B_i\wt{P}_iQ_j)=\de_{ij} B_j,$$
$$m_5( B_i\wt{P}_i X_i X_i)=\de_{ii}  X_i,$$
$$m_5( B_i\wt{P}_i X_j X_j)=\de_{ij} X_j,$$
$$m_5(Q_i A_i X_i X_i)=\de_{ii}  X_i,$$
$$m_5(Q_i\wt{P}_i B_j)=-\eps_{ij} B_j,$$
$$m_5(Q_i\wt{P}_i X_j)=-\eps_{ij} X_j,$$
$$m_5(\wt{P}_i B_i\wt{P}_i)=m_5( A_i B_i\wt{P}_i X_i)=-\de_{ii} A_i,$$
$$m_5( A_i B_i\wt{P}_i B_i)=-\de_{ii}e_{\OO_{p_i}},$$
$$m_5( A_jQ_i\wt{P}_i)=\eps_{ij} A_j,$$ 
$$m_5( A_j X_j B_i\wt{P}_i)=m_5( A_j B_i\wt{P}_i X_j)=m_5( Y_j A_j B_i\wt{P}_i)=-\de_{ij} A_j,$$
$$m_5( A_k X_j B_i\wt{P}_i)=\a_{ij}\ga_{jk} A_k,$$
$$m_5( A_i X_j B_i\wt{P}_i)=\a_{ij}\ga_{ji} A_i,$$
$$m_5( A_j B_i\wt{P}_i B_j)=-\de_{ij}e_{\OO_{p_i}},$$
$$m_5( X_i X_i B_i\wt{P}_i)=m_5( X_i B_i Y_i Y_i A_i)=\sum_{j\neq i}\eta_{ij} X_j,$$
$$m_5( X_j B_i\wt{P}_i B_j)=\de_{ij} B_j,$$
$$m_5( X_i B_i Y_i A_i X_i)=\sum_{j\neq i}\eta_{ij} X_j+\de_{ii} X_i,$$
$$m_5( X_i B_i\wt{P}_i B_i)=\de_{ii} B_i,$$ 
$$m_5( X_j B_i\wt{P}_i B_i)=-\a_{ij}\ga_{ji} B_i,$$
$$m_5( X_j B_i\wt{P}_i B_k)=-\a_{ij}\ga_{jk} B_k,$$
$$m_5( X_j B_i\wt{P}_i X_j)=\de_{ij} X_j,$$
$$m_5( X_k B_i\wt{P}_i X_j)=-\a_{ik}\ga_{kj} X_j,$$
$$m_5( X_j B_i Y_i A_i X_i)=\sum_{k\neq j}\b_{ij}\a_{jk} X_k-\eps_{ij} X_j-\a_{ij}\ga_{ji} X_i,$$
$$m_5( X_j B_i A_i X_i X_i)=-\a_{ij}\ga_{ji} X_i,$$
$$m_5( X_j X_j B_i\wt{P}_i)=\sum_{k\neq j}\a_{ij}\b_{jk} X_k,$$
$$m_5( X_j B_i Y_i Y_i A_i)=m_5( X_j X_i B_i\wt{P}_i)=\sum_{k\neq j}\b_{ij}\a_{jk} X_k-\eps_{ij} X_j,$$
$$m_5( X_i X_j B_i\wt{P}_i)=\sum_{l\neq i}\a_{ij}\a_{ji}\a_{il} X_l-\a_{ij}\ga_{ji} X_i,$$
$$m_5( X_k X_j B_i\wt{P}_i)=\sum_{l\neq k}\a_{ij}\a_{jk}\a_{kl} X_l-\a_{ij}\ga_{jk} X_k.$$

\end{document}